\documentclass[11pt,leqno]{article}
\usepackage{amsmath,amssymb,amsthm}
\usepackage[mathscr]{eucal}
\begin{document}

%
%
%
\ifx\figfortexisloaded\relax \else\let\figfortexisloaded=\relax\fi
\message{version 1.6.1}
\newif\iftextures
\catcode`\@=11
\newdimen\epsil@n\epsil@n=0.00005pt
\newdimen\Cepsil@n\Cepsil@n=0.005pt
\newdimen\dcq@\dcq@=254pt
\newdimen\PI@\PI@=3.141592pt
\newdimen\PI@deg\PI@deg=360pt
\newdimen\PI@sDdeg\PI@sDdeg=180pt
\chardef\t@n=10
\newcount\p@rtent \newcount\f@ctech \newcount\result@tent
\newdimen\v@lmin \newdimen\v@lmax \newdimen\v@leur
\newdimen\result@t\newdimen\result@@t
\newdimen\mili@u \newdimen\c@rre \newdimen\delt@
\def\degT@rd{0.017453} 
\def\rdT@deg{57.295779} 
{\catcode`p=12 \catcode`t=12 \gdef\v@leurseule#1pt{#1}}
\def\repdecn@mb#1{\expandafter\v@leurseule\the#1\space}
\def\arct@n#1(#2,#3){{\v@lmin=#2\v@lmax=#3%
    \maxim@m{\mili@u}{-\v@lmin}{\v@lmin}\maxim@m{\c@rre}{-\v@lmax}{\v@lmax}%
    \delt@=\mili@u\m@ech\mili@u%
    \ifdim\c@rre>9949\mili@u\v@lmax=0.5\v@lmax\c@lATAN\v@leur(\z@,\v@lmax)
    \else%
    \maxim@m{\mili@u}{-\v@lmin}{\v@lmin}\maxim@m{\c@rre}{-\v@lmax}{\v@lmax}%
    \m@ech\c@rre%
    \ifdim\mili@u>9949\c@rre\v@lmin=0.5\v@lmin
    \maxim@m{\mili@u}{-\v@lmin}{\v@lmin}\c@lATAN\v@leur(\mili@u,\z@)%
    \else\c@lATAN\v@leur(\delt@,\v@lmax)\fi\fi%
    \ifdim\v@lmin<\z@\v@leur=-\v@leur\ifdim\v@lmax<\z@\advance\v@leur-\PI@%
    \else\advance\v@leur\PI@\fi\fi%
    \global\result@t=\v@leur}#1=\result@t}
\def\m@ech#1{\ifdim#1>1.646pt\divide\mili@u\t@n\divide\c@rre\t@n\m@ech#1\fi}
\def\c@lATAN#1(#2,#3){{\v@lmin=#2\v@lmax=#3\v@leur=\z@\delt@=\tw@ pt%
    \un@iter{0.785398}{\v@lmax<}%
    \un@iter{0.463648}{\v@lmax<}%
    \un@iter{0.244979}{\v@lmax<}%
    \un@iter{0.124355}{\v@lmax<}%
    \un@iter{0.062419}{\v@lmax<}%
    \un@iter{0.031240}{\v@lmax<}%
    \un@iter{0.015624}{\v@lmax<}%
    \un@iter{0.007812}{\v@lmax<}%
    \un@iter{0.003906}{\v@lmax<}%
    \un@iter{0.001953}{\v@lmax<}%
    \un@iter{0.000976}{\v@lmax<}%
    \un@iter{0.000488}{\v@lmax<}%
    \un@iter{0.000244}{\v@lmax<}%
    \un@iter{0.000122}{\v@lmax<}%
    \un@iter{0.000061}{\v@lmax<}%
    \un@iter{0.000030}{\v@lmax<}%
    \un@iter{0.000015}{\v@lmax<}%
    \global\result@t=\v@leur}#1=\result@t}
\def\un@iter#1#2{%
    \divide\delt@\tw@\edef\dpmn@{\repdecn@mb{\delt@}}%
    \mili@u=\v@lmin%
    \ifdim#2\z@%
      \advance\v@lmin-\dpmn@\v@lmax\advance\v@lmax\dpmn@\mili@u%
      \advance\v@leur-#1pt%
    \else%
      \advance\v@lmin\dpmn@\v@lmax\advance\v@lmax-\dpmn@\mili@u%
      \advance\v@leur#1pt%
    \fi}
\def\c@ssin#1#2#3{{\v@leur=#3\maxim@m{\mili@u}{-\v@leur}{\v@leur}%
    \mili@u=\repdecn@mb\PI@\mili@u\divide\mili@u180 %
    \loop\ifdim\mili@u>2\PI@\advance\mili@u-2\PI@\repeat%
    \ifdim\mili@u<0.5\PI@\c@lCS{\v@lmin}{\v@lmax}{\mili@u}%
    \else\ifdim\mili@u<\PI@\advance\mili@u-0.5\PI@\c@lCS{\v@lmin}{\v@lmax}{\mili@u}%
         \mili@u=\v@lmin\v@lmin=-\v@lmax\v@lmax=\mili@u%
    \else\ifdim\mili@u<1.5\PI@\advance\mili@u-\PI@\c@lCS{\v@lmin}{\v@lmax}{\mili@u}%
         \v@lmin=-\v@lmin\v@lmax=-\v@lmax%
    \else\advance\mili@u-2\PI@\c@lCS{\v@lmin}{\v@lmax}{\mili@u}\fi\fi\fi%
    \ifdim\v@leur<\z@\v@lmax=-\v@lmax\fi%
    \global\result@t=\v@lmin\global\result@@t=\v@lmax}#1=\result@t#2=\result@@t}
\def\c@lCS#1#2#3{{\v@lmin=0.607253pt\v@lmax=\z@\v@leur=#3\delt@=\tw@ pt%
    \un@iter{0.785398}{\v@leur>}%
    \un@iter{0.463648}{\v@leur>}%
    \un@iter{0.244979}{\v@leur>}%
    \un@iter{0.124355}{\v@leur>}%
    \un@iter{0.062419}{\v@leur>}%
    \un@iter{0.031240}{\v@leur>}%
    \un@iter{0.015624}{\v@leur>}%
    \un@iter{0.007812}{\v@leur>}%
    \un@iter{0.003906}{\v@leur>}%
    \un@iter{0.001953}{\v@leur>}%
    \un@iter{0.000976}{\v@leur>}%
    \un@iter{0.000488}{\v@leur>}%
    \un@iter{0.000244}{\v@leur>}%
    \un@iter{0.000122}{\v@leur>}%
    \un@iter{0.000061}{\v@leur>}%
    \un@iter{0.000030}{\v@leur>}%
    \global\result@t=\v@lmin\global\result@@t=\v@lmax}#1=\result@t#2=\result@@t}
\def\invers@#1#2{{\v@leur=#2%
    \maxim@m{\v@lmax}{-\v@leur}{\v@leur}
    \f@ctech=\@ne\loop\ifdim\v@lmax<\p@\multiply\v@lmax\t@n\multiply\f@ctech\t@n\repeat%
    \multiply\v@leur\f@ctech\edef\v@lv@leur{\repdecn@mb{\v@leur}}%
    \p@rtentiere{\p@rtent}{\v@leur}\v@lmin=\p@\divide\v@lmin\p@rtent%
    \loop\v@lmax=-\v@lmin\v@lmax=\v@lv@leur\v@lmax%
    \advance\v@lmax\tw@ pt\v@lmax=\repdecn@mb{\v@lmin}\v@lmax%
    \delt@=\v@lmax\advance\delt@-\v@lmin\ifdim\delt@<\z@\delt@=-\delt@\fi%
    \ifdim\delt@>\epsil@n\v@lmin=\v@lmax\repeat\multiply\v@lmax\f@ctech%
    \global\result@t=\v@lmax}#1=\result@t}
\def\minim@m#1#2#3{\relax\ifdim#2<#3#1=#2\else#1=#3\fi}
\def\maxim@m#1#2#3{\relax\ifdim#2>#3#1=#2\else#1=#3\fi}
\def\p@rtentiere#1#2{#1=#2\divide#1by65536 }
\def\r@undint#1#2{{\v@leur=#2\v@leur=0.1\v@leur\p@rtentiere{\p@rtent}{\v@leur}%
    \v@leur=\p@rtent pt\global\result@t=10\v@leur}#1=\result@t}
\def\sqrt@#1#2{{\v@leur=#2%
    \minim@m{\v@lmin}{\p@}{\v@leur}\maxim@m{\v@lmax}{\p@}{\v@leur}%
    \f@ctech=\@ne\loop\ifdim\v@leur>\dcq@\divide\v@leur100 \v@lmax=\v@leur%
                   \multiply\f@ctech\t@n\repeat%
    \loop\mili@u=\v@lmin\advance\mili@u\v@lmax\divide\mili@u\tw@%
    \c@rre=\repdecn@mb{\mili@u}\mili@u%
    \ifdim\c@rre<\v@leur\v@lmin=\mili@u\else\v@lmax=\mili@u\fi%
    \delt@=\v@lmax\advance\delt@-\v@lmin%
    \ifdim\delt@>\epsil@n\repeat%
    \mili@u=\v@lmin\advance\mili@u\v@lmax\divide\mili@u\tw@\multiply\mili@u\f@ctech%
    \global\result@t=\mili@u}#1=\result@t}
\def\extrairelepremi@r#1\de#2{\expandafter\lepremi@r#2@#1#2}
\def\lepremi@r#1,#2@#3#4{\def#3{#1}\def#4{#2}\ignorespaces}
\def\@cfor#1:=#2\do#3{%
  \edef\@fortemp{#2}%
  \ifx\@fortemp\@empty\else\@cforloop#2,\@nil,\@nil\@@#1{#3}\fi}
\def\@cforloop#1,#2\@@#3#4{%
  \def#3{#1}%
  \ifx#3\@nnil\let\@nextwhile=\@fornoop\else#4\relax\let\@nextwhile=\@cforloop\fi%
  \@nextwhile#2\@@#3{#4}}

\def\@ecfor#1:=#2\do#3{%
  \def\@@cfor{\@cfor#1:=}%
  \edef\@@@cfor{#2}%
  \expandafter\@@cfor\@@@cfor\do{#3}}
\def\@empty{}
\def\@nnil{\@nil}
\def\@fornoop#1\@@#2#3{}
\newbox\b@xvisu
\newtoks\let@xte
\newif\ifitis@K
\newcount\s@mme
\newcount\l@mbd@un \newcount\l@mbd@de
\newcount\superc@ntr@l\superc@ntr@l=\@ne        
\newcount\typec@ntr@l\typec@ntr@l=\superc@ntr@l 
\newdimen\v@lX  \newdimen\v@lY  \newdimen\v@lZ
\newdimen\v@lXa \newdimen\v@lYa \newdimen\v@lZa
\newdimen\unit@\unit@=\p@ 
\def\unit@util{pt}
\def\ptT@ptps{0.99626400996264009962}
\def\ptpsT@pt{1.00375}
\def\ptT@unit@{1} 
\def\setunit@#1{\def\unit@util{#1}\setunit@@#1:\invers@{\result@t}{\unit@}%
    \edef\ptT@unit@{\repdecn@mb\result@t}}
\def\setunit@@#1#2:{\ifcat#1a\unit@=\@ne#1#2\else\unit@=#1#2\fi}
\newif\ifBdingB@x\BdingB@xtrue
\newdimen\c@@rdXmin \newdimen\c@@rdYmin  
\newdimen\c@@rdXmax \newdimen\c@@rdYmax
\def\b@undb@x#1#2{\ifBdingB@x%
    \minim@m{\global\c@@rdXmin}{\c@@rdXmin}{#1}%
    \minim@m{\global\c@@rdYmin}{\c@@rdYmin}{#2}%
    \maxim@m{\global\c@@rdXmax}{\c@@rdXmax}{#1}%
    \maxim@m{\global\c@@rdYmax}{\c@@rdYmax}{#2}\fi}
\def\b@undb@xP#1{{\figg@tXY{#1}\b@undb@x{\v@lX}{\v@lY}}}
\def\ellBB@x#1;#2,#3(#4,#5,#6){{\s@uvc@ntr@l\et@tellBB@x%
    \setc@ntr@l{2}\figptell-2::#1;#2,#3(#4,#6)\b@undb@xP{-2}%
    \figptell-2::#1;#2,#3(#5,#6)\b@undb@xP{-2}%
    \c@ssin{\v@lmin}{\v@lmax}{#6pt}%
    \mili@u=#3\v@lmin\delt@=#2\v@lmax\arct@n\v@leur(\delt@,\mili@u)%
    \mili@u=-#3\v@lmax\delt@=#2\v@lmin\arct@n\c@rre(\delt@,\mili@u)%
    \v@leur=\rdT@deg\v@leur\advance\v@leur-\PI@deg%
    \c@rre=\rdT@deg\c@rre\advance\c@rre-\PI@deg%
    \v@lmin=#4pt\v@lmax=#5pt%
    \loop\ifdim\v@leur<\v@lmax\ifdim\v@leur>\v@lmin%
    \edef\@ngle{\repdecn@mb\v@leur}\figptell-2::#1;#2,#3(\@ngle,#6)%
    \b@undb@xP{-2}\fi\advance\v@leur\PI@sDdeg\repeat%
    \loop\ifdim\c@rre<\v@lmax\ifdim\c@rre>\v@lmin%
    \edef\@ngle{\repdecn@mb\c@rre}\figptell-2::#1;#2,#3(\@ngle,#6)%
    \b@undb@xP{-2}\fi\advance\c@rre\PI@sDdeg\repeat%
    \resetc@ntr@l\et@tellBB@x}\ignorespaces}
\def\initb@undb@x{\c@@rdXmin=\maxdimen\c@@rdYmin=\maxdimen%
    \c@@rdXmax=-\maxdimen\c@@rdYmax=-\maxdimen}
\def\c@ntr@lnum#1{%
    \relax\ifnum\typec@ntr@l=\@ne%
    \ifnum#1<\z@%
    \immediate\write16{*** Forbidden point number (#1). Abort.}\end\fi\fi%
    \set@bjc@de{#1}}
\def\set@bjc@de#1{\edef\objc@de{@BJ\ifnum#1<\z@ M\romannumeral-#1\else\romannumeral#1\fi}}
\def\setc@ntr@l#1{\ifnum\superc@ntr@l>#1\typec@ntr@l=\superc@ntr@l%
    \else\typec@ntr@l=#1\fi}
\def\resetc@ntr@l#1{\global\superc@ntr@l=#1\setc@ntr@l{#1}}
\def\s@uvc@ntr@l#1{\edef#1{\the\superc@ntr@l}}
\def\c@lproscalDD#1[#2,#3]{{\figg@tXY{#2}%
    \edef\Xu@{\repdecn@mb{\v@lX}}\edef\Yu@{\repdecn@mb{\v@lY}}\figg@tXY{#3}%
    \global\result@t=\Xu@\v@lX\global\advance\result@t\Yu@\v@lY}#1=\result@t}
\def\c@lproscalTD#1[#2,#3]{{\figg@tXY{#2}\edef\Xu@{\repdecn@mb{\v@lX}}%
    \edef\Yu@{\repdecn@mb{\v@lY}}\edef\Zu@{\repdecn@mb{\v@lZ}}%
    \figg@tXY{#3}\global\result@t=\Xu@\v@lX\global\advance\result@t\Yu@\v@lY%
    \global\advance\result@t\Zu@\v@lZ}#1=\result@t}
\def\c@lprovec#1{%
    \det@rmC\v@lZa(\v@lX,\v@lY,\v@lmin,\v@lmax)%
    \det@rmC\v@lXa(\v@lY,\v@lZ,\v@lmax,\v@leur)%
    \det@rmC\v@lYa(\v@lZ,\v@lX,\v@leur,\v@lmin)%
    \figv@ctCreg#1(\v@lXa,\v@lYa,\v@lZa)}
\def\det@rm#1[#2,#3]{{\figg@tXY{#2}\figg@tXYa{#3}%
    \delt@=\repdecn@mb{\v@lX}\v@lYa\advance\delt@-\repdecn@mb{\v@lY}\v@lXa%
    \global\result@t=\delt@}#1=\result@t}
\def\det@rmC#1(#2,#3,#4,#5){{\global\result@t=\repdecn@mb{#2}#5%
    \global\advance\result@t-\repdecn@mb{#3}#4}#1=\result@t}
\def\getredf@ctDD#1(#2,#3){{\maxim@m{\v@lXa}{-#2}{#2}\maxim@m{\v@lYa}{-#3}{#3}%
    \maxim@m{\v@lXa}{\v@lXa}{\v@lYa}
    \ifdim\v@lXa>64pt\divide\v@lXa64 %
    \p@rtentiere{\p@rtent}{\v@lXa}\advance\p@rtent\@ne\else\p@rtent=\@ne\fi%
    \global\result@tent=\p@rtent}#1=\result@tent\ignorespaces}
\def\getredf@ctTD#1(#2,#3,#4){{\maxim@m{\v@lXa}{-#2}{#2}\maxim@m{\v@lYa}{-#3}{#3}%
    \maxim@m{\v@lZa}{-#4}{#4}\maxim@m{\v@lXa}{\v@lXa}{\v@lYa}%
    \maxim@m{\v@lXa}{\v@lXa}{\v@lZa}
    \ifdim\v@lXa>42pt\divide\v@lXa42 %
    \p@rtentiere{\p@rtent}{\v@lXa}\advance\p@rtent\@ne\else\p@rtent=\@ne\fi%
    \global\result@tent=\p@rtent}#1=\result@tent\ignorespaces}
\def\figptintercircB@zDD#1:#2[#3,#4,#5,#6]{{\s@uvc@ntr@l\et@tfigptintercircB@zDD%
    \setc@ntr@l{2}\figvectPDD-1[#3,#6]\figg@tXY{-1}\getredf@ctDD\f@ctech(\v@lX,\v@lY)%
    \v@lmin=\z@\v@lmax=\p@\mili@u=#2\unit@\divide\mili@u\f@ctech%
    \c@rre=\repdecn@mb{\mili@u}\mili@u\loop\mili@u=\v@lmin\advance\mili@u\v@lmax%
    \divide\mili@u\tw@\edef\T@{\repdecn@mb{\mili@u}}\figptBezierDD-2::\T@[#3,#4,#5,#6]%
    \figvectPDD-1[#3,-2]\n@rmeucCDD{\delt@}{-1}\ifdim\delt@>\c@rre\v@lmax=\mili@u%
    \else\v@lmin=\mili@u\fi\v@leur=\v@lmax\advance\v@leur-\v@lmin%
    \ifdim\v@leur>\epsil@n\repeat\figptcopyDD#1:/-2/%
    \resetc@ntr@l\et@tfigptintercircB@zDD}\ignorespaces}
\def\inters@cDD#1:#2[#3,#4;#5,#6]{{\s@uvc@ntr@l\et@tinters@cDD%
    \setc@ntr@l{2}\vecunit@{-1}{#4}\vecunit@{-2}{#6}%
    \figg@tXY{-1}\setc@ntr@l{1}\figg@tXYa{#3}%
    \edef\A@{\repdecn@mb{\v@lX}}\edef\B@{\repdecn@mb{\v@lY}}%
    \v@lmin=\B@\v@lXa\advance\v@lmin-\A@\v@lYa%
    \figg@tXYa{#5}\setc@ntr@l{2}\figg@tXY{-2}%
    \edef\C@{\repdecn@mb{\v@lX}}\edef\D@{\repdecn@mb{\v@lY}}%
    \v@lmax=\D@\v@lXa\advance\v@lmax-\C@\v@lYa%
    \delt@=\A@\v@lY\advance\delt@-\B@\v@lX%
    \invers@{\v@leur}{\delt@}\edef\v@ldelta{\repdecn@mb{\v@leur}}%
    \v@lXa=\A@\v@lmax\advance\v@lXa-\C@\v@lmin%
    \v@lYa=\B@\v@lmax\advance\v@lYa-\D@\v@lmin%
    \v@lXa=\v@ldelta\v@lXa\v@lYa=\v@ldelta\v@lYa%
    \setc@ntr@l{1}\figp@intregDD#1:{#2}(\v@lXa,\v@lYa)%
    \resetc@ntr@l\et@tinters@cDD}\ignorespaces}
\def\inters@cTD#1:#2[#3,#4;#5,#6]{{\s@uvc@ntr@l\et@tinters@cTD%
    \setc@ntr@l{2}\figvectNVTD-1[#4,#6]\figvectNVTD-2[#6,-1]\figvectPTD-1[#3,#5]%
    \r@pPSTD\v@leur[-2,-1,#4]\edef\v@lcoef{\repdecn@mb{\v@leur}}%
    \figpttraTD#1:{#2}=#3/\v@lcoef,#4/\resetc@ntr@l\et@tinters@cTD}\ignorespaces}
\def\r@pPSTD#1[#2,#3,#4]{{\figg@tXY{#2}\edef\Xu@{\repdecn@mb{\v@lX}}%
    \edef\Yu@{\repdecn@mb{\v@lY}}\edef\Zu@{\repdecn@mb{\v@lZ}}%
    \figg@tXY{#3}\v@lmin=\Xu@\v@lX\advance\v@lmin\Yu@\v@lY\advance\v@lmin\Zu@\v@lZ%
    \figg@tXY{#4}\v@lmax=\Xu@\v@lX\advance\v@lmax\Yu@\v@lY\advance\v@lmax\Zu@\v@lZ%
    \invers@{\v@leur}{\v@lmax}\global\result@t=\repdecn@mb{\v@leur}\v@lmin}%
    #1=\result@t}
\def\n@rminfDD#1#2{{\figg@tXY{#2}\maxim@m{\v@lX}{\v@lX}{-\v@lX}%
    \maxim@m{\v@lY}{\v@lY}{-\v@lY}\maxim@m{\global\result@t}{\v@lX}{\v@lY}}%
    #1=\result@t}
\def\n@rminfTD#1#2{{\figg@tXY{#2}\maxim@m{\v@lX}{\v@lX}{-\v@lX}%
    \maxim@m{\v@lY}{\v@lY}{-\v@lY}\maxim@m{\v@lZ}{\v@lZ}{-\v@lZ}%
    \maxim@m{\v@lX}{\v@lX}{\v@lY}\maxim@m{\global\result@t}{\v@lX}{\v@lZ}}%
    #1=\result@t}
\def\n@rmeucCDD#1#2{\figg@tXY{#2}\divide\v@lX\f@ctech\divide\v@lY\f@ctech%
    #1=\repdecn@mb{\v@lX}\v@lX\v@lX=\repdecn@mb{\v@lY}\v@lY\advance#1\v@lX}
\def\n@rmeucCTD#1#2{\figg@tXY{#2}%
    \divide\v@lX\f@ctech\divide\v@lY\f@ctech\divide\v@lZ\f@ctech%
    #1=\repdecn@mb{\v@lX}\v@lX\v@lX=\repdecn@mb{\v@lY}\v@lY\advance#1\v@lX%
    \v@lX=\repdecn@mb{\v@lZ}\v@lZ\advance#1\v@lX}
\def\n@rmeucSVDD#1#2{{\figg@tXY{#2}%
    \v@lXa=\repdecn@mb{\v@lX}\v@lX\v@lYa=\repdecn@mb{\v@lY}\v@lY%
    \advance\v@lXa\v@lYa\sqrt@{\global\result@t}{\v@lXa}}#1=\result@t}
\def\n@rmeucSVTD#1#2{{\figg@tXY{#2}\v@lXa=\repdecn@mb{\v@lX}\v@lX%
    \v@lYa=\repdecn@mb{\v@lY}\v@lY\v@lZa=\repdecn@mb{\v@lZ}\v@lZ%
    \advance\v@lXa\v@lYa\advance\v@lXa\v@lZa\sqrt@{\global\result@t}{\v@lXa}}#1=\result@t}
\def\n@rmeucDD#1#2{{\figg@tXY{#2}\getredf@ctDD\f@ctech(\v@lX,\v@lY)%
    \divide\v@lX\f@ctech\divide\v@lY\f@ctech%
    \v@lXa=\repdecn@mb{\v@lX}\v@lX\v@lYa=\repdecn@mb{\v@lY}\v@lY%
    \advance\v@lXa\v@lYa\sqrt@{\global\result@t}{\v@lXa}%
    \global\multiply\result@t\f@ctech}#1=\result@t}
\def\n@rmeucTD#1#2{{\figg@tXY{#2}\getredf@ctTD\f@ctech(\v@lX,\v@lY,\v@lZ)%
    \divide\v@lX\f@ctech\divide\v@lY\f@ctech\divide\v@lZ\f@ctech%
    \v@lXa=\repdecn@mb{\v@lX}\v@lX%
    \v@lYa=\repdecn@mb{\v@lY}\v@lY\v@lZa=\repdecn@mb{\v@lZ}\v@lZ%
    \advance\v@lXa\v@lYa\advance\v@lXa\v@lZa\sqrt@{\global\result@t}{\v@lXa}%
    \global\multiply\result@t\f@ctech}#1=\result@t}
\def\vecunit@DD#1#2{{\figg@tXY{#2}\getredf@ctDD\f@ctech(\v@lX,\v@lY)%
    \divide\v@lX\f@ctech\divide\v@lY\f@ctech%
    \figv@ctCreg#1(\v@lX,\v@lY)\n@rmeucSV{\v@lYa}{#1}%
    \invers@{\v@lXa}{\v@lYa}\edef\v@lv@lXa{\repdecn@mb{\v@lXa}}%
    \v@lX=\v@lv@lXa\v@lX\v@lY=\v@lv@lXa\v@lY%
    \figv@ctCreg#1(\v@lX,\v@lY)\multiply\v@lYa\f@ctech\global\result@t=\v@lYa}}
\def\vecunit@TD#1#2{{\figg@tXY{#2}\getredf@ctTD\f@ctech(\v@lX,\v@lY,\v@lZ)%
    \divide\v@lX\f@ctech\divide\v@lY\f@ctech\divide\v@lZ\f@ctech%
    \figv@ctCreg#1(\v@lX,\v@lY,\v@lZ)\n@rmeucSV{\v@lYa}{#1}%
    \invers@{\v@lXa}{\v@lYa}\edef\v@lv@lXa{\repdecn@mb{\v@lXa}}%
    \v@lX=\v@lv@lXa\v@lX\v@lY=\v@lv@lXa\v@lY\v@lZ=\v@lv@lXa\v@lZ%
    \figv@ctCreg#1(\v@lX,\v@lY,\v@lZ)\multiply\v@lYa\f@ctech\global\result@t=\v@lYa}}
\def\vecunitC@TD[#1,#2]{\figg@tXYa{#1}\figg@tXY{#2}%
    \advance\v@lX-\v@lXa\advance\v@lY-\v@lYa\advance\v@lZ-\v@lZa\c@lvecunitTD}
\def\vecunitCV@TD#1{\figg@tXY{#1}\c@lvecunitTD}
\def\c@lvecunitTD{\getredf@ctTD\f@ctech(\v@lX,\v@lY,\v@lZ)%
    \divide\v@lX\f@ctech\divide\v@lY\f@ctech\divide\v@lZ\f@ctech%
    \v@lXa=\repdecn@mb{\v@lX}\v@lX%
    \v@lYa=\repdecn@mb{\v@lY}\v@lY\v@lZa=\repdecn@mb{\v@lZ}\v@lZ%
    \advance\v@lXa\v@lYa\advance\v@lXa\v@lZa\sqrt@{\v@lYa}{\v@lXa}%
    \invers@{\v@lXa}{\v@lYa}\edef\v@lv@lXa{\repdecn@mb{\v@lXa}}%
    \v@lX=\v@lv@lXa\v@lX\v@lY=\v@lv@lXa\v@lY\v@lZ=\v@lv@lXa\v@lZ}
\def\figgetangleDD#1[#2,#3,#4]{\ifps@cri{\s@uvc@ntr@l\et@tfiggetangleDD\setc@ntr@l{2}%
    \figvectPDD-1[#2,#3]\figvectPDD-2[#2,#4]\vecunit@{-1}{-1}%
    \c@lproscalDD\delt@[-2,-1]\figvectNVDD-1[-1]\c@lproscalDD\v@leur[-2,-1]%
    \arct@n\v@lmax(\delt@,\v@leur)\v@lmax=\rdT@deg\v@lmax%
    \ifdim\v@lmax<\z@\advance\v@lmax\PI@deg\fi\xdef#1{\repdecn@mb{\v@lmax}}%
    \resetc@ntr@l\et@tfiggetangleDD}\ignorespaces\fi}
\def\figgetangleTD#1[#2,#3,#4,#5]{\ifps@cri{\s@uvc@ntr@l\et@tfiggetangleTD\setc@ntr@l{2}%
    \figvectPTD-1[#2,#3]\figvectPTD-2[#2,#5]\figvectNVTD-3[-1,-2]%
    \figvectPTD-2[#2,#4]\figvectNVTD-4[-3,-1]%
    \vecunit@{-1}{-1}\c@lproscalTD\delt@[-2,-1]\c@lproscalTD\v@leur[-2,-4]%
    \arct@n\v@lmax(\delt@,\v@leur)\v@lmax=\rdT@deg\v@lmax%
    \ifdim\v@lmax<\z@\advance\v@lmax\PI@deg\fi\xdef#1{\repdecn@mb{\v@lmax}}%
    \resetc@ntr@l\et@tfiggetangleTD}\ignorespaces\fi}    
\def\figgetdist#1[#2,#3]{\ifps@cri{\s@uvc@ntr@l\et@tfiggetdist\setc@ntr@l{2}%
    \figvectP-1[#2,#3]\n@rmeuc{\v@lX}{-1}\v@lX=\ptT@unit@\v@lX\xdef#1{\repdecn@mb{\v@lX}}%
    \resetc@ntr@l\et@tfiggetdist}\ignorespaces\fi}
\def\figg@tT#1{\c@ntr@lnum{#1}%
    {\expandafter\expandafter\expandafter\extr@ctT\csname\objc@de\endcsname:%
     \ifdim\wd\Gb@x=\z@\ptn@me{#1}\else\unhcopy\Gb@x\fi}}
\def\extr@ctT#1)#2:{\setbox\Gb@x=\hbox{#2}}
\def\figg@tXY#1{\c@ntr@lnum{#1}%
    \expandafter\expandafter\expandafter\extr@ctC\csname\objc@de\endcsname:}
\def\extr@ctCDD#1(#2,#3,#4:{\v@lX=#2\v@lY=#3}
\def\extr@ctCTD#1(#2,#3,#4)#5:{\v@lX=#2\v@lY=#3\v@lZ=#4}
\def\figg@tXYa#1{\c@ntr@lnum{#1}%
    \expandafter\expandafter\expandafter\extr@ctCa\csname\objc@de\endcsname:}
\def\extr@ctCaDD#1(#2,#3,#4:{\v@lXa=#2\v@lYa=#3}
\def\extr@ctCaTD#1(#2,#3,#4)#5:{\v@lXa=#2\v@lYa=#3\v@lZa=#4}
\def\figinit#1{\initpr@lim\figinit@#1,:\initpss@ttings\ignorespaces}
\def\figinit@#1,#2:{\setunit@{#1}\def\t@xt@{#2}\ifx\t@xt@\empty\else\figinit@@#2:\fi}
\def\figinit@@#1#2:{\if#12 \else\figs@tproj{#1}\initTD@\fi}
\newif\ifTr@isDim
\def\UnD@fined{UNDEFINED}
\def\ifundefined#1{\expandafter\ifx\csname#1\endcsname\relax}
\def\initpr@lim{\initb@undb@x\figsetmark{}\figsetptname{$A_{##1}$}\def\Sc@leFact{1}%
    \initDD@\figsetroundcoord{yes}\ps@critrue%
    \edef\disob@unit{\UnD@fined}}
\def\initDD@{\Tr@isDimfalse%
    \let\c@lDCUn=\c@lDCUnDD%
    \let\c@lDCDeux=\c@lDCDeuxDD%
    \let\c@ldefproj=\relax%
    \let\c@lproscal=\c@lproscalDD%
    \let\c@lprojSP=\relax%
    \let\extr@ctC=\extr@ctCDD%
    \let\extr@ctCa=\extr@ctCaDD%
    \let\extr@ctCF=\extr@ctCFDD%
    \let\figp@intreg=\figp@intregDD%
    \let\n@rmeucSV=\n@rmeucSVDD\let\n@rmeuc=\n@rmeucDD\let\n@rminf=\n@rminfDD%
    \let\vecunit@=\vecunit@DD%
    \let\figcoord=\figcoordDD%
    \let\figgetangle=\figgetangleDD%
    \let\figpt=\figptDD%
    \let\figptBezier=\figptBezierDD%
    \let\figptbary=\figptbaryDD%
    \let\figptbaryR=\figptbaryRDD%
    \let\figptcirc=\figptcircDD%
    \let\figptcircumcenter=\figptcircumcenterDD%
    \let\figptcopy=\figptcopyDD%
    \let\figptcurvcenter=\figptcurvcenterDD%
    \let\figptell=\figptellDD%
    \let\figptendnormal=\figptendnormalDD%
    \def\figptinterlineplane{\un@v@ilable{figptinterlineplane}}%
    \let\figptinterlines=\inters@cDD%
    \let\figptorthocenter=\figptorthocenterDD%
    \let\figptorthoprojline=\figptorthoprojlineDD%
    \def\figptorthoprojplane{\un@v@ilable{figptorthoprojplane}}%
    \let\figptrot=\figptrotDD%
    \let\figptscontrol=\figptscontrolDD%
    \let\figptsintercirc=\figptsintercircDD%
    \let\figptsorthoprojline=\figptsorthoprojlineDD%
    \let\figptsrot=\figptsrotDD%
    \let\figptssym=\figptssymDD%
    \let\figptstra=\figptstraDD%
    \let\figptsym=\figptsymDD%
    \let\figpttraC=\figpttraCDD%
    \let\figpttra=\figpttraDD%
    \def\figsetobdist{\un@v@ilable{figsetobdist}}%
    \def\figsettarget{\un@v@ilable{figsettarget}}%
    \def\figsetview{\un@v@ilable{figsetview}}%
    \let\figvectDBezier=\figvectDBezierDD%
    \let\figvectN=\figvectNDD%
    \let\figvectNV=\figvectNVDD%
    \let\figvectP=\figvectPDD%
    \let\figvectU=\figvectUDD%
    \let\psarccircP=\psarccircPDD%
    \let\psarccirc=\psarccircDD%
    \let\psarcell=\psarcellDD%
    \let\psarcellPA=\psarcellPADD%
    \let\psarrowBezier=\psarrowBezierDD%
    \let\psarrowcircP=\psarrowcircPDD%
    \let\psarrowcirc=\psarrowcircDD%
    \let\psarrowhead=\psarrowheadDD%
    \let\psarrow=\psarrowDD%
    \let\psBezier=\psBezierDD%
    \let\pscirc=\pscircDD%
    \let\pscurve=\pscurveDD%
    \let\psnormal=\psnormalDD%
    }
\def\initTD@{\Tr@isDimtrue\initb@undb@xTD\newt@rgetptfalse\newdis@bfalse%
    \let\c@lDCUn=\c@lDCUnTD%
    \let\c@lDCDeux=\c@lDCDeuxTD%
    \let\c@ldefproj=\c@ldefprojTD%
    \let\c@lproscal=\c@lproscalTD%
    \let\extr@ctC=\extr@ctCTD%
    \let\extr@ctCa=\extr@ctCaTD%
    \let\extr@ctCF=\extr@ctCFTD%
    \let\figp@intreg=\figp@intregTD%
    \let\n@rmeucSV=\n@rmeucSVTD\let\n@rmeuc=\n@rmeucTD\let\n@rminf=\n@rminfTD%
    \let\vecunit@=\vecunit@TD%
    \let\figcoord=\figcoordTD%
    \let\figgetangle=\figgetangleTD%
    \let\figpt=\figptTD%
    \let\figptBezier=\figptBezierTD%
    \let\figptbary=\figptbaryTD%
    \let\figptbaryR=\figptbaryRTD%
    \let\figptcirc=\figptcircTD%
    \let\figptcircumcenter=\figptcircumcenterTD%
    \let\figptcopy=\figptcopyTD%
    \let\figptcurvcenter=\figptcurvcenterTD%
    \let\figptinterlineplane=\figptinterlineplaneTD%
    \let\figptinterlines=\inters@cTD%
    \let\figptorthocenter=\figptorthocenterTD%
    \let\figptorthoprojline=\figptorthoprojlineTD%
    \let\figptorthoprojplane=\figptorthoprojplaneTD%
    \let\figptrot=\figptrotTD%
    \let\figptscontrol=\figptscontrolTD%
    \let\figptsorthoprojline=\figptsorthoprojlineTD%
    \let\figptsorthoprojplane=\figptsorthoprojplaneTD%
    \let\figptsrot=\figptsrotTD%
    \let\figptssym=\figptssymTD%
    \let\figptstra=\figptstraTD%
    \let\figptsym=\figptsymTD%
    \let\figpttraC=\figpttraCTD%
    \let\figpttra=\figpttraTD%
    \let\figsetobdist=\figsetobdistTD%
    \let\figsettarget=\figsettargetTD%
    \let\figsetview=\figsetviewTD%
    \let\figvectDBezier=\figvectDBezierTD%
    \let\figvectN=\figvectNTD%
    \let\figvectNV=\figvectNVTD%
    \let\figvectP=\figvectPTD%
    \let\figvectU=\figvectUTD%
    \let\psarccircP=\psarccircPTD%
    \let\psarccirc=\psarccircTD%
    \let\psarcell=\psarcellTD%
    \let\psarcellPA=\psarcellPATD%
    \let\psarrowBezier=\psarrowBezierTD%
    \let\psarrowcircP=\psarrowcircPTD%
    \let\psarrowcirc=\psarrowcircTD%
    \let\psarrowhead=\psarrowheadTD%
    \let\psarrow=\psarrowTD%
    \let\psBezier=\psBezierTD%
    \let\pscirc=\pscircTD%
    \let\pscurve=\pscurveTD%
    }
\def\un@v@ilable#1{\immediate\write16{*** The macro #1 is not available in the current context.}}
\def\figinsert#1{\figinsert@#1,:\ignorespaces}
\def\figinsert@#1,#2:{{\def\t@xt@{#2}\ifx\t@xt@\empty\xdef\Sc@leFact{1}\else%
    \def\Xarg@##1,{\def\@rgdeux{##1}}\Xarg@#2\xdef\Sc@leFact{\@rgdeux}\fi\@psfgetbb{#1}%
    \v@lX=\@psfllx\p@\v@lX=\ptpsT@pt\v@lX\v@lX=\Sc@leFact\v@lX%
    \v@lY=\@psflly\p@\v@lY=\ptpsT@pt\v@lY\v@lY=\Sc@leFact\v@lY%
    \b@undb@x{\v@lX}{\v@lY}%
    \v@lX=\@psfurx\p@\v@lX=\ptpsT@pt\v@lX\v@lX=\Sc@leFact\v@lX%
    \v@lY=\@psfury\p@\v@lY=\ptpsT@pt\v@lY\v@lY=\Sc@leFact\v@lY%
    \b@undb@x{\v@lX}{\v@lY}%
    \v@lX=100pt\v@lX=\Sc@leFact\v@lX\edef\F@ct{\repdecn@mb{\v@lX}}%
    \iftextures\special{postscriptfile #1 vscale=\F@ct\space hscale=\F@ct}%
    \else\includegraphics{#1}\fi%
    \message{[#1]}}\ignorespaces}
\def\figptDD#1:#2(#3,#4){\ifps@cri\c@ntr@lnum{#1}%
    {\v@lX=#3\unit@\v@lY=#4\unit@\let@xte={#2}
    \expandafter\xdef\csname\objc@de\endcsname{\ifitis@vect@r V\else P\fi%
    (\the\v@lX,\the\v@lY,\z@)\the\let@xte}}\ignorespaces\fi}
\def\figptTD#1:#2(#3,#4){\ifps@cri\c@ntr@lnum{#1}%
    \def\c@@rdYZ{#4,0,0}\extrairelepremi@r\c@@rdY\de\c@@rdYZ%
    \extrairelepremi@r\c@@rdZ\de\c@@rdYZ%
    {\v@lX=#3\unit@\v@lY=\c@@rdY\unit@\v@lZ=\c@@rdZ\unit@\let@xte={#2}%
    \expandafter\xdef\csname\objc@de\endcsname{\ifitis@vect@r V\else P\fi%
    (\the\v@lX,\the\v@lY,\the\v@lZ)\the\let@xte}%
    \b@undb@xTD{\v@lX}{\v@lY}{\v@lZ}}\ignorespaces\fi}
\def\figp@intregDD#1:#2(#3,#4){\c@ntr@lnum{#1}%
    {\result@t=#4\v@lX=#3\v@lY=\result@t\let@xte={#2}
    \expandafter\xdef\csname\objc@de\endcsname{\ifitis@vect@r V\else P\fi%
    (\the\v@lX,\the\v@lY,\z@)\the\let@xte}}\ignorespaces}
\def\figp@intregTD#1:#2(#3,#4){\c@ntr@lnum{#1}%
    \def\c@@rdYZ{#4,\z@,\z@}\extrairelepremi@r\c@@rdY\de\c@@rdYZ%
    \extrairelepremi@r\c@@rdZ\de\c@@rdYZ%
    {\v@lX=#3\v@lY=\c@@rdY\v@lZ=\c@@rdZ\let@xte={#2}%
    \expandafter\xdef\csname\objc@de\endcsname{\ifitis@vect@r V\else P\fi%
    (\the\v@lX,\the\v@lY,\the\v@lZ)\the\let@xte}%
    \b@undb@xTD{\v@lX}{\v@lY}{\v@lZ}}\ignorespaces}
\def\figptBezierDD#1:#2:#3[#4,#5,#6,#7]{\ifps@cri{\s@uvc@ntr@l\et@tfigptBezierDD%
    \figptBezier@#3[#4,#5,#6,#7]\figp@intregDD#1:{#2}(\v@lX,\v@lY)%
    \resetc@ntr@l\et@tfigptBezierDD}\ignorespaces\fi}
\def\figptBezierTD#1:#2:#3[#4,#5,#6,#7]{\ifps@cri{\s@uvc@ntr@l\et@tfigptBezierTD%
    \figptBezier@#3[#4,#5,#6,#7]\figp@intregTD#1:{#2}(\v@lX,\v@lY,\v@lZ)%
    \resetc@ntr@l\et@tfigptBezierTD}\ignorespaces\fi}
\def\figptBezier@#1[#2,#3,#4,#5]{\setc@ntr@l{2}%
    \edef\T@{#1}\v@leur=\p@\advance\v@leur-#1pt\edef\UNmT@{\repdecn@mb{\v@leur}}%
    \figptcopy-4:/#2/\figptcopy-3:/#3/\figptcopy-2:/#4/\figptcopy-1:/#5/%
    \s@mme=\z@\p@rtent=\z@\loop\ifnum\s@mme<\thr@@\advance\p@rtent\m@ne%
    \c@lDecast{-4}{-3}{\p@rtent}\advance\s@mme\@ne\repeat\figg@tXY{-4}}
\def\c@lDCUnDD#1#2{\figg@tXY{#1}\v@lX=\UNmT@\v@lX\v@lY=\UNmT@\v@lY%
    \figg@tXYa{#2}\advance\v@lX\T@\v@lXa\advance\v@lY\T@\v@lYa%
    \figp@intregDD#1:(\v@lX,\v@lY)}
\def\c@lDCUnTD#1#2{\figg@tXY{#1}\v@lX=\UNmT@\v@lX\v@lY=\UNmT@\v@lY\v@lZ=\UNmT@\v@lZ%
    \figg@tXYa{#2}\advance\v@lX\T@\v@lXa\advance\v@lY\T@\v@lYa\advance\v@lZ\T@\v@lZa%
    \figp@intregTD#1:(\v@lX,\v@lY,\v@lZ)}
\def\c@lDecast#1#2#3{\l@mbd@un=#1\l@mbd@de=#2%
    {\loop\ifnum\l@mbd@un<#3\c@lDCUn{\l@mbd@un}{\l@mbd@de}%
    \advance\l@mbd@un\@ne\advance\l@mbd@de\@ne\repeat}}
\def\figptbaryDD#1:#2[#3;#4]{\ifps@cri{\def\list@num{#3}\extrairelepremi@r\p@int\de\list@num%
    \s@mme=\z@\@ecfor\c@ef:=#4\do{\advance\s@mme\c@ef}%
    \def\listec@ef{#4,0}\extrairelepremi@r\c@ef\de\listec@ef%
    \figg@tXY{\p@int}\divide\v@lX\s@mme\divide\v@lY\s@mme%
    \multiply\v@lX\c@ef\multiply\v@lY\c@ef%
    \@ecfor\p@int:=\list@num\do{\extrairelepremi@r\c@ef\de\listec@ef%
           \figg@tXYa{\p@int}\divide\v@lXa\s@mme\divide\v@lYa\s@mme%
           \multiply\v@lXa\c@ef\multiply\v@lYa\c@ef%
           \advance\v@lX\v@lXa\advance\v@lY\v@lYa}%
    \figp@intregDD#1:{#2}(\v@lX,\v@lY)}\ignorespaces\fi}
\def\figptbaryTD#1:#2[#3;#4]{\ifps@cri{\def\list@num{#3}\extrairelepremi@r\p@int\de\list@num%
    \s@mme=\z@\@ecfor\c@ef:=#4\do{\advance\s@mme\c@ef}%
    \def\listec@ef{#4,0}\extrairelepremi@r\c@ef\de\listec@ef%
    \figg@tXY{\p@int}\divide\v@lX\s@mme\divide\v@lY\s@mme\divide\v@lZ\s@mme%
    \multiply\v@lX\c@ef\multiply\v@lY\c@ef\multiply\v@lZ\c@ef%
    \@ecfor\p@int:=\list@num\do{\extrairelepremi@r\c@ef\de\listec@ef%
           \figg@tXYa{\p@int}\divide\v@lXa\s@mme\divide\v@lYa\s@mme\divide\v@lZa\s@mme%
           \multiply\v@lXa\c@ef\multiply\v@lYa\c@ef\multiply\v@lZa\c@ef%
           \advance\v@lX\v@lXa\advance\v@lY\v@lYa\advance\v@lZ\v@lZa}%
    \figp@intregTD#1:{#2}(\v@lX,\v@lY,\v@lZ)}\ignorespaces\fi}
\def\figptbaryRDD#1:#2[#3;#4]{\ifps@cri{\def\list@num{#3}\extrairelepremi@r\p@int\de\list@num%
    \v@leur=\z@\@ecfor\c@ef:=#4\do{\advance\v@leur\c@ef pt}%
    \invers@\v@leur\v@leur\edef\invs@mme{\repdecn@mb{\v@leur}}%
    \def\listec@ef{#4,0}\extrairelepremi@r\c@ef\de\listec@ef%
    \figg@tXY{\p@int}\v@lX=\invs@mme\v@lX\v@lY=\invs@mme\v@lY%
    \v@lX=\c@ef\v@lX\v@lY=\c@ef\v@lY%
    \@ecfor\p@int:=\list@num\do{\extrairelepremi@r\c@ef\de\listec@ef%
           \figg@tXYa{\p@int}\v@lXa=\invs@mme\v@lXa\v@lYa=\invs@mme\v@lYa%
	   \v@lXa=\c@ef\v@lXa\v@lYa=\c@ef\v@lYa%
           \advance\v@lX\v@lXa\advance\v@lY\v@lYa}%
    \figp@intregDD#1:{#2}(\v@lX,\v@lY)}\ignorespaces\fi}
\def\figptbaryRTD#1:#2[#3;#4]{\ifps@cri{\def\list@num{#3}\extrairelepremi@r\p@int\de\list@num%
    \v@leur=\z@\@ecfor\c@ef:=#4\do{\advance\v@leur\c@ef pt}%
    \invers@\v@leur\v@leur\edef\invs@mme{\repdecn@mb{\v@leur}}%
    \def\listec@ef{#4,0}\extrairelepremi@r\c@ef\de\listec@ef%
    \figg@tXY{\p@int}\v@lX=\invs@mme\v@lX\v@lY=\invs@mme\v@lY\v@lZ=\invs@mme\v@lZ%
    \v@lX=\c@ef\v@lX\v@lY=\c@ef\v@lY\v@lZ=\c@ef\v@lZ%
    \@ecfor\p@int:=\list@num\do{\extrairelepremi@r\c@ef\de\listec@ef%
           \figg@tXYa{\p@int}%
	   \v@lXa=\invs@mme\v@lXa\v@lYa=\invs@mme\v@lYa\v@lZa=\invs@mme\v@lZa%
	   \v@lXa=\c@ef\v@lXa\v@lYa=\c@ef\v@lYa\v@lZa=\c@ef\v@lZa%
           \advance\v@lX\v@lXa\advance\v@lY\v@lYa\advance\v@lZ\v@lZa}%
    \figp@intregTD#1:{#2}(\v@lX,\v@lY,\v@lZ)}\ignorespaces\fi}
\def\figptcircDD#1:#2:#3;#4(#5){\ifps@cri{\s@uvc@ntr@l\et@tfigptcircDD%
    \c@lptellDD#1:{#2}:#3;#4,#4(#5)\resetc@ntr@l\et@tfigptcircDD}\ignorespaces\fi}
\def\figptcircTD#1:#2:#3,#4,#5;#6(#7){\ifps@cri{\s@uvc@ntr@l\et@tfigptcircTD%
    \setc@ntr@l{2}\c@lExtAxes#3,#4,#5(#6)\figptellP#1:{#2}:#3,-4,-5(#7)%
    \resetc@ntr@l\et@tfigptcircTD}\ignorespaces\fi}
\def\figptcircumcenterDD#1:#2[#3,#4,#5]{\ifps@cri{\s@uvc@ntr@l\et@tfigptcircumcenterDD%
    \setc@ntr@l{2}\figvectNDD-5[#3,#4]\figptbaryDD-3:[#3,#4;1,1]%
                  \figvectNDD-6[#4,#5]\figptbaryDD-4:[#4,#5;1,1]%
    \resetc@ntr@l{2}\inters@cDD#1:{#2}[-3,-5;-4,-6]%
    \resetc@ntr@l\et@tfigptcircumcenterDD}\ignorespaces\fi}
\def\figptcircumcenterTD#1:#2[#3,#4,#5]{\ifps@cri{\s@uvc@ntr@l\et@tfigptcircumcenterTD%
    \setc@ntr@l{2}\figvectNTD-1[#3,#4,#5]%
    \figvectPTD-3[#3,#4]\figvectNVTD-5[-1,-3]\figptbaryTD-3:[#3,#4;1,1]%
    \figvectPTD-4[#4,#5]\figvectNVTD-6[-1,-4]\figptbaryTD-4:[#4,#5;1,1]%
    \resetc@ntr@l{2}\inters@cTD#1:{#2}[-3,-5;-4,-6]%
    \resetc@ntr@l\et@tfigptcircumcenterTD}\ignorespaces\fi}
\def\figptcopyDD#1:#2/#3/{\ifps@cri{\figg@tXY{#3}%
    \figp@intregDD#1:{#2}(\v@lX,\v@lY)}\ignorespaces\fi}
\def\figptcopyTD#1:#2/#3/{\ifps@cri{\figg@tXY{#3}%
    \figp@intregTD#1:{#2}(\v@lX,\v@lY,\v@lZ)}\ignorespaces\fi}
\def\figptcurvcenterDD#1:#2:#3[#4,#5,#6,#7]{\ifps@cri{\s@uvc@ntr@l\et@tfigptcurvcenterDD%
    \setc@ntr@l{2}\c@lcurvradDD#3[#4,#5,#6,#7]\edef\Sprim@{\repdecn@mb{\result@t}}%
    \figptBezierDD-1::#3[#4,#5,#6,#7]\figpttraDD#1:{#2}=-1/\Sprim@,-5/%
    \resetc@ntr@l\et@tfigptcurvcenterDD}\ignorespaces\fi}
\def\figptcurvcenterTD#1:#2:#3[#4,#5,#6,#7]{\ifps@cri{\s@uvc@ntr@l\et@tfigptcurvcenterTD%
    \setc@ntr@l{2}\figvectDBezierTD -5:1,#3[#4,#5,#6,#7]%
    \figvectDBezierTD -6:2,#3[#4,#5,#6,#7]\vecunit@TD{-5}{-5}%
    \edef\Sprim@{\repdecn@mb{\result@t}}\figvectNVTD-1[-6,-5]%
    \figvectNVTD-5[-5,-1]\c@lproscalTD\v@leur[-6,-5]%
    \invers@{\v@leur}{\v@leur}\v@leur=\Sprim@\v@leur\v@leur=\Sprim@\v@leur%
    \figptBezierTD-1::#3[#4,#5,#6,#7]\edef\Sprim@{\repdecn@mb{\v@leur}}%
    \figpttraTD#1:{#2}=-1/\Sprim@,-5/\resetc@ntr@l\et@tfigptcurvcenterTD}\ignorespaces\fi}
\def\c@lcurvradDD#1[#2,#3,#4,#5]{{\figvectDBezierDD -5:1,#1[#2,#3,#4,#5]%
    \figvectDBezierDD -6:2,#1[#2,#3,#4,#5]\vecunit@DD{-5}{-5}%
    \edef\Sprim@{\repdecn@mb{\result@t}}\figvectNVDD-5[-5]\c@lproscalDD\v@leur[-6,-5]%
    \invers@{\v@leur}{\v@leur}\v@leur=\Sprim@\v@leur\v@leur=\Sprim@\v@leur%
    \global\result@t=\v@leur}}
\def\figptellDD#1:#2:#3;#4,#5(#6,#7){\ifps@cri{\s@uvc@ntr@l\et@tfigptell%
    \c@lptellDD#1::#3;#4,#5(#6)\figptrotDD#1:{#2}=#1/#3,#7/%
    \resetc@ntr@l\et@tfigptell}\ignorespaces\fi}
\def\c@lptellDD#1:#2:#3;#4,#5(#6){\c@ssin{\v@lmin}{\v@lmax}{#6pt}%
    \v@lmin=#4\v@lmin\v@lmax=#5\v@lmax%
    \edef\Xc@mp{\repdecn@mb{\v@lmin}}\edef\Yc@mp{\repdecn@mb{\v@lmax}}%
    \setc@ntr@l{2}\figvectC-1(\Xc@mp,\Yc@mp)\figpttraDD#1:{#2}=#3/1,-1/}
\def\figptellP#1:#2:#3,#4,#5(#6){\ifps@cri{\s@uvc@ntr@l\et@tfigptellP%
    \setc@ntr@l{2}\figvectP-1[#3,#4]\figvectP-2[#3,#5]%
    \v@leur=#6pt\c@lptellP{#3}{-1}{-2}\figptcopy#1:{#2}/-3/%
    \resetc@ntr@l\et@tfigptellP}\ignorespaces\fi}
\def\c@lptellP#1#2#3{\c@ssin{\v@lmin}{\v@lmax}{\v@leur}%
    \edef\C@{\repdecn@mb{\v@lmin}}\edef\S@{\repdecn@mb{\v@lmax}}%
    \figpttra-3:=#1/\C@,#2/\figpttra-3:=-3/\S@,#3/}
\def\figptendnormalDD#1:#2:#3,#4[#5,#6]{\ifps@cri{\s@uvc@ntr@l\et@tfigptendnormal%
    \figg@tXYa{#5}\figg@tXY{#6}%
    \advance\v@lX-\v@lXa\advance\v@lY-\v@lYa%
    \setc@ntr@l{2}\figv@ctCreg-1(\v@lX,\v@lY)\vecunit@{-1}{-1}\figg@tXY{-1}%
    \delt@=#3\unit@\maxim@m{\delt@}{\delt@}{-\delt@}\edef\l@ngueur{\repdecn@mb{\delt@}}%
    \v@lX=\l@ngueur\v@lX\v@lY=\l@ngueur\v@lY%
    \delt@=\p@\advance\delt@-#4pt\edef\l@ngueur{\repdecn@mb{\delt@}}%
    \figptbaryRDD-1:[#5,#6;#4,\l@ngueur]\figg@tXYa{-1}%
    \advance\v@lXa\v@lY\advance\v@lYa-\v@lX%
    \setc@ntr@l{1}\figp@intregDD#1:{#2}(\v@lXa,\v@lYa)\resetc@ntr@l\et@tfigptendnormal}%
    \ignorespaces\fi}
\def\figptinscribedcenter#1:#2[#3,#4,#5]{\ifps@cri{%
    \figgetdist\LA@[#4,#5]\figgetdist\LB@[#3,#5]\figgetdist\LC@[#3,#4]%
    \figptbaryR#1:{#2}[#3,#4,#5;\LA@,\LB@,\LC@]}\ignorespaces\fi}
\def\figptinterlineplaneTD#1:#2[#3,#4;#5,#6]{\ifps@cri{\s@uvc@ntr@l\et@tfigptinterlineplane%
    \setc@ntr@l{2}\figvectPTD-1[#3,#5]\vecunit@TD{-2}{#6}%
    \r@pPSTD\v@leur[-2,-1,#4]\edef\v@lcoef{\repdecn@mb{\v@leur}}%
    \figpttraTD#1:{#2}=#3/\v@lcoef,#4/\resetc@ntr@l\et@tfigptinterlineplane}\ignorespaces\fi}
\def\figptorthocenterDD#1:#2[#3,#4,#5]{\ifps@cri{\s@uvc@ntr@l\et@tfigptorthocenterDD%
    \setc@ntr@l{2}\figvectNDD-3[#3,#4]\figvectNDD-4[#4,#5]%
    \resetc@ntr@l{2}\inters@cDD#1:{#2}[#5,-3;#3,-4]%
    \resetc@ntr@l\et@tfigptorthocenterDD}\ignorespaces\fi}
\def\figptorthocenterTD#1:#2[#3,#4,#5]{\ifps@cri{\s@uvc@ntr@l\et@tfigptorthocenterTD%
    \setc@ntr@l{2}\figvectNTD-1[#3,#4,#5]%
    \figvectPTD-2[#3,#4]\figvectNVTD-3[-1,-2]%
    \figvectPTD-2[#4,#5]\figvectNVTD-4[-1,-2]%
    \resetc@ntr@l{2}\inters@cTD#1:{#2}[#5,-3;#3,-4]%
    \resetc@ntr@l\et@tfigptorthocenterTD}\ignorespaces\fi}
\def\figptorthoprojlineDD#1:#2=#3/#4,#5/{\ifps@cri{\s@uvc@ntr@l\et@tfigptorthoprojlineDD%
    \setc@ntr@l{2}\figvectPDD-3[#4,#5]\figvectNVDD-4[-3]\resetc@ntr@l{2}%
    \inters@cDD#1:{#2}[#3,-4;#4,-3]\resetc@ntr@l\et@tfigptorthoprojlineDD}\ignorespaces\fi}
\def\figptorthoprojlineTD#1:#2=#3/#4,#5/{\ifps@cri{\s@uvc@ntr@l\et@tfigptorthoprojlineTD%
    \setc@ntr@l{2}\figvectPTD-1[#4,#3]\figvectPTD-2[#4,#5]\vecunit@TD{-2}{-2}%
    \c@lproscalTD\v@leur[-1,-2]\edef\v@lcoef{\repdecn@mb{\v@leur}}%
    \figpttraTD#1:{#2}=#4/\v@lcoef,-2/\resetc@ntr@l\et@tfigptorthoprojlineTD}\ignorespaces\fi}
\def\figptorthoprojplaneTD#1:#2=#3/#4,#5/{\ifps@cri{\s@uvc@ntr@l\et@tfigptorthoprojplane%
    \setc@ntr@l{2}\figvectPTD-1[#3,#4]\vecunit@TD{-2}{#5}%
    \c@lproscalTD\v@leur[-1,-2]\edef\v@lcoef{\repdecn@mb{\v@leur}}%
    \figpttraTD#1:{#2}=#3/\v@lcoef,-2/\resetc@ntr@l\et@tfigptorthoprojplane}\ignorespaces\fi}
\def\figpthom#1:#2=#3/#4,#5/{\ifps@cri{\s@uvc@ntr@l\et@tfigpthom%
    \setc@ntr@l{2}\figvectP-1[#4,#3]\figpttra#1:{#2}=#4/#5,-1/%
    \resetc@ntr@l\et@tfigpthom}\ignorespaces\fi}
\def\figptrotDD#1:#2=#3/#4,#5/{\ifps@cri{\s@uvc@ntr@l\et@tfigptrotDD%
    \c@ssin{\v@lmin}{\v@lmax}{#5pt}%
    \edef\C@{\repdecn@mb{\v@lmin}}\edef\S@{\repdecn@mb{\v@lmax}}%
    \setc@ntr@l{2}\figvectPDD-1[#4,#3]\figg@tXY{-1}%
    \v@lXa=\C@\v@lX\v@lYa=-\S@\v@lY\advance\v@lXa\v@lYa%
    \v@lYa=\S@\v@lX\v@lX=\C@\v@lY\advance\v@lYa\v@lX%
    \figv@ctCreg-1(\v@lXa,\v@lYa)\figpttraDD#1:{#2}=#4/1,-1/%
    \resetc@ntr@l\et@tfigptrotDD}\ignorespaces\fi}
\def\figptrotTD#1:#2=#3/#4,#5,#6/{\ifps@cri{\s@uvc@ntr@l\et@tfigptrotTD%
    \c@ssin{\v@lmin}{\v@lmax}{#5pt}%
    \edef\C@{\repdecn@mb{\v@lmin}}\edef\S@{\repdecn@mb{\v@lmax}}%
    \setc@ntr@l{2}\figptorthoprojplaneTD-3:=#4/#3,#6/\figvectPTD-2[-3,#3]%
    \n@rmeucTD\v@leur{-2}\ifdim\v@leur<\Cepsil@n\figg@tXYa{#3}\else%
    \edef\v@lcoef{\repdecn@mb{\v@leur}}\figvectNVTD-1[#6,-2]%
    \figg@tXYa{-1}\v@lXa=\v@lcoef\v@lXa\v@lYa=\v@lcoef\v@lYa\v@lZa=\v@lcoef\v@lZa%
    \v@lXa=\S@\v@lXa\v@lYa=\S@\v@lYa\v@lZa=\S@\v@lZa\figg@tXY{-2}%
    \advance\v@lXa\C@\v@lX\advance\v@lYa\C@\v@lY\advance\v@lZa\C@\v@lZ%
    \figg@tXY{-3}\advance\v@lXa\v@lX\advance\v@lYa\v@lY\advance\v@lZa\v@lZ\fi%
    \figp@intregTD#1:{#2}(\v@lXa,\v@lYa,\v@lZa)\resetc@ntr@l\et@tfigptrotTD}\ignorespaces\fi}
\def\figptsymDD#1:#2=#3/#4,#5/{\ifps@cri{\s@uvc@ntr@l\et@tfigptsymDD%
    \resetc@ntr@l{2}\figptorthoprojlineDD-5:=#3/#4,#5/\figvectPDD-2[#3,-5]%
    \figpttraDD#1:{#2}=#3/2,-2/\resetc@ntr@l\et@tfigptsymDD}\ignorespaces\fi}
\def\figptsymTD#1:#2=#3/#4,#5/{\ifps@cri{\s@uvc@ntr@l\et@tfigptsymTD%
    \resetc@ntr@l{2}\figptorthoprojplaneTD-3:=#3/#4,#5/\figvectPTD-2[#3,-3]%
    \figpttraTD#1:{#2}=#3/2,-2/\resetc@ntr@l\et@tfigptsymTD}\ignorespaces\fi}
\def\figpttraDD#1:#2=#3/#4,#5/{\ifps@cri{\figg@tXYa{#5}\v@lXa=#4\v@lXa\v@lYa=#4\v@lYa%
    \figg@tXY{#3}\advance\v@lX\v@lXa\advance\v@lY\v@lYa%
    \figp@intregDD#1:{#2}(\v@lX,\v@lY)}\ignorespaces\fi}
\def\figpttraTD#1:#2=#3/#4,#5/{\ifps@cri{\figg@tXYa{#5}\v@lXa=#4\v@lXa\v@lYa=#4\v@lYa%
    \v@lZa=#4\v@lZa\figg@tXY{#3}\advance\v@lX\v@lXa\advance\v@lY\v@lYa%
    \advance\v@lZ\v@lZa\figp@intregTD#1:{#2}(\v@lX,\v@lY,\v@lZ)}\ignorespaces\fi}
\def\figpttraCDD#1:#2=#3/#4,#5/{\ifps@cri{\v@lXa=#4\unit@\v@lYa=#5\unit@%
    \figg@tXY{#3}\advance\v@lX\v@lXa\advance\v@lY\v@lYa%
    \figp@intregDD#1:{#2}(\v@lX,\v@lY)}\ignorespaces\fi}
\def\figpttraCTD#1:#2=#3/#4,#5,#6/{\ifps@cri{\v@lXa=#4\unit@\v@lYa=#5\unit@\v@lZa=#6\unit@%
    \figg@tXY{#3}\advance\v@lX\v@lXa\advance\v@lY\v@lYa\advance\v@lZ\v@lZa%
    \figp@intregTD#1:{#2}(\v@lX,\v@lY,\v@lZ)}\ignorespaces\fi}
\def\figptscontrolDD#1[#2,#3,#4,#5]{\ifps@cri{\s@uvc@ntr@l\et@tfigptscontrolDD%
    \setc@ntr@l{2}\figptcopyDD-2:/#2/\s@mme=#1\advance\s@mme\@ne%
    \figg@tXY{-2}\v@lX=8\v@lX\v@lY=8\v@lY\figp@intregDD-1:(\v@lX,\v@lY)%
    \figpttraDD-1:=-1/-27,#3/\figpttraDD-1:=-1/1,#5/%
    \figpttraDD-2:=-2/-27,#4/\figpttraDD-2:=-2/8,#5/%
    \c@lptCtrlDD#1[-2,-1]\c@lptCtrlDD\the\s@mme[-1,-2]%
    \resetc@ntr@l\et@tfigptscontrolDD}\ignorespaces\fi}
\def\figptscontrolTD#1[#2,#3,#4,#5]{\ifps@cri{\s@uvc@ntr@l\et@tfigptscontrolTD%
    \setc@ntr@l{2}\figptcopyTD-2:/#2/\s@mme=#1\advance\s@mme\@ne%
    \figg@tXY{-2}\v@lX=8\v@lX\v@lY=8\v@lY\v@lZ=8\v@lZ\figp@intregTD-1:(\v@lX,\v@lY,\v@lZ)%
    \figpttraTD-1:=-1/-27,#3/\figpttraTD-1:=-1/1,#5/%
    \figpttraTD-2:=-2/-27,#4/\figpttraTD-2:=-2/8,#5/%
    \c@lptCtrlTD#1[-2,-1]\c@lptCtrlTD\the\s@mme[-1,-2]%
    \resetc@ntr@l\et@tfigptscontrolTD}\ignorespaces\fi}
\def\figptsintercircDD#1[#2,#3;#4,#5]{\ifps@cri{\s@uvc@ntr@l\et@tfigptsintercirc%
    \setc@ntr@l{2}\figvectP-1[#2,#4]\n@rmeucDD{\delt@}{-1}\s@mme=#1\advance\s@mme\@ne%
    \figptcopyDD#1:/#2/\figptcopyDD\the\s@mme:/#4/%
    \ifdim\delt@=\z@\else%
    \v@lmin=#3\unit@\v@lmax=#5\unit@\v@leur=\v@lmin\advance\v@leur\v@lmax%
    \ifdim\v@leur>\delt@%
    \v@leur=\v@lmin\advance\v@leur-\v@lmax\maxim@m{\v@leur}{\v@leur}{-\v@leur}%
    \ifdim\v@leur<\delt@%
    \v@lmin=\repdecn@mb{\v@lmin}\v@lmin\v@lmax=\repdecn@mb{\v@lmax}\v@lmax%
    \invers@{\v@leur}{\delt@}\advance\v@lmax-\v@lmin%
    \v@lmax=-\repdecn@mb{\v@leur}\v@lmax\advance\delt@\v@lmax\delt@=.5\delt@%
    \setc@ntr@l{2}\figvectPDD-1[#2,#4]\vecunit@{-1}{-1}%
    \edef\t@ille{\repdecn@mb{\delt@}}\figpttraDD-2:=#2/\t@ille,-1/%
    \figg@tXY{-1}\figv@ctCreg-1(-\v@lY,\v@lX)%
    \delt@=\repdecn@mb{\delt@}\delt@\advance\v@lmin-\delt@%
    \sqrt@{\v@leur}{\v@lmin}\edef\t@ille{\repdecn@mb{\v@leur}}%
    \figpttraDD#1:=-2/-\t@ille,-1/\figpttraDD\the\s@mme:=-2/\t@ille,-1/\fi\fi\fi%
    \resetc@ntr@l\et@tfigptsintercirc}\ignorespaces\fi}
\def\figptsorthoprojlineDD#1=#2/#3,#4/{\ifps@cri{\s@uvc@ntr@l\et@tfigptsorthoprojlineDD%
    \setc@ntr@l{2}\figvectPDD-3[#3,#4]\figvectNVDD-4[-3]\resetc@ntr@l{2}%
    \def\list@num{#2}\s@mme=#1\@ecfor\p@int:=\list@num\do{%
    \inters@cDD\the\s@mme:[\p@int,-4;#3,-3]\advance\s@mme\@ne}%
    \resetc@ntr@l\et@tfigptsorthoprojlineDD}\ignorespaces\fi}
\def\figptsorthoprojlineTD#1=#2/#3,#4/{\ifps@cri{\s@uvc@ntr@l\et@tfigptsorthoprojlineTD%
    \setc@ntr@l{2}\figvectPTD-2[#3,#4]\vecunit@TD{-2}{-2}%
    \def\list@num{#2}\s@mme=#1\@ecfor\p@int:=\list@num\do{%
    \figvectPTD-1[#3,\p@int]\c@lproscalTD\v@leur[-1,-2]%
    \edef\v@lcoef{\repdecn@mb{\v@leur}}\figpttraTD\the\s@mme:=#3/\v@lcoef,-2/%
    \advance\s@mme\@ne}\resetc@ntr@l\et@tfigptsorthoprojlineTD}\ignorespaces\fi}
\def\figptsorthoprojplaneTD#1=#2/#3,#4/{\ifps@cri{\s@uvc@ntr@l\et@tfigptsorthoprojplane%
    \setc@ntr@l{2}\vecunit@TD{-2}{#4}%
    \def\list@num{#2}\s@mme=#1\@ecfor\p@int:=\list@num\do{\figvectPTD-1[\p@int,#3]%
    \c@lproscalTD\v@leur[-1,-2]\edef\v@lcoef{\repdecn@mb{\v@leur}}%
    \figpttraTD\the\s@mme:=\p@int/\v@lcoef,-2/\advance\s@mme\@ne}%
    \resetc@ntr@l\et@tfigptsorthoprojplane}\ignorespaces\fi}
\def\figptshom#1=#2/#3,#4/{\ifps@cri{\s@uvc@ntr@l\et@tfigptshom%
    \setc@ntr@l{2}\def\list@num{#2}\s@mme=#1%
    \@ecfor\p@int:=\list@num\do{\figvectP-1[#3,\p@int]%
    \figpttra\the\s@mme:=#3/#4,-1/\advance\s@mme\@ne}%
    \resetc@ntr@l\et@tfigptshom}\ignorespaces\fi}
\def\figptsrotDD#1=#2/#3,#4/{\ifps@cri{\s@uvc@ntr@l\et@tfigptsrotDD%
    \c@ssin{\v@lmin}{\v@lmax}{#4pt}%
    \edef\C@{\repdecn@mb{\v@lmin}}\edef\S@{\repdecn@mb{\v@lmax}}%
    \setc@ntr@l{2}\def\list@num{#2}\s@mme=#1%
    \@ecfor\p@int:=\list@num\do{\figvectPDD-1[#3,\p@int]\figg@tXY{-1}%
    \v@lXa=\C@\v@lX\v@lYa=-\S@\v@lY\advance\v@lXa\v@lYa%
    \v@lYa=\S@\v@lX\v@lX=\C@\v@lY\advance\v@lYa\v@lX%
    \figv@ctCreg-1(\v@lXa,\v@lYa)\figpttraDD\the\s@mme:=#3/1,-1/\advance\s@mme\@ne}%
    \resetc@ntr@l\et@tfigptsrotDD}\ignorespaces\fi}
\def\figptsrotTD#1=#2/#3,#4,#5/{\ifps@cri{\s@uvc@ntr@l\et@tfigptsrotTD%
    \c@ssin{\v@lmin}{\v@lmax}{#4pt}%
    \edef\C@{\repdecn@mb{\v@lmin}}\edef\S@{\repdecn@mb{\v@lmax}}%
    \setc@ntr@l{2}\def\list@num{#2}\s@mme=#1%
    \@ecfor\p@int:=\list@num\do{\figptorthoprojplaneTD-3:=#3/\p@int,#5/%
    \figvectPTD-2[-3,\p@int]%
    \figvectNVTD-1[#5,-2]\n@rmeucTD\v@leur{-2}\edef\v@lcoef{\repdecn@mb{\v@leur}}%
    \figg@tXYa{-1}\v@lXa=\v@lcoef\v@lXa\v@lYa=\v@lcoef\v@lYa\v@lZa=\v@lcoef\v@lZa%
    \v@lXa=\S@\v@lXa\v@lYa=\S@\v@lYa\v@lZa=\S@\v@lZa\figg@tXY{-2}%
    \advance\v@lXa\C@\v@lX\advance\v@lYa\C@\v@lY\advance\v@lZa\C@\v@lZ%
    \figg@tXY{-3}\advance\v@lXa\v@lX\advance\v@lYa\v@lY\advance\v@lZa\v@lZ%
    \figp@intregTD\the\s@mme:(\v@lXa,\v@lYa,\v@lZa)\advance\s@mme\@ne}%
    \resetc@ntr@l\et@tfigptsrotTD}\ignorespaces\fi}
\def\figptssymDD#1=#2/#3,#4/{\ifps@cri{\s@uvc@ntr@l\et@tfigptssymDD%
    \setc@ntr@l{2}\figvectPDD-3[#3,#4]\figg@tXY{-3}\figv@ctCreg-4(-\v@lY,\v@lX)%
    \resetc@ntr@l{2}\def\list@num{#2}\s@mme=#1%
    \@ecfor\p@int:=\list@num\do{\inters@cDD-5:[#3,-3;\p@int,-4]\figvectPDD-2[\p@int,-5]%
    \figpttraDD\the\s@mme:=\p@int/2,-2/\advance\s@mme\@ne}%
    \resetc@ntr@l\et@tfigptssymDD}\ignorespaces\fi}
\def\figptssymTD#1=#2/#3,#4/{\ifps@cri{\s@uvc@ntr@l\et@tfigptssymTD%
    \setc@ntr@l{2}\vecunit@TD{-2}{#4}\def\list@num{#2}\s@mme=#1%
    \@ecfor\p@int:=\list@num\do{\figvectPTD-1[\p@int,#3]%
    \c@lproscalTD\v@leur[-1,-2]\v@leur=2\v@leur\edef\v@lcoef{\repdecn@mb{\v@leur}}%
    \figpttraTD\the\s@mme:=\p@int/\v@lcoef,-2/\advance\s@mme\@ne}%
    \resetc@ntr@l\et@tfigptssymTD}\ignorespaces\fi}
\def\figptstraDD#1=#2/#3,#4/{\ifps@cri{\figg@tXYa{#4}\v@lXa=#3\v@lXa\v@lYa=#3\v@lYa%
    \def\list@num{#2}\s@mme=#1\@ecfor\p@int:=\list@num\do{\figg@tXY{\p@int}%
    \advance\v@lX\v@lXa\advance\v@lY\v@lYa%
    \figp@intregDD\the\s@mme:(\v@lX,\v@lY)\advance\s@mme\@ne}}\ignorespaces\fi}
\def\figptstraTD#1=#2/#3,#4/{\ifps@cri{\figg@tXYa{#4}\v@lXa=#3\v@lXa\v@lYa=#3\v@lYa%
    \v@lZa=#3\v@lZa\def\list@num{#2}\s@mme=#1\@ecfor\p@int:=\list@num\do{\figg@tXY{\p@int}%
    \advance\v@lX\v@lXa\advance\v@lY\v@lYa\advance\v@lZ\v@lZa%
    \figp@intregTD\the\s@mme:(\v@lX,\v@lY,\v@lZ)\advance\s@mme\@ne}}\ignorespaces\fi}
\def\figscan#1(#2,#3){{\s@uvc@ntr@l\et@tfigscan\@psfgetbb{#1}%
    \if@psfbbfound\else\def\@psfllx{0}\def\@psflly{50}\def\@psfurx{540}\def\@psfury{800}\fi%
    \unit@=\@ne bp\setc@ntr@l{2}\figsetmark{}%
    \def\minst@p{20pt}\v@lX=\@psfllx\p@\v@lX=\Sc@leFact\v@lX\r@undint\v@lX\v@lX%
    \v@lY=\@psflly\p@\v@lY=\Sc@leFact\v@lY%
    \delt@=\@psfury\p@\delt@=\Sc@leFact\delt@%
    \advance\delt@-\v@lY\v@lXa=\@psfurx\p@\v@lXa=\Sc@leFact\v@lXa\v@leur=\minst@p%
    \edef\valv@lY{\repdecn@mb{\v@lY}}\edef\LgTr@it{\the\delt@}%
    \loop\ifdim\v@lX<\v@lXa\edef\valv@lX{\repdecn@mb{\v@lX}}%
    \figptDD -1:(\valv@lX,\valv@lY)\figwriten -1:\hbox{\vrule height\LgTr@it}(0)%
    \ifdim\v@leur<\minst@p\else\figsetmark{$\scriptscriptstyle\triangle$}%
    \figwrites -1:\@ffichnb{0}{\valv@lX}(4)\v@leur=\z@\figsetmark{}\fi%
    \advance\v@leur#2pt\advance\v@lX#2pt\repeat%
    \def\minst@p{10pt}\v@lX=\@psfllx\p@\v@lX=\Sc@leFact\v@lX%
    \v@lY=\@psflly\p@\v@lY=\Sc@leFact\v@lY\r@undint\v@lY\v@lY%
    \delt@=\@psfurx\p@\delt@=\Sc@leFact\delt@%
    \advance\delt@-\v@lX\v@lYa=\@psfury\p@\v@lYa=\Sc@leFact\v@lYa\v@leur=\minst@p%
    \edef\valv@lX{\repdecn@mb{\v@lX}}\edef\LgTr@it{\the\delt@}%
    \loop\ifdim\v@lY<\v@lYa\edef\valv@lY{\repdecn@mb{\v@lY}}%
    \figptDD -1:(\valv@lX,\valv@lY)\figwritee -1:\vbox{\hrule width\LgTr@it}(0)%
    \ifdim\v@leur<\minst@p\else\figsetmark{$\triangleright$}%
    \figwritew -1:\@ffichnb{0}{\valv@lY}(4)\v@leur=\z@\figsetmark{}\fi%
    \advance\v@leur#3pt\advance\v@lY#3pt\repeat%
    \resetc@ntr@l\et@tfigscan}\ignorespaces}
\def\figshowpts[#1,#2]{{\figsetmark{$\bullet$}\figsetptname{\bf ##1}%
    \p@rtent=#2\relax\ifnum\p@rtent<\z@\p@rtent=\z@\fi%
    \s@mme=#1\relax\ifnum\s@mme<\z@\s@mme=\z@\fi%
    \loop\ifnum\s@mme<\p@rtent\pt@rvect{\s@mme}%
    \ifitis@K\figwriten{\the\s@mme}:(4pt)\fi\advance\s@mme\@ne\repeat%
    \pt@rvect{\s@mme}\ifitis@K\figwriten{\the\s@mme}:(4pt)\fi}\ignorespaces}
\def\pt@rvect#1{\set@bjc@de{#1}%
    \expandafter\expandafter\expandafter\inqpt@rvec\csname\objc@de\endcsname:}
\def\inqpt@rvec#1#2:{\if#1P\itis@vect@rtrue\else\itis@vect@rfalse\fi}
\def\figshowsettings{{%
    \immediate\write16{====================================================================}%
    \immediate\write16{ Current settings about:}%
    \immediate\write16{ --- GENERAL ---}%
    \immediate\write16{Scale factor and Unit = \unit@util\space (\the\unit@)
     \space -> \BS@ figinit{ScaleFactorUnit}}%
    \immediate\write16{Update mode = \ifpstestm@de yes\else no\fi
     \space-> \BS@ pssetupdate{yes/no}}%
    \immediate\write16{ --- PRINTING ---}%
    \immediate\write16{Implicit point name = \ptn@me{i} \space-> \BS@ figsetptname{Name}}%
    \immediate\write16{Point marker = \the\c@nsymb \space -> \BS@ figsetmark{Mark}}%
    \immediate\write16{Print rounded coordinates = \ifr@undcoord yes\else no\fi
     \space-> \BS@ figsetroundcoord{yes/no}}%
    \immediate\write16{ --- POSTSCRIPT ---}%
    \immediate\write16{Arrow-head:}%
    \immediate\write16{ (half-)Angle = \@rrowheadangle
     \space-> \BS@ pssetarrowheadangle{Angle}}%
    \immediate\write16{ Filling mode = \if@rrowhfill yes\else no\fi
     \space-> \BS@ pssetarrowheadfill{yes/no}}%
    \immediate\write16{ "Outside" = \if@rrowhout yes\else no\fi
     \space-> \BS@ pssetarrowheadout{yes/no}}%
    \immediate\write16{ Length = \@rrowheadlength
     \if@rrowratio\space(not active)\else\space(active)\fi
     \space-> \BS@ pssetarrowheadlength{Length}}%
    \immediate\write16{ Ratio = \@rrowheadratio
     \if@rrowratio\space(active)\else\space(not active)\fi
     \space-> \BS@ pssetarrowheadratio{Ratio}}%
    \immediate\write16{Color = \curr@ntcolor \space-> \BS@ pssetcmyk{\BS@ ColorName} or}%
    \immediate\write16{ \space\space\space \BS@ pssetrgb{\BS@ ColorName}
     or \BS@ pssetgray{GrayLevel}}%
    \immediate\write16{Curve roundness = \curv@roundness
     \space-> \BS@ pssetroundness{Roundness}}%
    \immediate\write16{Filling mode = \iffillm@de yes\else no\fi
     \space-> \BS@ pssetfillmode{yes/no}}%
    \immediate\write16{Line style = \curr@ntdash \space-> \BS@ pssetdash{Index/Pattern}}%
    \immediate\write16{Line width = \curr@ntwidth \space-> \BS@ pssetwidth{Width}}%
    \immediate\write16{Secondary settings:}%
    \immediate\write16{ Color = \sec@ndcolor \space-> \BS@ pssetsecondcmyk{\BS@ ColorName} or}%
    \immediate\write16{ \space\space\space \BS@ pssetsecondrgb{\BS@ ColorName}
     or \BS@ pssetsecondgray{GrayLevel}}%
    \immediate\write16{ Line style = \curr@ntseconddash \space-> \BS@ pssetseconddash{Index/Pattern}}%
    \immediate\write16{ Line width = \curr@ntsecondwidth \space-> \BS@ pssetsecondwidth{Width}}%
    \ifTr@isDim%
    \immediate\write16{ --- 3D to 2D PROJECTION ---}%
    \immediate\write16{Projection : \typ@proj \space-> \BS@ figinit{ScaleFactorUnit, ProjType}}%
    \immediate\write16{Psi = \v@lPsi \space-> \BS@ figsetview(Psi)}%
    \ifcase\curr@ntproj\let@xte={Lambda}\else\let@xte={Theta}\fi%
    \immediate\write16{\the\let@xte\space = \v@lTheta \space-> \BS@ figsetview(Psi, \the\let@xte)}%
    \ifnum\curr@ntproj=\tw@%
    \immediate\write16{Observation distance = \disob@unit \space-> \BS@ figsetobdist (Dist)}%
     \v@lX=\ptT@unit@\Xt@rget\v@lY=\ptT@unit@\Yt@rget\v@lZ=\ptT@unit@\Zt@rget%
    \immediate\write16{Target point = (\repdecn@mb{\v@lX}, \repdecn@mb{\v@lY},
     \repdecn@mb{\v@lZ}) \space-> \BS@ figsettarget [Pt]}\fi%
    \fi%
    \immediate\write16{====================================================================}%
    \ignorespaces}}
{\catcode`\/=0 \catcode`/\=12 /gdef/BS@{\}}
\newif\ifitis@vect@r
\def\figvectC#1(#2,#3){{\itis@vect@rtrue\figpt#1:(#2,#3)}\ignorespaces}
\def\figv@ctCreg#1(#2,#3){{\itis@vect@rtrue\figp@intreg#1:(#2,#3)}\ignorespaces}
\def\figvectDBezierDD#1:#2,#3[#4,#5,#6,#7]{\ifps@cri{\s@uvc@ntr@l\et@tfigvectDBezierDD%
    \figvectDBezier@#2,#3[#4,#5,#6,#7]\v@lX=\c@ef\v@lX\v@lY=\c@ef\v@lY%
    \figv@ctCreg#1(\v@lX,\v@lY)\resetc@ntr@l\et@tfigvectDBezierDD}\ignorespaces\fi}
\def\figvectDBezierTD#1:#2,#3[#4,#5,#6,#7]{\ifps@cri{\s@uvc@ntr@l\et@tfigvectDBezierTD%
    \figvectDBezier@#2,#3[#4,#5,#6,#7]\v@lX=\c@ef\v@lX\v@lY=\c@ef\v@lY\v@lZ=\c@ef\v@lZ%
    \figv@ctCreg#1(\v@lX,\v@lY,\v@lZ)\resetc@ntr@l\et@tfigvectDBezierTD}\ignorespaces\fi}
\def\figvectDBezier@#1,#2[#3,#4,#5,#6]{\setc@ntr@l{2}%
    \edef\T@{#2}\v@leur=\p@\advance\v@leur-#2pt\edef\UNmT@{\repdecn@mb{\v@leur}}%
    \ifnum#1=\tw@\def\c@ef{6}\else\def\c@ef{3}\fi%
    \figptcopy-4:/#3/\figptcopy-3:/#4/\figptcopy-2:/#5/\figptcopy-1:/#6/%
    \c@lDecast{-4}{-3}{-1}\ifnum#1=\tw@\c@lDCDeux{-4}{-3}\c@lDCDeux{-3}{-2}%
    \c@lDCDeux{-4}{-3}\else\c@lDecast{-4}{-3}{-2}\c@lDCDeux{-4}{-3}\fi\figg@tXY{-4}}
\def\c@lDCDeuxDD#1#2{\figg@tXY{#2}\figg@tXYa{#1}%
    \advance\v@lX-\v@lXa\advance\v@lY-\v@lYa\figp@intregDD#1:(\v@lX,\v@lY)}
\def\c@lDCDeuxTD#1#2{\figg@tXY{#2}\figg@tXYa{#1}\advance\v@lX-\v@lXa%
    \advance\v@lY-\v@lYa\advance\v@lZ-\v@lZa\figp@intregTD#1:(\v@lX,\v@lY,\v@lZ)}
\def\figvectNDD#1[#2,#3]{\ifps@cri{\figg@tXYa{#2}\figg@tXY{#3}%
    \advance\v@lX-\v@lXa\advance\v@lY-\v@lYa%
    \figv@ctCreg#1(-\v@lY,\v@lX)}\ignorespaces\fi}
\def\figvectNTD#1[#2,#3,#4]{\ifps@cri{\vecunitC@TD[#2,#4]\v@lmin=\v@lX\v@lmax=\v@lY%
    \v@leur=\v@lZ\vecunitC@TD[#2,#3]\c@lprovec{#1}}\ignorespaces\fi}
\def\figvectNVDD#1[#2]{\ifps@cri{\figg@tXY{#2}\figv@ctCreg#1(-\v@lY,\v@lX)}\ignorespaces\fi}
\def\figvectNVTD#1[#2,#3]{\ifps@cri{\vecunitCV@TD{#3}\v@lmin=\v@lX\v@lmax=\v@lY%
    \v@leur=\v@lZ\vecunitCV@TD{#2}\c@lprovec{#1}}\ignorespaces\fi}
\def\figvectPDD#1[#2,#3]{\ifps@cri{\figg@tXYa{#2}\figg@tXY{#3}%
    \advance\v@lX-\v@lXa\advance\v@lY-\v@lYa%
    \figv@ctCreg#1(\v@lX,\v@lY)}\ignorespaces\fi}
\def\figvectPTD#1[#2,#3]{\ifps@cri{\figg@tXYa{#2}\figg@tXY{#3}%
    \advance\v@lX-\v@lXa\advance\v@lY-\v@lYa\advance\v@lZ-\v@lZa%
    \figv@ctCreg#1(\v@lX,\v@lY,\v@lZ)}\ignorespaces\fi}
\def\figvectUDD#1[#2]{\ifps@cri{\n@rmeuc\v@leur{#2}\invers@\v@leur\v@leur%
    \delt@=\repdecn@mb{\v@leur}\unit@\edef\v@ldelt@{\repdecn@mb{\delt@}}%
    \figg@tXY{#2}\v@lX=\v@ldelt@\v@lX\v@lY=\v@ldelt@\v@lY%
    \figv@ctCreg#1(\v@lX,\v@lY)}\ignorespaces\fi}
\def\figvectUTD#1[#2]{\ifps@cri{\n@rmeuc\v@leur{#2}\invers@\v@leur\v@leur%
    \delt@=\repdecn@mb{\v@leur}\unit@\edef\v@ldelt@{\repdecn@mb{\delt@}}%
    \figg@tXY{#2}\v@lX=\v@ldelt@\v@lX\v@lY=\v@ldelt@\v@lY\v@lZ=\v@ldelt@\v@lZ%
    \figv@ctCreg#1(\v@lX,\v@lY,\v@lZ)}\ignorespaces\fi}
\def\figvisu#1#2#3{\c@ldefproj\initb@undb@x\setbox\b@xvisu=\hbox{\ignorespaces#3}%
    \v@lXa=-\c@@rdYmin\v@lYa=\c@@rdYmax\advance\v@lYa-\c@@rdYmin%
    \v@lX=\c@@rdXmax\advance\v@lX-\c@@rdXmin%
    \setbox#1=\hbox{#2}\v@lY=-\v@lX\maxim@m{\v@lX}{\v@lX}{\wd#1}%
    \advance\v@lY\v@lX\v@lY=0.5\v@lY\advance\v@lY-\c@@rdXmin%
    \setbox#1=\vbox{\parindent0mm\hsize=\v@lX\vskip\v@lYa%
    \rlap{\hskip\v@lY\smash{\raise\v@lXa\box\b@xvisu}}%
    \def\t@xt@{#2}\ifx\t@xt@\empty\else\medskip\centerline{#2}\fi}\wd#1=\v@lX}
\newdimen\Xt@rget  \newdimen\Yt@rget  \newdimen\Zt@rget
\newdimen\XminTD@ \newdimen\YminTD@ \newdimen\ZminTD@ 
\newdimen\XmaxTD@ \newdimen\YmaxTD@ \newdimen\ZmaxTD@
\newif\ifnewt@rgetpt\newif\ifnewdis@b
\def\b@undb@xTD#1#2#3{%
    \minim@m{\global\XminTD@}{\XminTD@}{#1}%
    \minim@m{\global\YminTD@}{\YminTD@}{#2}%
    \minim@m{\global\ZminTD@}{\ZminTD@}{#3}%
    \maxim@m{\global\XmaxTD@}{\XmaxTD@}{#1}%
    \maxim@m{\global\YmaxTD@}{\YmaxTD@}{#2}%
    \maxim@m{\global\ZmaxTD@}{\ZmaxTD@}{#3}}
\def\c@ldefdisob{{\ifdim\XminTD@<\maxdimen\v@leur=\XmaxTD@\advance\v@leur-\XminTD@%
    \delt@=\YmaxTD@\advance\delt@-\YminTD@\maxim@m{\v@leur}{\v@leur}{\delt@}%
    \delt@=\ZmaxTD@\advance\delt@-\ZminTD@\maxim@m{\v@leur}{\v@leur}{\delt@}%
    \v@leur=5\v@leur\else\v@leur=800pt\fi\c@ldefdisob@{\v@leur}}}
\def\c@ldefdisob@#1{{\v@leur=#1\ifdim\v@leur<\p@\v@leur=800pt\fi%
    \delt@=\ptT@unit@\v@leur\xdef\disob@unit{\repdecn@mb{\delt@}}%
    \f@ctech=\@ne\loop\ifdim\v@leur>10pt\divide\v@leur\t@n\multiply\f@ctech\t@n\repeat%
    \xdef\disob@{\repdecn@mb{\v@leur}}\xdef\divf@ctproj{\the\f@ctech}}%
    \global\newdis@btrue}
\def\c@ldeft@rgetpt{\global\newt@rgetpttrue%
    \v@leur=\XmaxTD@\advance\v@leur-\XminTD@\Xt@rget=\XminTD@\advance\Xt@rget.5\v@leur%
    \v@leur=\YmaxTD@\advance\v@leur-\YminTD@\Yt@rget=\YminTD@\advance\Yt@rget.5\v@leur%
    \v@leur=\ZmaxTD@\advance\v@leur-\ZminTD@\Zt@rget=\ZminTD@\advance\Zt@rget.5\v@leur}
\def\c@ldefprojTD{\ifnewt@rgetpt\else\c@ldeft@rgetpt\fi\ifnewdis@b\else\c@ldefdisob\fi}
\def\c@lprojcav{
    \v@leur=\cxa@\v@lY\advance\v@lX\v@leur%
    \v@leur=\cxb@\v@lY\v@lY=\v@lZ\advance\v@lY\v@leur\ignorespaces}
\def\c@lprojrea{
    \advance\v@lX-\Xt@rget\advance\v@lY-\Yt@rget\advance\v@lZ-\Zt@rget%
    \v@leur=\cza@\v@lX\advance\v@leur\czb@\v@lY\advance\v@leur\czc@\v@lZ%
    \divide\v@leur\divf@ctproj\advance\v@leur\disob@ pt\invers@{\v@leur}{\v@leur}%
    \v@leur=\disob@\v@leur\edef\v@lcoef{\repdecn@mb{\v@leur}}%
    \v@lXa=\cxa@\v@lX\advance\v@lXa\cxb@\v@lY\v@lXa=\v@lcoef\v@lXa%
    \v@lY=\cyb@\v@lY\advance\v@lY\cya@\v@lX\advance\v@lY\cyc@\v@lZ%
    \v@lY=\v@lcoef\v@lY\v@lX=\v@lXa\ignorespaces}
\def\c@lprojort{
    \v@lXa=\cxa@\v@lX\advance\v@lXa\cxb@\v@lY%
    \v@lY=\cyb@\v@lY\advance\v@lY\cya@\v@lX\advance\v@lY\cyc@\v@lZ%
    \v@lX=\v@lXa\ignorespaces}
\def\figptpr@j#1:#2/#3/{{\figg@tXY{#3}\superc@lprojSP%
    \figp@intregDD#1:{#2}(\v@lX,\v@lY)}\ignorespaces}
\def\figsetobdistTD(#1){{\ifcurr@ntPS%
    \immediate\write16{*** \BS@ figsetobdist is ignored inside a
     \BS@ psbeginfig-\BS@ psendfig block.}%
    \else\v@leur=#1\unit@\c@ldefdisob@{\v@leur}\fi}\ignorespaces}
\def\figs@tproj#1{%
    \if#13 \d@faultproj\else\if#1c\d@faultproj%
    \else\if#1o\xdef\curr@ntproj{1}\xdef\typ@proj{orthogonal}\figsetviewTD(40,25)%
         \global\let\c@lprojSP=\c@lprojort\global\let\superc@lprojSP=\c@lprojort%
    \else\if#1r\xdef\curr@ntproj{2}\xdef\typ@proj{realistic}\figsetviewTD(40,25)%
         \global\let\c@lprojSP=\c@lprojrea\global\let\superc@lprojSP=\c@lprojrea%
    \else\d@faultproj\message{*** Unknown projection. Cavalier projection assumed.}%
    \fi\fi\fi\fi}
\def\d@faultproj{\xdef\curr@ntproj{0}\xdef\typ@proj{cavalier}\figsetviewTD(40,0.5)%
         \global\let\c@lprojSP=\c@lprojcav\global\let\superc@lprojSP=\c@lprojcav}
\def\figsettargetTD[#1]{{\ifcurr@ntPS%
    \immediate\write16{*** \BS@ figsettarget is ignored inside a
     \BS@ psbeginfig-\BS@ psendfig block.}%
    \else\global\newt@rgetpttrue\figg@tXY{#1}\global\Xt@rget=\v@lX%
    \global\Yt@rget=\v@lY\global\Zt@rget=\v@lZ\fi}\ignorespaces}
\def\figsetviewTD(#1){\ifcurr@ntPS%
     \immediate\write16{*** \BS@ figsetview is ignored inside a
     \BS@ psbeginfig-\BS@ psendfig block.}\else\figsetview@#1,:\fi\ignorespaces}
\def\figsetview@#1,#2:{{\c@ssin{\v@lmin}{\v@lmax}{#1pt}\xdef\v@lPsi{#1}\def\t@xt@{#2}%
    \ifx\t@xt@\empty\def\@rgdeux{\v@lTheta}\else\def\Xarg@##1,{\def\@rgdeux{##1}}\Xarg@#2\fi%
    \xdef\costhet@{\repdecn@mb{\v@lmin}}\xdef\sinthet@{\repdecn@mb{\v@lmax}}%
    \ifcase\curr@ntproj%
    \v@leur=\@rgdeux\v@lmin\xdef\cxa@{\repdecn@mb{\v@leur}}%
    \v@leur=\@rgdeux\v@lmax\xdef\cxb@{\repdecn@mb{\v@leur}}\v@leur=\@rgdeux pt%
    \relax\ifdim\v@leur>\p@\message{*** Lambda too large ! See \BS@ figsetview !}\fi%
    \else%
    \v@lmax=-\v@lmax\xdef\cxa@{\repdecn@mb{\v@lmax}}\xdef\cxb@{\costhet@}%
    \ifx\t@xt@\empty\def\@rgdeux{25}\fi\c@ssin{\v@lmin}{\v@lmax}{\@rgdeux pt}%
    \v@lmax=-\v@lmax%
    \v@leur=\v@lmax\v@leur=\costhet@\v@leur\xdef\cya@{\repdecn@mb{\v@leur}}%
    \v@leur=\v@lmax\v@leur=\sinthet@\v@leur\xdef\cyb@{\repdecn@mb{\v@leur}}%
    \xdef\cyc@{\repdecn@mb{\v@lmin}}%
    \v@lmin=-\v@lmin%
    \v@leur=\v@lmin\v@leur=\costhet@\v@leur\xdef\cza@{\repdecn@mb{\v@leur}}%
    \v@leur=\v@lmin\v@leur=\sinthet@\v@leur\xdef\czb@{\repdecn@mb{\v@leur}}%
    \xdef\czc@{\repdecn@mb{\v@lmax}}\fi%
    \xdef\v@lTheta{\@rgdeux}}}
\def\initb@undb@xTD{\XminTD@=\maxdimen\YminTD@=\maxdimen\ZminTD@=\maxdimen%
    \XmaxTD@=-\maxdimen\YmaxTD@=-\maxdimen\ZmaxTD@=-\maxdimen}
\newbox\Gb@x      
\newbox\Gb@xSC    
\newtoks\c@nsymb  
\newdimen\V@l@ur  
\newif\ifr@undcoord\newif\ifunitpr@sent
\def\unssqrttw@{0.70710678}
\def\Ast{\raise-5.05pt\hbox{$\ast$}}
\def\Bullet{\raise-5.05pt\hbox{$\bullet$}}
\def\Circ{\raise-5.05pt\hbox{$\circ$}}
\def\Diamond{\raise-5.05pt\hbox{$\diamond$}}%
\def\c@nterpt{\ignorespaces%
    \kern-.5\wd\Gb@xSC%
    \raise-.5\ht\Gb@xSC\rlap{\hbox{\raise.5\dp\Gb@xSC\hbox{\copy\Gb@xSC}}}%
    \kern .5\wd\Gb@xSC\ignorespaces}
\def\b@undb@xSC#1#2{{\v@lXa=#1\v@lYa=#2%
    \V@l@ur=\ht\Gb@xSC\advance\V@l@ur\dp\Gb@xSC%
    \advance\v@lXa-.5\wd\Gb@xSC\advance\v@lYa-.5\V@l@ur\b@undb@x{\v@lXa}{\v@lYa}%
    \advance\v@lXa\wd\Gb@xSC\advance\v@lYa\V@l@ur\b@undb@x{\v@lXa}{\v@lYa}}}
\def\@keldist#1#2{\edef\Dist@n{#2}\y@tiunit{\Dist@n}%
    \ifunitpr@sent#1=\Dist@n\else#1=\Dist@n\unit@\fi}
\def\y@tiunit#1{\unitpr@sentfalse\expandafter\y@tiunit@#1:}
\def\y@tiunit@#1#2:{\ifcat#1a\unitpr@senttrue\else\def\l@suite{#2}%
    \ifx\l@suite\empty\else\y@tiunit@#2:\fi\fi}
\def\figcoordDD#1{{\v@lX=\ptT@unit@\v@lX\v@lY=\ptT@unit@\v@lY%
    \ifr@undcoord\ifcase#1\v@leur=0.5pt\or\v@leur=0.05pt\or\v@leur=0.005pt%
    \or\v@leur=0.0005pt\else\v@leur=\z@\fi%
    \ifdim\v@lX<\z@\advance\v@lX-\v@leur\else\advance\v@lX\v@leur\fi%
    \ifdim\v@lY<\z@\advance\v@lY-\v@leur\else\advance\v@lY\v@leur\fi\fi%
    (\@ffichnb{#1}{\repdecn@mb{\v@lX}},\ifmmode\else\thinspace\fi%
    \@ffichnb{#1}{\repdecn@mb{\v@lY}})}}
\def\@ffichnb#1#2{\def\@@ffich{\@ffich#1(}\edef\n@mbre{#2}%
    \expandafter\@@ffich\n@mbre)}
\def\@ffich#1(#2.#3){{#2\ifnum#1>\z@.\fi\def\dig@ts{#3}\s@mme=\z@%
    \loop\ifnum\s@mme<#1\expandafter\@ffichdec\dig@ts:\advance\s@mme\@ne\repeat}}
\def\@ffichdec#1#2:{\relax#1\def\dig@ts{#20}}
\def\figcoordTD#1{{\v@lX=\ptT@unit@\v@lX\v@lY=\ptT@unit@\v@lY\v@lZ=\ptT@unit@\v@lZ%
    \ifr@undcoord\ifcase#1\v@leur=0.5pt\or\v@leur=0.05pt\or\v@leur=0.005pt%
    \or\v@leur=0.0005pt\else\v@leur=\z@\fi%
    \ifdim\v@lX<\z@\advance\v@lX-\v@leur\else\advance\v@lX\v@leur\fi%
    \ifdim\v@lY<\z@\advance\v@lY-\v@leur\else\advance\v@lY\v@leur\fi%
    \ifdim\v@lZ<\z@\advance\v@lZ-\v@leur\else\advance\v@lZ\v@leur\fi\fi%
    (\@ffichnb{#1}{\repdecn@mb{\v@lX}},\ifmmode\else\thinspace\fi%
     \@ffichnb{#1}{\repdecn@mb{\v@lY}},\ifmmode\else\thinspace\fi%
     \@ffichnb{#1}{\repdecn@mb{\v@lZ}})}}
\def\figsetroundcoord#1{\figsetr@undcoord#1:}
\def\figsetr@undcoord#1#2:{\if#1n\r@undcoordfalse\else\r@undcoordtrue\fi}
\def\figsetmark#1{\c@nsymb={#1}\setbox\Gb@xSC=\hbox{\the\c@nsymb}\ignorespaces}
\def\figsetptname#1{\def\ptn@me##1{#1}\ignorespaces}
\def\figwrite[#1]#2{{\ignorespaces\def\list@num{#1}\@ecfor\p@int:=\list@num\do{%
    \setbox\Gb@x=\hbox{\def\t@xt@{#2}\ifx\t@xt@\empty\figg@tT{\p@int}\else#2\fi}%
    \figwrit@{\p@int}}}\ignorespaces}
\def\figwrit@#1{\figg@tXY{#1}\c@lprojSP%
    \rlap{\kern\v@lX\raise\v@lY\hbox{\unhcopy\Gb@x}}\V@l@ur=\v@lY%
    \advance\v@lY\ht\Gb@x\b@undb@x{\v@lX}{\v@lY}\advance\v@lX\wd\Gb@x%
    \v@lY=\V@l@ur\advance\v@lY-\dp\Gb@x\b@undb@x{\v@lX}{\v@lY}}
\def\figwritec[#1]#2{{\ignorespaces\def\list@num{#1}%
    \@ecfor\p@int:=\list@num\do{\figwrit@c{\p@int}{#2}}}\ignorespaces}
\def\figwrit@c#1#2{\figg@tXY{#1}\c@lprojSP%
    \setbox\Gb@x=\hbox{\def\t@xt@{#2}\ifx\t@xt@\empty\figg@tT{#1}\else#2\fi}%
    \rlap{\kern\v@lX\raise\v@lY\hbox{\rlap{\kern-.5\wd\Gb@x%
    \raise-.5\ht\Gb@x\hbox{\raise.5\dp\Gb@x\hbox{\unhcopy\Gb@x}}}}}%
    \V@l@ur=\ht\Gb@x\advance\V@l@ur\dp\Gb@x%
    \advance\v@lX-.5\wd\Gb@x\advance\v@lY-.5\V@l@ur\b@undb@x{\v@lX}{\v@lY}%
    \advance\v@lX\wd\Gb@x\advance\v@lY\V@l@ur\b@undb@x{\v@lX}{\v@lY}}
\def\figwritep[#1]{{\ignorespaces\def\list@num{#1}\setbox\Gb@x=\hbox{\c@nterpt}%
    \@ecfor\p@int:=\list@num\do{\figwrit@{\p@int}}}\ignorespaces}
\def\figwritew#1:#2(#3){{\ignorespaces\@keldist\V@l@ur{#3}\delt@=\z@\def\list@num{#1}%
    \@ecfor\p@int:=\list@num\do{\figwrit@gcw{\p@int}{#2}}}\ignorespaces}
\def\figwritee#1:#2(#3){{\ignorespaces\@keldist\V@l@ur{#3}\delt@=\z@\def\list@num{#1}%
    \@ecfor\p@int:=\list@num\do{\figwrit@gce{\p@int}{#2}}}\ignorespaces}
\def\figwriten#1:#2(#3){{\ignorespaces\@keldist\V@l@ur{#3}\def\list@num{#1}%
    \@ecfor\p@int:=\list@num\do{\figwrit@NS\p@int:{#2}(\@ne)}}\ignorespaces}
\def\figwrites#1:#2(#3){{\ignorespaces\@keldist\V@l@ur{#3}\V@l@ur=-\V@l@ur\def\list@num{#1}%
    \@ecfor\p@int:=\list@num\do{\figwrit@NS\p@int:{#2}(\m@ne)}}\ignorespaces}
\def\figwritenw#1:#2(#3){{\ignorespaces\@keldist\V@l@ur{#3}%
    \V@l@ur=\unssqrttw@\V@l@ur\delt@=\V@l@ur\ifdim\delt@=\z@\delt@=\epsil@n\fi%
    \def\list@num{#1}\@ecfor\p@int:=\list@num\do{\figwrit@gw{\p@int}{#2}}}\ignorespaces}
\def\figwritesw#1:#2(#3){{\ignorespaces\@keldist\V@l@ur{#3}%
    \V@l@ur=\unssqrttw@\V@l@ur\delt@=-\V@l@ur\ifdim\delt@=\z@\delt@=-\epsil@n\fi%
    \def\list@num{#1}\@ecfor\p@int:=\list@num\do{\figwrit@gw{\p@int}{#2}}}\ignorespaces}
\def\figwritene#1:#2(#3){{\ignorespaces\@keldist\V@l@ur{#3}%
    \V@l@ur=\unssqrttw@\V@l@ur\delt@=\V@l@ur\ifdim\delt@=\z@\delt@=\epsil@n\fi%
    \def\list@num{#1}\@ecfor\p@int:=\list@num\do{\figwrit@ge{\p@int}{#2}}}\ignorespaces}
\def\figwritese#1:#2(#3){{\ignorespaces\@keldist\V@l@ur{#3}%
    \V@l@ur=\unssqrttw@\V@l@ur\delt@=-\V@l@ur\ifdim\delt@=\z@\delt@=-\epsil@n\fi%
    \def\list@num{#1}\@ecfor\p@int:=\list@num\do{\figwrit@ge{\p@int}{#2}}}\ignorespaces}
\def\figwritegcw#1:#2(#3,#4){{\ignorespaces\@keldist\V@l@ur{#3}\@keldist\delt@{#4}%
    \def\list@num{#1}\@ecfor\p@int:=\list@num\do{\figwrit@gcw{\p@int}{#2}}}\ignorespaces}
\def\figwrit@gcw#1#2{\figg@tXY{#1}\c@lprojSP%
    \setbox\Gb@x=\hbox{\def\t@xt@{#2}\ifx\t@xt@\empty\figg@tT{#1}\else#2\fi}%
    \rlap{\kern\v@lX\raise\v@lY\hbox{\rlap{\kern-\wd\Gb@x\kern-\V@l@ur%
          \raise-.5\ht\Gb@x\hbox{\raise\delt@\hbox{\raise.5\dp\Gb@x\hbox{\unhcopy\Gb@x}}}}%
          \c@nterpt}}%
    \ifdim\wd\Gb@xSC>\z@\b@undb@xSC{\v@lX}{\v@lY}\fi
    \advance\v@lX-\wd\Gb@x\advance\v@lX-\V@l@ur\v@lZ=\ht\Gb@x\advance\v@lZ\dp\Gb@x%
    \advance\v@lY\delt@\advance\v@lY.5\v@lZ\b@undb@x{\v@lX}{\v@lY}%
    \advance\v@lX\wd\Gb@x\advance\v@lY-\v@lZ\b@undb@x{\v@lX}{\v@lY}}
\def\figwritegce#1:#2(#3,#4){{\ignorespaces\@keldist\V@l@ur{#3}\@keldist\delt@{#4}%
    \def\list@num{#1}\@ecfor\p@int:=\list@num\do{\figwrit@gce{\p@int}{#2}}}\ignorespaces}
\def\figwrit@gce#1#2{\figg@tXY{#1}\c@lprojSP%
    \setbox\Gb@x=\hbox{\def\t@xt@{#2}\ifx\t@xt@\empty\figg@tT{#1}\else#2\fi}%
    \rlap{\kern\v@lX\raise\v@lY\hbox{\c@nterpt\kern\V@l@ur%
          \raise-.5\ht\Gb@x\hbox{\raise\delt@\hbox{\raise.5\dp\Gb@x\hbox{\unhcopy\Gb@x}}}}}%
    \ifdim\wd\Gb@xSC>\z@\b@undb@xSC{\v@lX}{\v@lY}\fi
    \advance\v@lX\wd\Gb@x\advance\v@lX\V@l@ur\v@lZ=\ht\Gb@x\advance\v@lZ\dp\Gb@x%
    \advance\v@lY\delt@\advance\v@lY.5\v@lZ\b@undb@x{\v@lX}{\v@lY}%
    \advance\v@lX-\wd\Gb@x\advance\v@lY-\v@lZ\b@undb@x{\v@lX}{\v@lY}}
\def\figwritegw#1:#2(#3,#4){{\ignorespaces\@keldist\V@l@ur{#3}\@keldist\delt@{#4}%
    \def\list@num{#1}\@ecfor\p@int:=\list@num\do{\figwrit@gw{\p@int}{#2}}}\ignorespaces}
\def\figwrit@gw#1#2{\figg@tXY{#1}\c@lprojSP%
    \setbox\Gb@x=\hbox{\def\t@xt@{#2}\ifx\t@xt@\empty\figg@tT{#1}\else#2\fi}%
    \v@lXa=\z@\v@lYa=\ht\Gb@x\advance\v@lYa\dp\Gb@x%
    \ifdim\delt@>\z@\relax%
    \rlap{\kern\v@lX\raise\v@lY\hbox{\rlap{\kern-\wd\Gb@x\kern-\V@l@ur%
          \raise\delt@\hbox{\raise\dp\Gb@x\hbox{\unhcopy\Gb@x}}}\c@nterpt}}%
    \else\ifdim\delt@<\z@\relax\v@lYa=-\v@lYa%
    \rlap{\kern\v@lX\raise\v@lY\hbox{\rlap{\kern-\wd\Gb@x\kern-\V@l@ur%
          \raise\delt@\hbox{\raise-\ht\Gb@x\hbox{\unhcopy\Gb@x}}}\c@nterpt}}%
    \else\v@lXa=-.5\v@lYa%
    \rlap{\kern\v@lX\raise\v@lY\hbox{\rlap{\kern-\wd\Gb@x\kern-\V@l@ur%
          \raise-.5\ht\Gb@x\hbox{\raise.5\dp\Gb@x\hbox{\unhcopy\Gb@x}}}\c@nterpt}}%
    \fi\fi%
    \ifdim\wd\Gb@xSC>\z@\b@undb@xSC{\v@lX}{\v@lY}\fi
    \advance\v@lY\delt@%
    \advance\v@lX-\V@l@ur\advance\v@lY\v@lXa\b@undb@x{\v@lX}{\v@lY}%
    \advance\v@lX-\wd\Gb@x\advance\v@lY\v@lYa\b@undb@x{\v@lX}{\v@lY}}
\def\figwritege#1:#2(#3,#4){{\ignorespaces\@keldist\V@l@ur{#3}\@keldist\delt@{#4}%
    \def\list@num{#1}\@ecfor\p@int:=\list@num\do{\figwrit@ge{\p@int}{#2}}}\ignorespaces}
\def\figwrit@ge#1#2{\figg@tXY{#1}\c@lprojSP%
    \setbox\Gb@x=\hbox{\def\t@xt@{#2}\ifx\t@xt@\empty\figg@tT{#1}\else#2\fi}%
    \v@lXa=\z@\v@lYa=\ht\Gb@x\advance\v@lYa\dp\Gb@x%
    \ifdim\delt@>\z@\relax%
    \rlap{\kern\v@lX\raise\v@lY\hbox{\c@nterpt\kern\V@l@ur%
          \raise\delt@\hbox{\raise\dp\Gb@x\hbox{\unhcopy\Gb@x}}}}%
    \else\ifdim\delt@<\z@\relax\v@lYa=-\v@lYa%
    \rlap{\kern\v@lX\raise\v@lY\hbox{\c@nterpt\kern\V@l@ur%
          \raise\delt@\hbox{\raise-\ht\Gb@x\hbox{\unhcopy\Gb@x}}}}%
    \else\v@lXa=-.5\v@lYa%
    \rlap{\kern\v@lX\raise\v@lY\hbox{\c@nterpt\kern\V@l@ur%
          \raise-.5\ht\Gb@x\hbox{\raise.5\dp\Gb@x\hbox{\unhcopy\Gb@x}}}}%
    \fi\fi%
    \ifdim\wd\Gb@xSC>\z@\b@undb@xSC{\v@lX}{\v@lY}\fi
    \advance\v@lY\delt@%
    \advance\v@lX\V@l@ur\advance\v@lY\v@lXa\b@undb@x{\v@lX}{\v@lY}%
    \advance\v@lX\wd\Gb@x\advance\v@lY\v@lYa\b@undb@x{\v@lX}{\v@lY}}
\def\figwrit@NS#1:#2(#3){\figg@tXY{#1}\c@lprojSP%
    \setbox\Gb@x=\hbox{\def\t@xt@{#2}\ifx\t@xt@\empty\figg@tT{#1}\else#2\fi}%
    \v@lZ=\V@l@ur\relax\ifdim#3pt<\z@\advance\v@lZ-\ht\Gb@x%
    \else\advance\v@lZ\dp\Gb@x\fi%
    \rlap{\kern\v@lX\raise\v@lY\hbox{\rlap{\kern-.5\wd\Gb@x%
          \raise\v@lZ\hbox{\unhcopy\Gb@x}}\c@nterpt}}%
    \ifdim\wd\Gb@xSC>\z@\b@undb@xSC{\v@lX}{\v@lY}\fi
    \advance\v@lY\v@lZ%
    \advance\v@lY-\dp\Gb@x\advance\v@lX-.5\wd\Gb@x\b@undb@x{\v@lX}{\v@lY}%
    \advance\v@lY\ht\Gb@x\advance\v@lY\dp\Gb@x%
    \advance\v@lX\wd\Gb@x\b@undb@x{\v@lX}{\v@lY}}
\newread\frf@g  \newwrite\fwf@g
\newif\ifcurr@ntPS
\newif\ifps@cri
\newif\ifUse@llipse
\newif\ifpsdebugmode \psdebugmodefalse 
\newif\ifmored@ta
\def\c@pypsfile#1#2{\def\blankline{\par}{\catcode`\%=12
    \loop\ifeof#2\mored@tafalse\else\mored@tatrue\immediate\read#2 to\tr@c
    \ifx\tr@c\blankline\else\immediate\write#1{\tr@c}\fi\fi\ifmored@ta\repeat}}
\let\@psffilein=\frf@g 
\newif\if@psffileok    
\newif\if@psfbbfound   
\newif\if@psfverbose   
\@psfverbosetrue
\def\@psfgetbb#1{\global\@psfbbfoundfalse%
\global\def\@psfllx{0}\global\def\@psflly{0}%
\global\def\@psfurx{30}\global\def\@psfury{30}%
\openin\@psffilein=#1
\ifeof\@psffilein\errmessage{I couldn't open #1, will ignore it}\else
   \edef\setcolonc@tcode{\catcode`\noexpand\:\the\catcode`\:\relax}%
   {\@psffileoktrue \chardef\other=12
    \def\do##1{\catcode`##1=\other}\dospecials \catcode`\ =10 \setcolonc@tcode
    \loop
       \read\@psffilein to \@psffileline
       \ifeof\@psffilein\@psffileokfalse\else
          \expandafter\@psfaux\@psffileline:. \\%
       \fi
   \if@psffileok\repeat
   \if@psfbbfound\else
    \if@psfverbose\message{No bounding box comment in #1; using defaults}\fi\fi
   }\closein\@psffilein\fi}%
{\catcode`\%=12 \global\let\@psfpercent=
\long\def\@psfaux#1#2:#3\\{\ifx#1\@psfpercent
   \def\testit{#2}\ifx\testit\@psfbblit
      \@psfgrab #3 . . . \\%
      \@psffileokfalse
      \global\@psfbbfoundtrue
   \fi\else\ifx#1\par\else\@psffileokfalse\fi\fi}%
\def\@psfempty{}%
\def\@psfgrab #1 #2 #3 #4 #5\\{%
\global\def\@psfllx{#1}\ifx\@psfllx\@psfempty
      \@psfgrab #2 #3 #4 #5 .\\\else
   \global\def\@psflly{#2}%
   \global\def\@psfurx{#3}\global\def\@psfury{#4}\fi}%
\newif\iffillm@de
\def\pssetfillmode#1{\setfillm@de#1:}
\def\setfillm@de#1#2:{\if#1n\fillm@defalse\else\fillm@detrue\fi}
\newif\ifpstestm@de
\def\pssetupdate#1{\ifcurr@ntPS\immediate\write16{*** \BS@ pssetupdate is ignored inside a
     \BS@ psbeginfig-\BS@ psendfig block.}%
    \immediate\write16{*** It must be called before \BS@ psbeginfig.}\else\setupd@te#1:\fi}
\def\setupd@te#1#2:{\if#1n\pstestm@defalse\else\pstestm@detrue\fi}
\def\psaltitude#1[#2,#3,#4]{{\ifps@cri%
    \psc@mment{psaltitude Square Dim=#1, Triangle=[#2 / #3,#4]}%
    \s@uvc@ntr@l\et@tpsaltitude\resetc@ntr@l{2}\figptorthoprojline-5:=#2/#3,#4/%
    \figvectP -1[#3,#4]\n@rminf{\v@leur}{-1}%
    \figvectP -1[-5,#3]\n@rminf{\v@lmin}{-1}\figvectP -2[-5,#4]\n@rminf{\v@lmax}{-2}%
    \ifdim\v@lmin<\v@lmax\vecunit@{-3}{-2}\s@mme=#3%
     \else\v@lmax=\v@lmin\vecunit@{-3}{-1}\s@mme=#4\fi%
    \figvectP -4[-5,#2]\vecunit@{-4}{-4}\delt@=#1\unit@%
    \edef\t@ille{\repdecn@mb{\delt@}}\figpttra-1:=-5/\t@ille,-4/%
    \figpttra-2:=-1/\t@ille,-3/\figpttra-3:=-2/-\t@ille,-4/%
    \psline[#2,-5]\psline[-1,-2,-3]%
    \ifdim\v@leur<\v@lmax\pss@tsecondSt\psline[-5,\the\s@mme]\psrest@reSt\fi%
    \psc@mment{End psaltitude}\resetc@ntr@l\et@tpsaltitude\fi}}
\def\ps@rcerc#1;#2(#3,#4){\ellBB@x#1;#2,#2(#3,#4,0)%
    \immediate\write\fwf@g{newpath}{\delt@=#2\unit@\delt@=\ptT@ptps\delt@%
    \BdingB@xfalse%
    \pswrit@cmd{#1}{\repdecn@mb{\delt@}\space #3\space #4\space arc}{\fwf@g}}}
\def\psarccircDD#1;#2(#3,#4){\ifps@cri%
    \psc@mment{psarccircDD Center=#1 ; Radius=#2 (Ang1=#3, Ang2=#4)}%
    \iffillm@de\immediate\write\fwf@g{currentrgbcolor}\ps@rcerc#1;#2(#3,#4)%
    \immediate\write\fwf@g{\curr@ntcolor\setfillc@md\space setrgbcolor}%
    \else\ps@rcerc#1;#2(#3,#4)\immediate\write\fwf@g{stroke}\fi%
    \psc@mment{End psarccircDD}\fi}
\def\psarccircTD#1,#2,#3;#4(#5,#6){{\ifps@cri\s@uvc@ntr@l\et@tpsarccircTD%
    \psc@mment{psarccircTD Center=#1,P1=#2,P2=#3 ; Radius=#4 (Ang1=#5, Ang2=#6)}%
    \setc@ntr@l{2}\c@lExtAxes#1,#2,#3(#4)\psarcellPATD#1,-4,-5(#5,#6)%
    \psc@mment{End psarccircTD}\resetc@ntr@l\et@tpsarccircTD\fi}}
\def\c@lExtAxes#1,#2,#3(#4){%
    \figvectPTD-5[#1,#2]\vecunit@{-5}{-5}\figvectNTD-4[#1,#2,#3]\vecunit@{-4}{-4}%
    \figvectNVTD-3[-4,-5]\delt@=#4\unit@\edef\r@yon{\repdecn@mb{\delt@}}%
    \figpttra-4:=#1/\r@yon,-5/\figpttra-5:=#1/\r@yon,-3/}
\def\psarccircPDD#1;#2[#3,#4]{{\ifps@cri\s@uvc@ntr@l\et@tpsarccircPDD%
    \psc@mment{psarccircPDD Center=#1; Radius=#2, [P1=#3, P2=#4]}%
    \ps@ngleparam#1;#2[#3,#4]\ifdim\v@lmin>\v@lmax\advance\v@lmax\PI@deg\fi%
    \edef\@ngdeb{\repdecn@mb{\v@lmin}}\edef\@ngfin{\repdecn@mb{\v@lmax}}%
    \psarccirc#1;\r@dius(\@ngdeb,\@ngfin)%
    \psc@mment{End psarccircPDD}\resetc@ntr@l\et@tpsarccircPDD\fi}}
\def\psarccircPTD#1;#2[#3,#4,#5]{{\ifps@cri\s@uvc@ntr@l\et@tpsarccircPTD%
    \psc@mment{psarccircPTD Center=#1; Radius=#2, [P1=#3, P2=#4, P3=#5]}%
    \setc@ntr@l{2}\c@lExtAxes#1,#3,#5(#2)\psarcellPP#1,-4,-5[#3,#4]%
    \psc@mment{End psarccircPTD}\resetc@ntr@l\et@tpsarccircPTD\fi}}
\def\ps@ngleparam#1;#2[#3,#4]{\setc@ntr@l{2}%
    \figvectPDD-1[#1,#3]\vecunit@{-1}{-1}\figg@tXY{-1}\arct@n\v@lmin(\v@lX,\v@lY)%
    \figvectPDD-2[#1,#4]\vecunit@{-2}{-2}\figg@tXY{-2}\arct@n\v@lmax(\v@lX,\v@lY)%
    \v@lmin=\rdT@deg\v@lmin\v@lmax=\rdT@deg\v@lmax%
    \v@leur=#2pt\maxim@m{\mili@u}{-\v@leur}{\v@leur}%
    \edef\r@dius{\repdecn@mb{\mili@u}}}
\def\ps@rell#1;#2,#3(#4,#5,#6){\ellBB@x#1;#2,#3(#4,#5,#6)%
    \immediate\write\fwf@g{newpath}{\v@lmin=#2\unit@\v@lmin=\ptT@ptps\v@lmin%
    \v@lmax=#3\unit@\v@lmax=\ptT@ptps\v@lmax\BdingB@xfalse%
    \pswrit@cmd{#1}%
    {#6\space\repdecn@mb{\v@lmin}\space\repdecn@mb{\v@lmax}\space #4\space #5\space ellipse}{\fwf@g}}%
    \global\Use@llipsetrue}
\def\psarcellDD#1;#2,#3(#4,#5,#6){{\ifps@cri%
    \psc@mment{psarcellDD Center=#1 ; XRad=#2, YRad=#3 (Ang1=#4, Ang2=#5, Inclination=#6)}%
    \iffillm@de\immediate\write\fwf@g{currentrgbcolor}\ps@rell#1;#2,#3(#4,#5,#6)%
    \immediate\write\fwf@g{\curr@ntcolor\setfillc@md\space setrgbcolor}%
    \else\ps@rell#1;#2,#3(#4,#5,#6)\immediate\write\fwf@g{stroke}\fi%
    \psc@mment{End psarcellDD}\fi}}
\def\psarcellTD#1;#2,#3(#4,#5,#6){{\ifps@cri\s@uvc@ntr@l\et@tpsarcellTD%
    \psc@mment{psarcellTD Center=#1 ; XRad=#2, YRad=#3 (Ang1=#4, Ang2=#5, Inclination=#6)}%
    \setc@ntr@l{2}\figpttraC -8:=#1/#2,0,0/\figpttraC -7:=#1/0,#3,0/%
    \figvectC -4(0,0,1)\figptsrot -8=-8,-7/#1,#6,-4/\psarcellPATD#1,-8,-7(#4,#5)%
    \psc@mment{End psarcellTD}\resetc@ntr@l\et@tpsarcellTD\fi}}
\def\psarcellPADD#1,#2,#3(#4,#5){{\ifps@cri\s@uvc@ntr@l\et@tpsarcellPADD%
    \psc@mment{psarcellPADD Center=#1,PtAxis1=#2,PtAxis2=#3 (Ang1=#4, Ang2=#5)}%
    \setc@ntr@l{2}\figvectPDD-1[#1,#2]\vecunit@DD{-1}{-1}\v@lX=\ptT@unit@\result@t%
    \edef\XR@d{\repdecn@mb{\v@lX}}\figg@tXY{-1}\arct@n\v@lmin(\v@lX,\v@lY)%
    \v@lmin=\rdT@deg\v@lmin\edef\Inclin@{\repdecn@mb{\v@lmin}}%
    \figgetdist\YR@d[#1,#3]\psarcellDD#1;\XR@d,\YR@d(#4,#5,\Inclin@)%
    \psc@mment{End psarcellPADD}\resetc@ntr@l\et@tpsarcellPADD\fi}}
\def\psarcellPATD#1,#2,#3(#4,#5){{\ifps@cri\s@uvc@ntr@l\et@tpsarcellPATD%
    \psc@mment{psarcellPATD Center=#1,PtAxis1=#2,PtAxis2=#3 (Ang1=#4, Ang2=#5)}%
    \iffillm@de\immediate\write\fwf@g{currentrgbcolor}\ps@rellPATD#1,#2,#3(#4,#5)%
    \immediate\write\fwf@g{\curr@ntcolor\setfillc@md\space setrgbcolor}%
    \else\ps@rellPATD#1,#2,#3(#4,#5)\immediate\write\fwf@g{stroke}\fi%
    \psc@mment{End psarcellPATD}\resetc@ntr@l\et@tpsarcellPATD\fi}}
\def\ps@rellPATD#1,#2,#3(#4,#5){\let\c@lprojSP=\relax%
    \setc@ntr@l{2}\figvectPTD-1[#1,#2]\figvectPTD-2[#1,#3]\c@lNbarcs{#4}{#5}%
    \v@leur=#4pt\c@lptellP{#1}{-1}{-2}\figptpr@j-5:/-3/%
    \immediate\write\fwf@g{newpath}\pswrit@cmdS{-5}{moveto}{\fwf@g}{\X@un}{\Y@un}%
    \s@mme=\z@\loop\ifnum\s@mme<\p@rtent\advance\s@mme\@ne%
    \figg@tXY{-5}\v@lX=8\v@lX\v@lY=8\v@lY\figp@intregDD-4:(\v@lX,\v@lY)%
    \advance\v@leur\delt@\c@lptellP{#1}{-1}{-2}\figptpr@j-6:/-3/\figpttraDD-4:=-4/-27,-6/%
    \advance\v@leur\delt@\c@lptellP{#1}{-1}{-2}\figptpr@j-6:/-3/\figpttraDD-5:=-5/-27,-6/%
    \advance\v@leur\delt@\c@lptellP{#1}{-1}{-2}\figptpr@j-6:/-3/%
    \figpttraDD-4:=-4/1,-6/\figpttraDD-5:=-5/8,-6/\BdingB@xfalse%
    \c@lptCtrlDD-3[-5,-4]\pswrit@cmdS{-3}{}{\fwf@g}{\X@de}{\Y@de}%
    \c@lptCtrlDD-3[-4,-5]\pswrit@cmdS{-3}{}{\fwf@g}{\X@tr}{\Y@tr}%
    \BdingB@xtrue\pswrit@cmdS{-6}{curveto}{\fwf@g}{\X@qu}{\Y@qu}%
    \B@zierBB@x{1}{\Y@un}(\X@un,\X@de,\X@tr,\X@qu)%
    \B@zierBB@x{2}{\X@un}(\Y@un,\Y@de,\Y@tr,\Y@qu)%
    \edef\X@un{\X@qu}\edef\Y@un{\Y@qu}\figptcopyDD-5:/-6/\repeat}
\def\c@lNbarcs#1#2{%
    \delt@=#2pt\advance\delt@-#1pt\maxim@m{\v@lmax}{\delt@}{-\delt@}%
    \v@leur=\v@lmax\divide\v@leur45 \p@rtentiere{\p@rtent}{\v@leur}\advance\p@rtent\@ne%
    \s@mme=\p@rtent\multiply\s@mme\thr@@\divide\delt@\s@mme}
\def\c@lptCtrlDD#1[#2,#3]{\figptcopyDD#1:/#2/\figpttraDD#1:=#1/-2,#3/%
    \figg@tXY{#1}\divide\v@lX18 \divide\v@lY18 \figp@intregDD#1:(\v@lX,\v@lY)}
\def\c@lptCtrlTD#1[#2,#3]{\figptcopyTD#1:/#2/\figpttraTD#1:=#1/-2,#3/%
    \figg@tXY{#1}\divide\v@lX18 \divide\v@lY18 \divide\v@lZ18 %
    \figp@intregTD#1:(\v@lX,\v@lY,\v@lZ)}
\def\psarcellPP#1,#2,#3[#4,#5]{{\ifps@cri\s@uvc@ntr@l\et@tpsarcellPP%
    \psc@mment{psarcellPP Center=#1,PtAxis1=#2,PtAxis2=#3 [Point1=#4, Point2=#5]}%
    \setc@ntr@l{2}\figvectP-2[#1,#3]\vecunit@{-2}{-2}\v@lmin=\result@t%
    \invers@{\v@lmax}{\v@lmin}%
    \figvectP-1[#1,#2]\vecunit@{-1}{-1}\v@leur=\result@t%
    \v@leur=\repdecn@mb{\v@lmax}\v@leur\edef\AsB@{\repdecn@mb{\v@leur}}
    \c@lAngle{#1}{#4}{\v@lmin}\edef\@ngdeb{\repdecn@mb{\v@lmin}}%
    \c@lAngle{#1}{#5}{\v@lmax}\ifdim\v@lmin>\v@lmax\advance\v@lmax\PI@deg\fi%
    \edef\@ngfin{\repdecn@mb{\v@lmax}}\psarcellPA#1,#2,#3(\@ngdeb,\@ngfin)%
    \psc@mment{End psarcellPP}\resetc@ntr@l\et@tpsarcellPP\fi}}
\def\c@lAngle#1#2#3{\figvectP-3[#1,#2]%
    \c@lproscal\delt@[-3,-1]\c@lproscal\v@leur[-3,-2]%
    \v@leur=\AsB@\v@leur\arct@n#3(\delt@,\v@leur)#3=\rdT@deg#3}
\newif\if@rrowratio\@rrowratiotrue
\newif\if@rrowhfill
\newif\if@rrowhout
\def\pssetarrowheadangle#1{\edef\@rrowheadangle{#1}{\c@ssin{\v@lmin}{\v@lmax}{#1pt}%
    \xdef\C@AHANG{\repdecn@mb{\v@lmin}}\xdef\S@AHANG{\repdecn@mb{\v@lmax}}%
    \invers@{\v@leur}{\v@lmax}\maxim@m{\v@leur}{\v@leur}{-\v@leur}%
    \xdef\UNSS@N{\the\v@leur}}}
\def\pssetarrowheadfill#1{\set@rrowhfill#1:}
\def\set@rrowhfill#1#2:{\if#1n\@rrowhfillfalse\else\@rrowhfilltrue\fi}
\def\pssetarrowheadout#1{\set@rrowhout#1:}
\def\set@rrowhout#1#2:{\if#1n\@rrowhoutfalse\else\@rrowhouttrue\fi}
\def\pssetarrowheadlength#1{\edef\@rrowheadlength{#1}\@rrowratiofalse}
\def\pssetarrowheadratio#1{\edef\@rrowheadratio{#1}\@rrowratiotrue}
\def\psresetarrowhead{%
    \pssetarrowheadangle{\defaultarrowheadangle}%
    \pssetarrowheadfill{no}\pssetarrowheadout{no}%
    \pssetarrowheadratio{\defaultarrowheadratio}%
    \c@ldef@hlength\pssetarrowheadlength{\defaultarrowheadlength}}
\def\defaultarrowheadratio{0.1}
\def\defaultarrowheadangle{20}
\def\c@ldef@hlength{{\v@leur=10pt\v@leur=\ptT@unit@\v@leur%
    \xdef\defaultarrowheadlength{\repdecn@mb{\v@leur}}}}
\def\psarrowDD[#1,#2]{{\ifps@cri\s@uvc@ntr@l\et@tpsarrow%
    \psc@mment{psarrowDD [Pt1,Pt2]=[#1,#2]}\pssetfillmode{no}%
    \psarrowheadDD[#1,#2]\setc@ntr@l{2}\psline[#1,-3]%
    \psc@mment{End psarrowDD}\resetc@ntr@l\et@tpsarrow\fi}}
\def\psarrowTD[#1,#2]{{\ifps@cri\s@uvc@ntr@l\et@tpsarrowTD%
    \psc@mment{psarrowTD [Pt1,Pt2]=[#1,#2]}\resetc@ntr@l{2}%
    \figptpr@j-5:/#1/\figptpr@j-6:/#2/\let\c@lprojSP=\relax\psarrowDD[-5,-6]%
    \psc@mment{End psarrowTD}\resetc@ntr@l\et@tpsarrowTD\fi}}
\def\psarrowheadDD[#1,#2]{{\ifps@cri\s@uvc@ntr@l\et@tpsarrowheadDD%
    \if@rrowhfill\def\@hangle{-\@rrowheadangle}\else\def\@hangle{\@rrowheadangle}\fi%
    \if@rrowratio%
    \if@rrowhout\def\@hratio{-\@rrowheadratio}\else\def\@hratio{\@rrowheadratio}\fi%
    \psc@mment{psarrowheadDD Ratio=\@hratio, Angle=\@hangle, [Pt1,Pt2]=[#1,#2]}%
    \ps@rrowhead\@hratio,\@hangle[#1,#2]%
    \else%
    \if@rrowhout\def\@hlength{-\@rrowheadlength}\else\def\@hlength{\@rrowheadlength}\fi%
    \psc@mment{psarrowheadDD Length=\@hlength, Angle=\@hangle, [Pt1,Pt2]=[#1,#2]}%
    \ps@rrowheadfd\@hlength,\@hangle[#1,#2]%
    \fi%
    \psc@mment{End psarrowheadDD}\resetc@ntr@l\et@tpsarrowheadDD\fi}}
\def\psarrowheadTD[#1,#2]{{\ifps@cri\s@uvc@ntr@l\et@tpsarrowheadTD%
    \psc@mment{psarrowheadTD [Pt1,Pt2]=[#1,#2]}\resetc@ntr@l{2}%
    \figptpr@j-5:/#1/\figptpr@j-6:/#2/\let\c@lprojSP=\relax\psarrowheadDD[-5,-6]%
    \psc@mment{End psarrowheadTD}\resetc@ntr@l\et@tpsarrowheadTD\fi}}
\def\ps@rrowhead#1,#2[#3,#4]{{%
    \psc@mment{ps@rrowhead Ratio=#1, Angle=#2, [Pt1,Pt2]=[#3,#4]}\v@leur=\UNSS@N%
    \v@leur=\curr@ntwidth\v@leur\v@leur=\ptpsT@pt\v@leur\delt@=.5\v@leur
    \setc@ntr@l{2}\figvectPDD-3[#4,#3]%
    \figg@tXY{-3}\v@lX=#1\v@lX\v@lY=#1\v@lY\figv@ctCreg-3(\v@lX,\v@lY)%
    \vecunit@{-4}{-3}\mili@u=\result@t%
    \ifdim#2pt>\z@\v@lXa=-\C@AHANG\delt@%
     \edef\c@ef{\repdecn@mb{\v@lXa}}\figpttraDD-3:=-3/\c@ef,-4/\fi%
    \edef\c@ef{\repdecn@mb{\delt@}}%
    \v@lXa=\mili@u\v@lXa=\C@AHANG\v@lXa%
    \v@lYa=\ptpsT@pt\p@\v@lYa=\curr@ntwidth\v@lYa\v@lYa=\sDcc@ngle\v@lYa%
    \advance\v@lXa-\v@lYa\gdef\sDcc@ngle{0}%
    \ifdim\v@lXa>\v@leur\edef\c@efendpt{\repdecn@mb{\v@leur}}%
    \else\edef\c@efendpt{\repdecn@mb{\v@lXa}}\fi%
    \figg@tXY{-3}\v@lmin=\v@lX\v@lmax=\v@lY%
    \v@lXa=\C@AHANG\v@lmin\v@lYa=\S@AHANG\v@lmax\advance\v@lXa\v@lYa%
    \v@lYa=-\S@AHANG\v@lmin\v@lX=\C@AHANG\v@lmax\advance\v@lYa\v@lX%
    \setc@ntr@l{1}\figg@tXY{#4}\advance\v@lX\v@lXa\advance\v@lY\v@lYa%
    \setc@ntr@l{2}\figp@intregDD-2:(\v@lX,\v@lY)%
    \v@lXa=\C@AHANG\v@lmin\v@lYa=-\S@AHANG\v@lmax\advance\v@lXa\v@lYa%
    \v@lYa=\S@AHANG\v@lmin\v@lX=\C@AHANG\v@lmax\advance\v@lYa\v@lX%
    \setc@ntr@l{1}\figg@tXY{#4}\advance\v@lX\v@lXa\advance\v@lY\v@lYa%
    \setc@ntr@l{2}\figp@intregDD-1:(\v@lX,\v@lY)%
    \ifdim#2pt<\z@\fillm@detrue\psline[-2,#4,-1]
    \else\figptstraDD-3=#4,-2,-1/\c@ef,-4/\psline[-2,-3,-1]\fi
    \ifdim#1pt>\z@\figpttraDD-3:=#4/\c@efendpt,-4/\else\figptcopyDD-3:/#4/\fi%
    \psc@mment{End ps@rrowhead}}}
\def\sDcc@ngle{0}
\def\ps@rrowheadfd#1,#2[#3,#4]{{%
    \psc@mment{ps@rrowheadfd Length=#1, Angle=#2, [Pt1,Pt2]=[#3,#4]}%
    \setc@ntr@l{2}\figvectPDD-1[#3,#4]\n@rmeucDD{\v@leur}{-1}\v@leur=\ptT@unit@\v@leur%
    \invers@{\v@leur}{\v@leur}\v@leur=#1\v@leur\edef\R@tio{\repdecn@mb{\v@leur}}%
    \ps@rrowhead\R@tio,#2[#3,#4]\psc@mment{End ps@rrowheadfd}}}
\def\psarrowBezierDD[#1,#2,#3,#4]{{\ifps@cri\s@uvc@ntr@l\et@tpsarrowBezierDD%
    \psc@mment{psarrowBezierDD Control points=#1,#2,#3,#4}\setc@ntr@l{2}%
    \if@rrowratio\c@larclengthDD\v@leur,10[#1,#2,#3,#4]\else\v@leur=\z@\fi%
    \ps@rrowB@zDD\v@leur[#1,#2,#3,#4]%
    \psc@mment{End psarrowBezierDD}\resetc@ntr@l\et@tpsarrowBezierDD\fi}}
\def\psarrowBezierTD[#1,#2,#3,#4]{{\ifps@cri\s@uvc@ntr@l\et@tpsarrowBezierTD%
    \psc@mment{psarrowBezierTD Control points=#1,#2,#3,#4}\resetc@ntr@l{2}%
    \figptpr@j-7:/#1/\figptpr@j-8:/#2/\figptpr@j-9:/#3/\figptpr@j-10:/#4/%
    \let\c@lprojSP=\relax\ifnum\curr@ntproj<\tw@\psarrowBezierDD[-7,-8,-9,-10]%
    \else\immediate\write\fwf@g{newpath}\pswrit@cmd{-7}{moveto}{\fwf@g}%
    \if@rrowratio\c@larclengthDD\mili@u,10[-7,-8,-9,-10]\else\mili@u=\z@\fi%
    \p@rtent=\NBz@rcs\advance\p@rtent\m@ne\subB@zierTD\p@rtent[#1,#2,#3,#4]%
    \immediate\write\fwf@g{stroke}%
    \advance\v@lmin\p@rtent\delt@
    \v@leur=\v@lmin\advance\v@leur0.33333\delt@\edef\unti@rs{\repdecn@mb{\v@leur}}%
    \v@leur=\v@lmin\advance\v@leur0.66666\delt@\edef\deti@rs{\repdecn@mb{\v@leur}}%
    \figptcopyDD-8:/-10/\c@lsubBzarc\unti@rs,\deti@rs[#1,#2,#3,#4]%
    \figptcopyDD-8:/-4/\figptcopyDD-9:/-3/\ps@rrowB@zDD\mili@u[-7,-8,-9,-10]\fi%
    \psc@mment{End psarrowBezierTD}\resetc@ntr@l\et@tpsarrowBezierTD\fi}}
\def\c@larclengthDD#1,#2[#3,#4,#5,#6]{{\p@rtent=#2\figptcopyDD-5:/#3/%
    \delt@=\p@\divide\delt@\p@rtent\c@rre=\z@\v@leur=\z@\s@mme=\z@%
    \loop\ifnum\s@mme<\p@rtent\advance\s@mme\@ne\advance\v@leur\delt@%
    \edef\T@{\repdecn@mb{\v@leur}}\figptBezierDD-6::\T@[#3,#4,#5,#6]%
    \figvectPDD-1[-5,-6]\n@rmeucDD{\mili@u}{-1}\advance\c@rre\mili@u%
    \figptcopyDD-5:/-6/\repeat\global\result@t=\ptT@unit@\c@rre}#1=\result@t}
\def\ps@rrowB@zDD#1[#2,#3,#4,#5]{{\pssetfillmode{no}%
    \if@rrowratio\delt@=\@rrowheadratio#1\else\delt@=\@rrowheadlength pt\fi%
    \v@leur=\C@AHANG\delt@\edef\R@dius{\repdecn@mb{\v@leur}}%
    \figptintercircB@zDD-5:\R@dius[#5,#4,#3,#2]%
    \pssetarrowheadlength{\repdecn@mb{\delt@}}\psarrowheadDD[-5,#5]%
    \let\n@rmeuc=\n@rmeucDD\figgetdist\R@dius[#5,-3]\figptBezierDD-5::0.33333[#5,#4,#3,#2]%
    \figptBezierDD-6::0.66666[#5,#4,#3,#2]\figptintercircB@zDD-3:\R@dius[#5,#4,#3,#2]%
    \figptscontrolDD-5[-3,-5,-6,#2]\psBezierDD1[-3,-5,-4,#2]}}
\def\psarrowcircDD#1;#2(#3,#4){{\ifps@cri\s@uvc@ntr@l\et@tpsarrowcircDD%
    \psc@mment{psarrowcircDD Center=#1 ; Radius=#2 (Ang1=#3,Ang2=#4)}%
    \pssetfillmode{no}\pscirc@rrowhead#1;#2(#3,#4)%
    \setc@ntr@l{2}\figvectPDD -4[#1,-3]\vecunit@{-4}{-4}%
    \figg@tXY{-4}\arct@n\v@lmin(\v@lX,\v@lY)%
    \v@lmin=\rdT@deg\v@lmin\v@leur=#4pt\advance\v@leur-\v@lmin%
    \maxim@m{\v@leur}{\v@leur}{-\v@leur}%
    \ifdim\v@leur>90pt\relax\ifdim\v@lmin<#4pt\advance\v@lmin\PI@deg%
    \else\advance\v@lmin-\PI@deg\fi\fi\edef\ar@ngle{\repdecn@mb{\v@lmin}}%
    \ifdim#3pt<#4pt\psarccirc#1;#2(#3,\ar@ngle)\else\psarccirc#1;#2(\ar@ngle,#3)\fi%
    \psc@mment{End psarrowcircDD}\resetc@ntr@l\et@tpsarrowcircDD\fi}}
\def\psarrowcircTD#1,#2,#3;#4(#5,#6){{\ifps@cri\s@uvc@ntr@l\et@tpsarrowcircTD%
    \psc@mment{psarrowcircTD Center=#1,P1=#2,P2=#3 ; Radius=#4 (Ang1=#5, Ang2=#6)}%
    \resetc@ntr@l{2}\c@lExtAxes#1,#2,#3(#4)\let\c@lprojSP=\relax%
    \figvectPTD-11[#1,-4]\figvectPTD-12[#1,-5]\c@lNbarcs{#5}{#6}%
    \if@rrowratio\v@lmax=\degT@rd\v@lmax\edef\D@lpha{\repdecn@mb{\v@lmax}}\fi%
    \advance\p@rtent\m@ne\mili@u=\z@%
    \v@leur=#5pt\c@lptellP{#1}{-11}{-12}\figptpr@j-9:/-3/%
    \immediate\write\fwf@g{newpath}\pswrit@cmdS{-9}{moveto}{\fwf@g}{\X@un}{\Y@un}%
    \s@mme=\z@\loop\ifnum\s@mme<\p@rtent\advance\s@mme\@ne%
    \advance\v@leur\delt@\c@lptellP{#1}{-11}{-12}\figptpr@j-5:/-3/%
    \advance\v@leur\delt@\c@lptellP{#1}{-11}{-12}\figptpr@j-6:/-3/%
    \advance\v@leur\delt@\c@lptellP{#1}{-11}{-12}\figptpr@j-10:/-3/%
    \figptscontrolDD-8[-9,-5,-6,-10]\BdingB@xfalse%
    \pswrit@cmdS{-8}{}{\fwf@g}{\X@de}{\Y@de}\pswrit@cmdS{-7}{}{\fwf@g}{\X@tr}{\Y@tr}%
    \BdingB@xtrue\pswrit@cmdS{-10}{curveto}{\fwf@g}{\X@qu}{\Y@qu}%
    \if@rrowratio\c@lcurvradDD0.5[-9,-8,-7,-10]\advance\mili@u\result@t\fi%
    \B@zierBB@x{1}{\Y@un}(\X@un,\X@de,\X@tr,\X@qu)%
    \B@zierBB@x{2}{\X@un}(\Y@un,\Y@de,\Y@tr,\Y@qu)%
    \edef\X@un{\X@qu}\edef\Y@un{\Y@qu}\figptcopyDD-9:/-10/\repeat%
    \immediate\write\fwf@g{stroke}%
    \advance\v@leur\delt@\c@lptellP{#1}{-11}{-12}\figptpr@j-5:/-3/%
    \advance\v@leur\delt@\c@lptellP{#1}{-11}{-12}\figptpr@j-6:/-3/%
    \advance\v@leur\delt@\c@lptellP{#1}{-11}{-12}\figptpr@j-10:/-3/%
    \figptscontrolDD-8[-9,-5,-6,-10]%
    \if@rrowratio\c@lcurvradDD0.5[-9,-8,-7,-10]\advance\mili@u\result@t%
    \maxim@m{\mili@u}{\mili@u}{-\mili@u}\mili@u=\ptT@unit@\mili@u%
    \mili@u=\D@lpha\mili@u\advance\p@rtent\@ne\divide\mili@u\p@rtent\fi%
    \ps@rrowB@zDD\mili@u[-9,-8,-7,-10]%
    \psc@mment{End psarrowcircTD}\resetc@ntr@l\et@tpsarrowcircTD\fi}}
\def\pscirc@rrowhead#1;#2(#3,#4){{%
    \psc@mment{pscirc@rrowhead Center=#1 ; Radius=#2 (Ang1=#3,Ang2=#4)}%
    \v@leur=#2\unit@\edef\s@glen{\repdecn@mb{\v@leur}}\v@lY=\z@\v@lX=\v@leur%
    \resetc@ntr@l{2}\figv@ctCreg-3(\v@lX,\v@lY)\figpttraDD-5:=#1/1,-3/%
    \figptrotDD-5:=-5/#1,#4/%
    \figvectPDD-3[#1,-5]\figg@tXY{-3}\v@leur=\v@lX%
    \ifdim#3pt<#4pt\v@lX=\v@lY\v@lY=-\v@leur\else\v@lX=-\v@lY\v@lY=\v@leur\fi%
    \figv@ctCreg-3(\v@lX,\v@lY)\vecunit@{-3}{-3}%
    \if@rrowratio\v@leur=#4pt\advance\v@leur-#3pt\maxim@m{\mili@u}{-\v@leur}{\v@leur}%
    \mili@u=\degT@rd\mili@u\v@leur=\s@glen\mili@u\edef\s@glen{\repdecn@mb{\v@leur}}%
    \mili@u=#2\mili@u\mili@u=\@rrowheadratio\mili@u\else\mili@u=\@rrowheadlength pt\fi%
    \figpttraDD-6:=-5/\s@glen,-3/\v@leur=#2pt\v@leur=2\v@leur%
    \invers@{\v@leur}{\v@leur}\c@rre=\repdecn@mb{\v@leur}\mili@u
    \mili@u=\c@rre\mili@u=\repdecn@mb{\c@rre}\mili@u%
    \v@leur=\p@\advance\v@leur-\mili@u
    \invers@{\mili@u}{2\v@leur}\delt@=\c@rre\delt@=\repdecn@mb{\mili@u}\delt@%
    \xdef\sDcc@ngle{\repdecn@mb{\delt@}}
    \sqrt@{\mili@u}{\v@leur}\arct@n\v@leur(\mili@u,\c@rre)%
    \v@leur=\rdT@deg\v@leur
    \ifdim#3pt<#4pt\v@leur=-\v@leur\fi%
    \if@rrowhout\v@leur=-\v@leur\fi\edef\cor@ngle{\repdecn@mb{\v@leur}}%
    \figptrotDD-6:=-6/-5,\cor@ngle/\psarrowheadDD[-6,-5]%
    \psc@mment{End pscirc@rrowhead}}}
\def\psarrowcircPDD#1;#2[#3,#4]{{\ifps@cri%
    \psc@mment{psarrowcircPDD Center=#1; Radius=#2, [P1=#3,P2=#4]}%
    \s@uvc@ntr@l\et@tpsarrowcircPDD\ps@ngleparam#1;#2[#3,#4]%
    \ifdim\v@leur>\z@\ifdim\v@lmin>\v@lmax\advance\v@lmax\PI@deg\fi%
    \else\ifdim\v@lmin<\v@lmax\advance\v@lmin\PI@deg\fi\fi%
    \edef\@ngdeb{\repdecn@mb{\v@lmin}}\edef\@ngfin{\repdecn@mb{\v@lmax}}%
    \psarrowcirc#1;\r@dius(\@ngdeb,\@ngfin)%
    \psc@mment{End psarrowcircPDD}\resetc@ntr@l\et@tpsarrowcircPDD\fi}}
\def\psarrowcircPTD#1;#2[#3,#4,#5]{{\ifps@cri\s@uvc@ntr@l\et@tpsarrowcircPTD%
    \psc@mment{psarrowcircPTD Center=#1; Radius=#2, [P1=#3,P2=#4,P3=#5]}%
    \figgetangleTD\@ngfin[#1,#3,#4,#5]\v@leur=#2pt%
    \maxim@m{\mili@u}{-\v@leur}{\v@leur}\edef\r@dius{\repdecn@mb{\mili@u}}%
    \ifdim\v@leur<\z@\v@lmax=\@ngfin pt\advance\v@lmax-\PI@deg%
    \edef\@ngfin{\repdecn@mb{\v@lmax}}\fi\psarrowcircTD#1,#3,#5;\r@dius(0,\@ngfin)%
    \psc@mment{End psarrowcircPTD}\resetc@ntr@l\et@tpsarrowcircPTD\fi}}
\def\psbeginfig#1{\edef\psfilen@me{#1}\edef\auxfilen@me{\jobname.anx}%
    \ifpstestm@de\ps@critrue\else\openin\frf@g=\psfilen@me\relax%
    \ifeof\frf@g\ps@critrue\else\ps@crifalse\fi\closein\frf@g\fi%
    \ifps@cri\initb@undb@x\initpss@ttings%
    \immediate\openout\fwf@g=\auxfilen@me%
    \pssetgray{0}\pssetwidth{\defaultwidth}\fi%
    \curr@ntPStrue\c@ldefproj}
\def\initpss@ttings{\pssetfillmode{no}\pssetupdate{no}%
    \psresetarrowhead\pssetroundness{\defaultroundness}\pssetmeshdiag{\defaultmeshdiag}%
    \def\curr@ntcolor{0}\def\curr@ntcolorc@md{setgray}%
    \edef\curr@ntdash{\defaultdash}\edef\curr@ntwidth{\defaultwidth}%
    \psresetsecondsettings\Use@llipsefalse}
\def\B@zierBB@x#1#2(#3,#4,#5,#6){{\c@rre=\t@n\epsil@n
    \v@lmax=#4\advance\v@lmax-#5\v@lmax=\thr@@\v@lmax\advance\v@lmax#6\advance\v@lmax-#3%
    \mili@u=#4\mili@u=-\tw@\mili@u\advance\mili@u#3\advance\mili@u#5%
    \v@lmin=#4\advance\v@lmin-#3\maxim@m{\v@leur}{-\v@lmax}{\v@lmax}%
    \maxim@m{\delt@}{-\mili@u}{\mili@u}\maxim@m{\v@leur}{\v@leur}{\delt@}%
    \maxim@m{\delt@}{-\v@lmin}{\v@lmin}\maxim@m{\v@leur}{\v@leur}{\delt@}%
    \ifdim\v@leur>\c@rre\invers@{\v@leur}{\v@leur}\edef\Uns@rM@x{\repdecn@mb{\v@leur}}%
    \v@lmax=\Uns@rM@x\v@lmax\mili@u=\Uns@rM@x\mili@u\v@lmin=\Uns@rM@x\v@lmin%
    \maxim@m{\v@leur}{-\v@lmax}{\v@lmax}\ifdim\v@leur<\c@rre%
    \maxim@m{\v@leur}{-\mili@u}{\mili@u}\ifdim\v@leur<\c@rre\else%
    \invers@{\mili@u}{\mili@u}\v@leur=-0.5\v@lmin%
    \v@leur=\repdecn@mb{\mili@u}\v@leur\m@jBBB@x{\v@leur}{#1}{#2}(#3,#4,#5,#6)\fi%
    \else\delt@=\repdecn@mb{\mili@u}\mili@u\v@leur=\repdecn@mb{\v@lmax}\v@lmin%
    \advance\delt@-\v@leur\ifdim\delt@<\z@\else\invers@{\v@lmax}{\v@lmax}%
    \edef\Uns@rAp{\repdecn@mb{\v@lmax}}\sqrt@{\delt@}{\delt@}%
    \v@leur=-\mili@u\advance\v@leur\delt@\v@leur=\Uns@rAp\v@leur%
    \m@jBBB@x{\v@leur}{#1}{#2}(#3,#4,#5,#6)%
    \v@leur=-\mili@u\advance\v@leur-\delt@\v@leur=\Uns@rAp\v@leur%
    \m@jBBB@x{\v@leur}{#1}{#2}(#3,#4,#5,#6)\fi\fi\fi}}
\def\m@jBBB@x#1#2#3(#4,#5,#6,#7){{\relax\ifdim#1>\z@\ifdim#1<\p@%
    \edef\T@{\repdecn@mb{#1}}\v@lX=\p@\advance\v@lX-#1\edef\UNmT@{\repdecn@mb{\v@lX}}%
    \v@lX=#4\v@lY=#5\v@lZ=#6\v@lXa=#7\v@lX=\UNmT@\v@lX\advance\v@lX\T@\v@lY%
    \v@lY=\UNmT@\v@lY\advance\v@lY\T@\v@lZ\v@lZ=\UNmT@\v@lZ\advance\v@lZ\T@\v@lXa%
    \v@lX=\UNmT@\v@lX\advance\v@lX\T@\v@lY\v@lY=\UNmT@\v@lY\advance\v@lY\T@\v@lZ%
    \v@lX=\UNmT@\v@lX\advance\v@lX\T@\v@lY%
    \ifcase#2\or\v@lY=#3\or\v@lY=\v@lX\v@lX=#3\fi\b@undb@x{\v@lX}{\v@lY}\fi\fi}}
\def\psB@zier#1[#2]{{\immediate\write\fwf@g{newpath}%
    \s@mme=\z@\def\list@num{#2,0}%
    \extrairelepremi@r\p@int\de\list@num\pswrit@cmdS{\p@int}{moveto}{\fwf@g}{\X@un}{\Y@un}%
    \loop\ifnum\s@mme<#1\advance\s@mme\@ne\BdingB@xfalse%
    \extrairelepremi@r\p@int\de\list@num\pswrit@cmdS{\p@int}{}{\fwf@g}{\X@de}{\Y@de}%
    \extrairelepremi@r\p@int\de\list@num\pswrit@cmdS{\p@int}{}{\fwf@g}{\X@tr}{\Y@tr}%
    \BdingB@xtrue%
    \extrairelepremi@r\p@int\de\list@num\pswrit@cmdS{\p@int}{curveto}{\fwf@g}{\X@qu}{\Y@qu}%
    \B@zierBB@x{1}{\Y@un}(\X@un,\X@de,\X@tr,\X@qu)%
    \B@zierBB@x{2}{\X@un}(\Y@un,\Y@de,\Y@tr,\Y@qu)%
    \edef\X@un{\X@qu}\edef\Y@un{\Y@qu}\repeat}}
\def\psBezierDD#1[#2]{\ifps@cri\psc@mment{psBezierDD N arcs=#1, Control points=#2}%
    \iffillm@de\immediate\write\fwf@g{currentrgbcolor}\psB@zier#1[#2]%
    \immediate\write\fwf@g{\curr@ntcolor\setfillc@md\space setrgbcolor}%
    \else\psB@zier#1[#2]\immediate\write\fwf@g{stroke}\fi%
    \psc@mment{End psBezierDD}\fi}
\def\psBezierTD#1[#2]{\ifps@cri\s@uvc@ntr@l\et@tpsBezierTD%
    \psc@mment{psBezierTD N arcs=#1, Control points=#2}%
    \iffillm@de\immediate\write\fwf@g{currentrgbcolor}\psB@zierTD#1[#2]%
    \immediate\write\fwf@g{\curr@ntcolor\setfillc@md\space setrgbcolor}%
    \else\psB@zierTD#1[#2]\immediate\write\fwf@g{stroke}\fi%
    \psc@mment{End psBezierTD}\resetc@ntr@l\et@tpsBezierTD\fi}
\def\psB@zierTD#1[#2]{\ifnum\curr@ntproj<\tw@\psB@zier#1[#2]\else\psB@zier@TD#1[#2]\fi}
\def\psB@zier@TD#1[#2]{{\immediate\write\fwf@g{newpath}%
    \s@mme=\z@\def\list@num{#2,0}\extrairelepremi@r\p@int\de\list@num%
    \let\c@lprojSP=\relax\setc@ntr@l{2}\figptpr@j-7:/\p@int/%
    \pswrit@cmd{-7}{moveto}{\fwf@g}%
    \loop\ifnum\s@mme<#1\advance\s@mme\@ne\extrairelepremi@r\p@intun\de\list@num%
    \extrairelepremi@r\p@intde\de\list@num\extrairelepremi@r\p@inttr\de\list@num%
    \subB@zierTD\NBz@rcs[\p@int,\p@intun,\p@intde,\p@inttr]\edef\p@int{\p@inttr}\repeat}}
\def\subB@zierTD#1[#2,#3,#4,#5]{\delt@=\p@\divide\delt@\NBz@rcs\v@lmin=\z@%
    {\figg@tXY{-7}\edef\X@un{\the\v@lX}\edef\Y@un{\the\v@lY}%
    \s@mme=\z@\loop\ifnum\s@mme<#1\advance\s@mme\@ne%
    \v@leur=\v@lmin\advance\v@leur0.33333\delt@\edef\unti@rs{\repdecn@mb{\v@leur}}%
    \v@leur=\v@lmin\advance\v@leur0.66666\delt@\edef\deti@rs{\repdecn@mb{\v@leur}}%
    \advance\v@lmin\delt@\edef\trti@rs{\repdecn@mb{\v@lmin}}%
    \figptBezierTD-8::\trti@rs[#2,#3,#4,#5]\figptpr@j-8:/-8/%
    \c@lsubBzarc\unti@rs,\deti@rs[#2,#3,#4,#5]\BdingB@xfalse%
    \pswrit@cmdS{-4}{}{\fwf@g}{\X@de}{\Y@de}\pswrit@cmdS{-3}{}{\fwf@g}{\X@tr}{\Y@tr}%
    \BdingB@xtrue\pswrit@cmdS{-8}{curveto}{\fwf@g}{\X@qu}{\Y@qu}%
    \B@zierBB@x{1}{\Y@un}(\X@un,\X@de,\X@tr,\X@qu)%
    \B@zierBB@x{2}{\X@un}(\Y@un,\Y@de,\Y@tr,\Y@qu)%
    \edef\X@un{\X@qu}\edef\Y@un{\Y@qu}\figptcopyDD-7:/-8/\repeat}}
\def\NBz@rcs{2}
\def\c@lsubBzarc#1,#2[#3,#4,#5,#6]{\figptBezierTD-5::#1[#3,#4,#5,#6]%
    \figptBezierTD-6::#2[#3,#4,#5,#6]\figptpr@j-4:/-5/\figptpr@j-5:/-6/%
    \figptscontrolDD-4[-7,-4,-5,-8]}
\def\pscircDD#1(#2){\ifps@cri\psc@mment{pscircDD Center=#1 (Radius=#2)}%
    \psarccircDD#1;#2(0,360)\psc@mment{End pscircDD}\fi}
\def\pscircTD#1,#2,#3(#4){\ifps@cri%
    \psc@mment{pscircTD Center=#1,P1=#2,P2=#3 (Radius=#4)}%
    \psarccircTD#1,#2,#3;#4(0,360)\psc@mment{End pscircTD}\fi}
{\catcode`\%=12\gdef\p@urcent{
\def\psc@mment#1{\ifpsdebugmode\immediate\write\fwf@g{\p@urcent\space#1}\fi}
{\catcode`\[=1\catcode`\{=12\gdef\acc@louv[{}}
{\catcode`\]=2\catcode`\}=12\gdef\acc@lfer{}]]
\def\psdict@{\ifUse@llipse%
    \immediate\write\fwf@g{/ellipsedict 9 dict def ellipsedict /mtrx matrix put}%
    \immediate\write\fwf@g{/ellipse \acc@louv ellipsedict begin}%
    \immediate\write\fwf@g{ /endangle exch def /startangle exch def}%
    \immediate\write\fwf@g{ /yrad exch def /xrad exch def}%
    \immediate\write\fwf@g{ /rotangle exch def /y exch def /x exch def}%
    \immediate\write\fwf@g{ /savematrix mtrx currentmatrix def}%
    \immediate\write\fwf@g{ x y translate rotangle rotate xrad yrad scale}%
    \immediate\write\fwf@g{ 0 0 1 startangle endangle arc}%
    \immediate\write\fwf@g{ savematrix setmatrix end\acc@lfer def}%
    \fi\pshe@der{EndProlog}}
\def\pssetroundness#1{\edef\curv@roundness{#1}}
\def\defaultroundness{0.2}         
\def\pscurveDD[#1]{{\ifps@cri\psc@mment{pscurveDD Points=#1}\s@uvc@ntr@l\et@tpscurveDD%
    \iffillm@de\immediate\write\fwf@g{currentrgbcolor}\psc@rveDD\curv@roundness[#1]%
    \immediate\write\fwf@g{\curr@ntcolor\setfillc@md\space setrgbcolor}%
    \else\psc@rveDD\curv@roundness[#1]\immediate\write\fwf@g{stroke}\fi%
    \psc@mment{End pscurveDD}\resetc@ntr@l\et@tpscurveDD\fi}}
\def\pscurveTD[#1]{{\ifps@cri%
    \psc@mment{pscurveTD Points=#1}\s@uvc@ntr@l\et@tpscurveTD\let\c@lprojSP=\relax%
    \iffillm@de\immediate\write\fwf@g{currentrgbcolor}\psc@rveTD\curv@roundness[#1]%
    \immediate\write\fwf@g{\curr@ntcolor\setfillc@md\space setrgbcolor}%
    \else\psc@rveTD\curv@roundness[#1]\immediate\write\fwf@g{stroke}\fi%
    \psc@mment{End pscurveTD}\resetc@ntr@l\et@tpscurveTD\fi}}
\def\psc@rveDD#1[#2]{%
    \def\list@num{#2}\extrairelepremi@r\Ak@\de\list@num%
    \extrairelepremi@r\Ai@\de\list@num\extrairelepremi@r\Aj@\de\list@num%
    \immediate\write\fwf@g{newpath}\pswrit@cmdS{\Ai@}{moveto}{\fwf@g}{\X@un}{\Y@un}%
    \setc@ntr@l{2}\figvectPDD -1[\Ak@,\Aj@]%
    \@ecfor\Ak@:=\list@num\do{\figpttraDD-2:=\Ai@/#1,-1/\BdingB@xfalse%
       \pswrit@cmdS{-2}{}{\fwf@g}{\X@de}{\Y@de}%
       \figvectPDD -1[\Ai@,\Ak@]\figpttraDD-2:=\Aj@/-#1,-1/%
       \pswrit@cmdS{-2}{}{\fwf@g}{\X@tr}{\Y@tr}\BdingB@xtrue%
       \pswrit@cmdS{\Aj@}{curveto}{\fwf@g}{\X@qu}{\Y@qu}%
       \B@zierBB@x{1}{\Y@un}(\X@un,\X@de,\X@tr,\X@qu)%
       \B@zierBB@x{2}{\X@un}(\Y@un,\Y@de,\Y@tr,\Y@qu)%
       \edef\X@un{\X@qu}\edef\Y@un{\Y@qu}\edef\Ai@{\Aj@}\edef\Aj@{\Ak@}}}
\def\psc@rveTD#1[#2]{\setc@ntr@l{2}%
    \def\list@num{#2}\extrairelepremi@r\Ak@\de\list@num\figptpr@j-5:/\Ak@/%
    \extrairelepremi@r\Ai@\de\list@num\figptpr@j-3:/\Ai@/%
    \extrairelepremi@r\Aj@\de\list@num\figptpr@j-4:/\Aj@/%
    \immediate\write\fwf@g{newpath}\pswrit@cmdS{-3}{moveto}{\fwf@g}{\X@un}{\Y@un}%
    \figvectPDD -1[-5,-4]%
    \@ecfor\Ak@:=\list@num\do{\figptpr@j-5:/\Ak@/\figpttraDD-2:=-3/#1,-1/%
       \BdingB@xfalse\pswrit@cmdS{-2}{}{\fwf@g}{\X@de}{\Y@de}%
       \figvectPDD -1[-3,-5]\figpttraDD-2:=-4/-#1,-1/%
       \pswrit@cmdS{-2}{}{\fwf@g}{\X@tr}{\Y@tr}\BdingB@xtrue%
       \pswrit@cmdS{-4}{curveto}{\fwf@g}{\X@qu}{\Y@qu}%
       \B@zierBB@x{1}{\Y@un}(\X@un,\X@de,\X@tr,\X@qu)%
       \B@zierBB@x{2}{\X@un}(\Y@un,\Y@de,\Y@tr,\Y@qu)%
       \edef\X@un{\X@qu}\edef\Y@un{\Y@qu}\figptcopyDD-3:/-4/\figptcopyDD-4:/-5/}}
\def\psendfig{\ifps@cri\immediate\closeout\fwf@g%
    \immediate\openout\fwf@g=\psfilen@me%
    \immediate\write\fwf@g{\p@urcent\string!PS-Adobe-2.0 EPSF-2.0}%
    \pshe@der{Creator\string: TeX (fig4tex.tex)}%
    \pshe@der{Title\string: \psfilen@me}%
    \pshe@der{CreationDate\string: \the\day/\the\month/\the\year}%
    {\v@lX=\ptT@ptps\c@@rdXmin\v@lY=\ptT@ptps\c@@rdYmin%
     \v@lXa=\ptT@ptps\c@@rdXmax\v@lYa=\ptT@ptps\c@@rdYmax%
    \pshe@der{BoundingBox\string: \repdecn@mb{\v@lX}\space\repdecn@mb{\v@lY}%
              \space\repdecn@mb{\v@lXa}\space\repdecn@mb{\v@lYa}}}%
    \pshe@der{EndComments}\psdict@%
    \immediate\write\fwf@g{currentrgbcolor}%
    \openin\frf@g=\auxfilen@me\c@pypsfile\fwf@g\frf@g\closein\frf@g%
    \immediate\write\fwf@g{setrgbcolor\space showpage}%
    \psc@mment{End of file.}\immediate\closeout\fwf@g%
    \immediate\openout\fwf@g=\auxfilen@me\immediate\closeout\fwf@g%
    \immediate\write16{File \psfilen@me\space created.}\fi\curr@ntPSfalse\ps@critrue}
\def\pshe@der#1{\immediate\write\fwf@g{\p@urcent\p@urcent#1}}
\def\psline[#1]{{\ifps@cri\psc@mment{psline Points=#1}%
    \let\pslign@=\pslign@P\pslin@{#1}\psc@mment{End psline}\fi}}
\def\pslineF#1{{\ifps@cri\psc@mment{pslineF Filename=#1}%
    \let\pslign@=\pslign@F\pslin@{#1}\psc@mment{End pslineF}\fi}}
\def\pslin@#1{\iffillm@de\immediate\write\fwf@g{currentrgbcolor}\pslign@{#1}%
    \immediate\write\fwf@g{\curr@ntcolor\setfillc@md\space setrgbcolor}%
    \else\pslign@{#1}\ifx\derp@int\premp@int%
    \immediate\write\fwf@g{closepath\space stroke}%
    \else\immediate\write\fwf@g{stroke}\fi\fi}
\def\pslign@P#1{\def\list@num{#1}\extrairelepremi@r\p@int\de\list@num%
    \edef\premp@int{\p@int}\immediate\write\fwf@g{newpath}%
    \pswrit@cmd{\p@int}{moveto}{\fwf@g}%
    \@ecfor\p@int:=\list@num\do{\pswrit@cmd{\p@int}{lineto}{\fwf@g}%
    \edef\derp@int{\p@int}}}
\def\pslign@F#1{\s@uvc@ntr@l\et@tpslign@F\setc@ntr@l{2}\openin\frf@g=#1\relax%
    \ifeof\frf@g\message{*** File #1 not found !}\end\else%
    \read\frf@g to\tr@c\edef\premp@int{\tr@c}\expandafter\extr@ctCF\tr@c:%
    \immediate\write\fwf@g{newpath}\pswrit@cmd{-1}{moveto}{\fwf@g}%
    \loop\read\frf@g to\tr@c\ifeof\frf@g\mored@tafalse\else\mored@tatrue\fi%
    \ifmored@ta\expandafter\extr@ctCF\tr@c:\pswrit@cmd{-1}{lineto}{\fwf@g}%
    \edef\derp@int{\tr@c}\repeat\fi\closein\frf@g\resetc@ntr@l\et@tpslign@F}
\def\extr@ctCFDD#1 #2:{\v@lX=#1\unit@\v@lY=#2\unit@\figp@intregDD-1:(\v@lX,\v@lY)}
\def\extr@ctCFTD#1 #2 #3:{\v@lX=#1\unit@\v@lY=#2\unit@\v@lZ=#3\unit@%
    \figp@intregTD-1:(\v@lX,\v@lY,\v@lZ)}
\newcount\c@ntrolmesh
\def\pssetmeshdiag#1{\c@ntrolmesh=#1}
\def\defaultmeshdiag{0}          
\def\psmesh#1,#2[#3,#4,#5,#6]{{\ifps@cri%
    \psc@mment{psmesh N1=#1, N2=#2, Quadrangle=[#3,#4,#5,#6]}%
    \s@uvc@ntr@l\et@tpsmesh\pss@tsecondSt\setc@ntr@l{2}%
    \ifnum#1>\@ne\psmeshp@rt#1[#3,#4,#5,#6]\fi%
    \ifnum#2>\@ne\psmeshp@rt#2[#4,#5,#6,#3]\fi%
    \ifnum\c@ntrolmesh>\z@\psmeshdi@g#1,#2[#3,#4,#5,#6]\fi%
    \ifnum\c@ntrolmesh<\z@\psmeshdi@g#2,#1[#4,#5,#6,#3]\fi\psrest@reSt%
    \psline[#3,#4,#5,#6,#3]\psc@mment{End psmesh}\resetc@ntr@l\et@tpsmesh\fi}}
\def\psmeshp@rt#1[#2,#3,#4,#5]{{\l@mbd@un=\@ne\l@mbd@de=#1\loop%
    \ifnum\l@mbd@un<#1\advance\l@mbd@de\m@ne\figptbary-1:[#2,#3;\l@mbd@de,\l@mbd@un]%
    \figptbary-2:[#5,#4;\l@mbd@de,\l@mbd@un]\psline[-1,-2]\advance\l@mbd@un\@ne\repeat}}
\def\psmeshdi@g#1,#2[#3,#4,#5,#6]{\figptcopy-2:/#3/\figptcopy-3:/#6/%
    \l@mbd@un=\z@\l@mbd@de=#1\loop\ifnum\l@mbd@un<#1%
    \advance\l@mbd@un\@ne\advance\l@mbd@de\m@ne\figptcopy-1:/-2/\figptcopy-4:/-3/%
    \figptbary-2:[#3,#4;\l@mbd@de,\l@mbd@un]%
    \figptbary-3:[#6,#5;\l@mbd@de,\l@mbd@un]\psmeshdi@gp@rt#2[-1,-2,-3,-4]\repeat}
\def\psmeshdi@gp@rt#1[#2,#3,#4,#5]{{\l@mbd@un=\z@\l@mbd@de=#1\loop%
    \ifnum\l@mbd@un<#1\figptbary-5:[#2,#5;\l@mbd@de,\l@mbd@un]%
    \advance\l@mbd@de\m@ne\advance\l@mbd@un\@ne%
    \figptbary-6:[#3,#4;\l@mbd@de,\l@mbd@un]\psline[-5,-6]\repeat}}
\def\psnormalDD#1,#2[#3,#4]{{\ifps@cri%
    \psc@mment{psnormal Length=#1, Lambda=#2 [Pt1,Pt2]=[#3,#4]}%
    \s@uvc@ntr@l\et@tpsnormal\resetc@ntr@l{2}\figptendnormal-6::#1,#2[#3,#4]%
    \figptcopyDD-5:/-1/\psarrow[-5,-6]%
    \psc@mment{End psnormal}\resetc@ntr@l\et@tpsnormal\fi}}
\def\pssetcmyk#1{\ifps@cri\psc@mment{pssetcmyk Color=#1}%
    \def\curr@ntcolor{#1}\def\setfillc@md{\space setcmykcolor\space fill}%
    \def\curr@ntcolorc@md{setcmykcolor}\pssetsecondcmyk{#1}%
    \immediate\write\fwf@g{#1\space setcmykcolor}\fi}
\def\pssetrgb#1{\ifps@cri\psc@mment{pssetrgb Color=#1}%
    \def\curr@ntcolor{#1}\def\setfillc@md{\space setrgbcolor \space fill}%
    \def\curr@ntcolorc@md{setrgbcolor}\pssetsecondrgb{#1}%
    \immediate\write\fwf@g{#1\space setrgbcolor}\fi}
\def\pssetgray#1{\ifps@cri\psc@mment{pssetgray Gray level=#1}%
    \def\curr@ntcolor{#1}\def\setfillc@md{\space setgray\space fill}%
    \def\curr@ntcolorc@md{setgray}\pssetsecondgray{#1}%
    \immediate\write\fwf@g{#1\space setgray}\fi}
\def\Blackcmyk{0 0 0 1}
\def\Whitecmyk{0 0 0 0}
\def\Cyancmyk{1 0 0 0}
\def\Magentacmyk{0 1 0 0}
\def\Yellowcmyk{0 0 1 0}
\def\Redcmyk{0 1 1 0}
\def\Greencmyk{1 0 1 0}
\def\Bluecmyk{1 1 0 0}
\def\Graycmyk{0 0 0 0.50}
\def\BrickRedcmyk{0 0.89 0.94 0.28} 
\def\Browncmyk{0 0.81 1 0.60} 
\def\ForestGreencmyk{0.91 0 0.88 0.12} 
\def\Goldenrodcmyk{ 0 0.10 0.84 0} 
\def\Marooncmyk{0 0.87 0.68 0.32} 
\def\Orangecmyk{0 0.61 0.87 0} 
\def\Purplecmyk{0.45 0.86 0 0} 
\def\RoyalBluecmyk{1. 0.50 0 0} 
\def\Violetcmyk{0.79 0.88 0 0} 
\def\Blackrgb{0 0 0}
\def\Whitergb{1 1 1}
\def\Redrgb{1 0 0}
\def\Greenrgb{0 1 0}
\def\Bluergb{0 0 1}
\def\Cyanrgb{0 1 1}
\def\Magentargb{1 0 1}
\def\Yellowrgb{1 1 0}
\def\Grayrgb{0.5 0.5 0.5}
\def\Chocolatergb{0.824 0.412 0.118}
\def\DarkGoldenrodrgb{0.722 0.525 0.043}
\def\DarkOrangergb{1 0.549 0}
\def\Firebrickrgb{0.698 0.133 0.133}
\def\ForestGreenrgb{0.133 0.545 0.133}
\def\Goldrgb{1 0.843 0}
\def\HotPinkrgb{1 0.412 0.706}
\def\Maroonrgb{0.690 0.188 0.376}
\def\Pinkrgb{1 0.753 0.796}
\def\RoyalBluergb{0.255 0.412 0.882}
\def\s@uvdash#1{\edef#1{\curr@ntdash}}
\def\defaultdash{1}             
\def\pssetdash#1{\ifps@cri\pssetd@sh#1 :\edef\curr@ntdash{#1}\fi}
\def\pssetd@shI#1{\psc@mment{pssetdash Index=#1}\ifcase#1%
    \or\immediate\write\fwf@g{[] 0 setdash}
    \or\immediate\write\fwf@g{[6 2] 0 setdash}
    \or\immediate\write\fwf@g{[4 2] 0 setdash}
    \or\immediate\write\fwf@g{[2 2] 0 setdash}
    \or\immediate\write\fwf@g{[1 2] 0 setdash}
    \or\immediate\write\fwf@g{[2 4] 0 setdash}
    \or\immediate\write\fwf@g{[3 5] 0 setdash}
    \or\immediate\write\fwf@g{[3 3] 0 setdash}
    \or\immediate\write\fwf@g{[3 5 1 5] 0 setdash}
    \or\immediate\write\fwf@g{[6 4 2 4] 0 setdash}
    \fi}
\def\pssetd@sh#1 #2:{{\def\t@xt@{#1}\ifx\t@xt@\empty\pssetd@sh#2:
    \else\def\t@xt@{#2}\ifx\t@xt@\empty\pssetd@shI{#1}\else\s@mme=\@ne\def\debutp@t{#1}%
    \an@lysd@sh#2:\ifodd\s@mme\edef\debutp@t{\debutp@t\space\finp@t}\def\finp@t{0}\fi%
    \psc@mment{pssetdash Pattern=#1 #2}%
    \immediate\write\fwf@g{[\debutp@t] \finp@t\space setdash}\fi\fi}}
\def\an@lysd@sh#1 #2:{\def\t@xt@{#2}\ifx\t@xt@\empty\def\finp@t{#1}\else%
    \edef\debutp@t{\debutp@t\space#1}\advance\s@mme\@ne\an@lysd@sh#2:\fi}
\def\s@uvwidth#1{\edef#1{\curr@ntwidth}}
\def\defaultwidth{0.4}            
\def\pssetwidth#1{\ifps@cri\psc@mment{pssetwidth Width=#1}%
    \immediate\write\fwf@g{#1 setlinewidth}\edef\curr@ntwidth{#1}\fi}
\def\pssetseconddash#1{\edef\curr@ntseconddash{#1}}
\def\defaultseconddash{4} 
\def\pssetsecondwidth#1{\edef\curr@ntsecondwidth{#1}}
\edef\defaultsecondwidth{\defaultwidth} 
\def\psresetsecondsettings{%
    \pssetseconddash{\defaultseconddash}\pssetsecondwidth{\defaultsecondwidth}%
    \def\sec@ndcolor{\curr@ntcolor}\def\sec@ndcolorc@md{\curr@ntcolorc@md}}
\def\pssetsecondcmyk#1{\def\sec@ndcolor{#1}\def\sec@ndcolorc@md{setcmykcolor}}
\def\pssetsecondrgb#1{\def\sec@ndcolor{#1}\def\sec@ndcolorc@md{setrgbcolor}}
\def\pssetsecondgray#1{\def\sec@ndcolor{#1}\def\sec@ndcolorc@md{setgray}}
\def\pss@tsecondSt{%
    \s@uvdash{\typ@dash}\pssetdash{\curr@ntseconddash}%
    \s@uvwidth{\typ@width}\pssetwidth{\curr@ntsecondwidth}%
    \immediate\write\fwf@g{\sec@ndcolor\space\sec@ndcolorc@md}}
\def\psrest@reSt{\pssetwidth{\typ@width}\pssetdash{\typ@dash}%
    \immediate\write\fwf@g{\curr@ntcolor\space\curr@ntcolorc@md}}
\def\pstrimesh#1[#2,#3,#4]{{\ifps@cri\psc@mment{pstrimesh Type=#1, Triangle=[#2,#3,#4]}%
    \s@uvc@ntr@l\et@tpstrimesh\ifnum#1>\@ne\pss@tsecondSt\setc@ntr@l{2}%
    \pstrimeshp@rt#1[#2,#3,#4]\pstrimeshp@rt#1[#3,#4,#2]%
    \pstrimeshp@rt#1[#4,#2,#3]\psrest@reSt\fi\psline[#2,#3,#4,#2]%
    \psc@mment{End pstrimesh}\resetc@ntr@l\et@tpstrimesh\fi}}
\def\pstrimeshp@rt#1[#2,#3,#4]{{\l@mbd@un=\@ne\l@mbd@de=#1\loop\ifnum\l@mbd@de>\@ne%
    \advance\l@mbd@de\m@ne\figptbary-1:[#2,#3;\l@mbd@de,\l@mbd@un]%
    \figptbary-2:[#2,#4;\l@mbd@de,\l@mbd@un]\psline[-1,-2]%
    \advance\l@mbd@un\@ne\repeat}}
\def\pswrit@cmd#1#2#3{{\figg@tXY{#1}\c@lprojSP\b@undb@x{\v@lX}{\v@lY}%
    \v@lX=\ptT@ptps\v@lX\v@lY=\ptT@ptps\v@lY%
    \immediate\write#3{\repdecn@mb{\v@lX}\space\repdecn@mb{\v@lY}\space#2}}}
\def\pswrit@cmdS#1#2#3#4#5{{\figg@tXY{#1}\c@lprojSP\b@undb@x{\v@lX}{\v@lY}%
    \global\result@t=\v@lX\global\result@@t=\v@lY%
    \v@lX=\ptT@ptps\v@lX\v@lY=\ptT@ptps\v@lY%
    \immediate\write#3{\repdecn@mb{\v@lX}\space\repdecn@mb{\v@lY}\space#2}}%
    \edef#4{\the\result@t}\edef#5{\the\result@@t}}
\initpr@lim\initpss@ttings
\catcode`\@=12

\input epsf
\newtheorem{thmalp}{Theorem}
\renewcommand{\thethmalp}{\Alph{thmalp}}
\newtheorem{lemalp}{Lemma}
\renewcommand{\thelemalp}{\Alph{lemalp}}

\newtheorem{thm}{Theorem}
\newtheorem{lemma}{Lemma}
\newtheorem{corr}[thm]{Corollary}
\newtheorem{remark}{Remark}
\newtheorem{deff}[thm]{Definition}
\newtheorem{case}[thm]{Case}
\newtheorem{prop}{Proposition}
\numberwithin{equation}{section}
\newcommand{\uG}{\underline{G}}
\newcommand{\uT}{\underline{T}}
\newcommand{\uB}{\underline{B}}
\newcommand{\uU}{\underline{U}}
\newcommand{\uA}{\underline{A}}
\newcommand{\uM}{\underline{M}}
\newcommand{\uN}{\underline{N}}
\newcommand{\uP}{\underline{P}}
\newcommand{\sN}{\cal N}
\newcommand{\bR}{\mathbb R}
\newcommand{\bC}{\Bbb C}
\newcommand{\ds}{\displaystyle}
\newcommand{\pa}{{\cal P}_{\alpha}}
\date{}

\begin{titlepage}
\title{\bf Sharp Integrability for Brownian Motion in  Parabola-shaped Regions}
\author{Rodrigo Ba\~nuelos\thanks{Supported in part by NSF Grant
\# 9700585-DMS}\\Department of  Mathematics
\\Purdue University\\West Lafayette, IN
47906\and Tom Carroll\\ Department of  Mathematics\\National
University of Ireland\\Cork, Ireland} \maketitle
\begin{abstract}
\noindent We study the sharp order of integrability of the exit
position of Brownian motion from  the planar domains ${\cal
P}_\alpha = \{(x,y)\in \bR\times \bR\colon x> 0,\ |y| <
Ax^{\alpha}\}$, $0<\alpha<1$. Together with some simple
good-$\lambda$ type arguments, this implies the  order of
integrability for the exit time of these domains; a result first
proved for $\alpha =1/2$ by Ba\~nuelos, DeBlassie and Smits
\cite{ba} and for general $\alpha$ by Li \cite{li}. A sharp version
 of this result  is also proved in
higher dimensions.

\bigskip
\centerline{\bf  Contents}

\begin{itemize}
\item[{\S1.}] {\sl Introduction}
\item[\S2.] {\sl Proofs of Theorems 1 and 4}
\item[\S3.] {\sl Parabola-shaped domains in the plane}
\begin{itemize}
\item[\S3.1]{\sl\footnotesize Harmonic measure estimates}
\item[\S3.2]{\sl\footnotesize Proof of (3.3) by conformal mapping}
\item[\S3.3]{\sl\footnotesize Proof of Theorem 2}
\end{itemize}
\item[\S4.] {\sl Parabola-shaped regions in $\bR^n$}
\begin{itemize}
\item[\S4.1]{\sl\footnotesize Carleman method: Upper bound for harmonic measure}
\item[\S4.2]{\sl\footnotesize Conformal mapping method: Lower bound for harmonic measure}
\begin{itemize}
\item[\S4.2.1]{\sl\footnotesize From a parabolic--shaped region in $\bR^n$ to a planar
strip}
\item[\S4.2.2]{\sl\footnotesize Asymptotic estimates for the conformal mapping $g$}
\item[\S4.2.3]{\sl\footnotesize Asymptotic form of the differential operator}
\item[\S4.2.4]{\sl\footnotesize Sub solutions and a maximimum principle}
\item[\S4.2.5]{\sl\footnotesize Rate of  exponential decay of solutions}
\item[\S4.2.6]{\sl\footnotesize Lower bound for harmonic measure}
\item[\S4.2.7]{\sl\footnotesize Concluding remarks}
\end{itemize}
\end{itemize}
\end{itemize}

\end{abstract}

\end{titlepage}

\section{Introduction}
For $\gamma$ with $0 < \gamma \leq \pi$, we denote by $\Gamma_\gamma$ the right
circular cone in $\bR^n$ of angle $\gamma$. We let $\{B_t: t\geq 0\}$ be
$n$-dimensional Brownian motion and denote by $E_x$ and $P_x$ the expectation and
the probability associated with this motion starting at $x$. We write
$T_\gamma=\inf\{t > 0: B_t\not\in\Gamma_\gamma\}$, so that $T_\gamma$ is the first
exit time of the Brownian motion from $\Gamma_\gamma$. In 1977, Burkholder
\cite{bu} found the sharp order of integrability of $T_\gamma$. More precisely, he
showed the existence of a critical constant $p(\gamma,n)$ (given explicitly in
terms of  zeros of certain hypergeometric functions) such that, for
$z\in\Gamma_\gamma$,
\begin{equation}
E_z T_\gamma^{p/2} < \infty,
\label{1.1}
\end{equation}
if and only if $p < p(\gamma,n)$. In particular, for a cone in two dimensions (a
case that had already been solved in \cite{Sp}) $E_z T_\gamma^{p/2} < \infty$ if
and only if $p<\pi/(2\gamma)$. Notice that by making the angle of the cone
arbitrarily small we can make $p$ arbitrarily large.

The cone in two dimensions can be thought of as the domain above
the graph of the function $y=a|x|$. It is then natural to study
the order of integrability of exit times from other unbounded
regions,  in particular parabolas. Since any parabola is contained
in a cone of arbitrarily small angle, the exit time for a parabola
has finite moments of all orders. On the other hand, by comparing
with rectangles, it is also easy to show that the exit time  is
not exponentially integrable. Ba\~nuelos, DeBlassie and Smits
showed in \cite{ba} that if $\tau_{\cal P}$ is the exit time of
the Brownian motion from  the parabola ${\cal P}=\{(x,y)\colon x>
0,\ |y| < A\sqrt{x}\}$ and $z \in {\cal P}$, then there exist
positive constants $A_1$ and $A_2$ such that
\begin{eqnarray}
-A_1&\leq&\liminf_{t\to\infty} t^{-\frac{1}{3}}
\log \left[ P_{z}\{\tau_{\cal P} > t\}\right]\label{1.2}\\
&&\quad\leq\ \,\limsup_{t\to\infty}
t^{-\frac{1}{3}}\log \left[P_{z}\{\tau_{\cal P} > t\}\right]\
\leq\  -A_2.\nonumber
\end{eqnarray}

Thus if $a < 1/3$, then   $E_{z} [\exp(b \tau_{\cal P} ^{a})] <
\infty$ for each $b>0$ and, if $a>{1/3}$, then   $E_{z} [\exp(b
\tau_{\cal P} ^{a})] =\infty$ for each $b>0$. This result was
extended to all dimensions and to other unbounded regions by W.~Li
\cite{li}. More recently, M.~Lifshits and Z.~Shi \cite{lif} found
the limit $-l$ of $ t^{-\frac{1}{3}}\log \left[ P_{z}\{\tau_{\cal
P} > t\}\right]$, as $t\to\infty$. Then, if $b <l$ then $ E_{z}
[\exp (b \tau_{\cal P}^{1/3})] < \infty $ and if $b>l$ then  $
E_{z} [\exp (b \tau_{\cal P}^{1/3})] = \infty $, (it is not known
whether or not $ E_{z} [\exp (l \tau_{\cal P}^{1/3})]$ is finite
or infinite), so that Lifshits and Shi determine the sharp order
of integrability of the exit time $\tau_{\cal P}$ (see the proof
of Theorem 2). Their result holds for the more general parabolic
regions studied in \cite{li} in any dimension. In \cite{van}, M.
van den Berg used the sharp results of Lifshits and Shi to obtain
analogues of these results for the Dirichlet heat kernel. Finally,
similar results have been obtained recently by D. DeBlassie and R.
Smits \cite{DS} for  ``twisted parabolas" in two dimensions.

In the case of the cone $\Gamma_{\gamma}$ in $\bR^n$, the sharp order of
integrability of the exit position $B_{T_{\gamma}}$ is known. For any domain $D$
in $\bR^n$ we write $B^*_{\tau_D}$ for the maximum distance of Brownian motion
from the origin up to the exit time $\tau_D$ for $D$. More precisely, we set
$$
B^*_{\tau_D}=\sup\{|B_t|\colon 0 \leq t <\tau_D\}.
$$
Then, by Burkholder's inequality (Theorem 2.1 in \cite{bu}), for
any finite, positive $p$ there are constants $C^1_{p,n}$ and
$C^2_{p,n}$ such that
\begin{equation}
C^1_{p,n}E_z\left[n\tau_D+|z|^2\right]^{p/2} \leq
E_z[B_{\tau_D}^*]^p \leq
C^2_{p,n}E_z\left[n\tau_D+|z|^2\right]^{p/2}.
\label{1.3}
\end{equation}
Thus it follows from (\ref{1.1}) that
\begin{equation}
E_z\,[B^*_{T_\gamma}]^{p}<\infty
\label{1.4}
\end{equation}
if and only if $p < p(\gamma,n)$.  We always have $|B_{T_\gamma}|\leq
B^*_{T_\gamma}$. On the other hand, if $p > 1$ and
\begin{equation}
E_z |B_{T_\gamma}|^{p} < \infty
\label{1.5}
\end{equation}
then by Doob's maximal inequality (1.4) holds, and hence so does
(1.1). This gives the sharp order of integrability $\pi/2\gamma$ of
$B_{T_\gamma}$ in two dimensions. For  $\bR^n$,  $n\geq 3$, since
$T_\gamma$ is finite a.s., it is a consequence of Burkholder's Theorem 2.2 in \cite{bu} that
$E_z [B^*_{T_\gamma}]^p$ is finite if $E_z |B_{T_\gamma}|^{p}$ is finite and this gives the
sharp order of integrability of $B_{T_\gamma}$ in all dimensions.
 The results of Burkholder
generated considerable interest amongst probabilists and analysts. In particular, M.\ Ess\'en
and K.\ Haliste
\cite{ess} used harmonic measure techniques to obtain some generalizations of Burkholder's
results.

The problem of obtaining the sharp order of integrability of the
exit position $B_{\tau_{\cal P}}$ for the parabola ${\cal P}$, and
of the related random variable $B_{\tau_{\cal P}}^*$, suggest
themselves -- a problem that we address in this paper. As we shall
see below, the order of integrability of $B_{\tau_{\cal P}}^*$ and
$\tau_{\cal P}$ are, as in the case of cones, also closely
related.  We consider, as in \cite{li} and \cite{lif}, more
general regions in $\bR^n$ of the form
\begin{equation}
{\cal P}_\alpha=\{(x,Y)\in\bR\times\bR^{n-1}\colon x > 0,\, |Y| < Ax^\alpha\},
\label{Palpha}
\end{equation}
with  $0 < \alpha < 1$ and $A>0$. (The case $\alpha=1$ is the cone
for which, as described above, we know everything.) We write
$\tau_\alpha$ for the exit time from ${\cal P}_{\alpha}$.

We begin with a very simple result which shows that the random
variables $B_{\tau_{\cal P}}$ and $B_{\tau_{\cal P}}^*$ share the
same integrability properties. More precisely,
\begin{thm}
Suppose that $a$ and $b$ are positive constants and that $z \in \pa$. Then
$E_{z}\left[\exp [ {b\,|{B_{\tau_\alpha}}|^a}]\right]<\infty$ if
and only if $E_{z}\left[\exp
[{b\,({B^*_{\tau_\alpha}})^a}]\right]<\infty$.
\end{thm}

Naturally, we will need to estimate the distribution functions of
$B_{\tau_{\cal P}}$ and $B_{\tau_{\cal P}}^*$ for large $t$ in the
manner of (\ref{1.2}). The probability that the Brownian motion
exits a parabola-shaped domain $\pa$ outside the ball of radius
$t$, that is $P_{z}\{\vert B_{\tau_{\alpha}} \vert > t\}$, is the
harmonic measure of that part of the boundary of $\pa$ lying
outside the ball $B(0,t)$ of radius $t$, taken w.r.t.\ the domain
$\pa$. The larger quantity $P_{z}\{ B^*_{\tau_{\alpha}}  > t\}$ is
the harmonic measure of the intersection of $\pa$ with the sphere
of radius $t$ taken w.r.t.\ the intersection of $\pa$ with the
ball $B(0,t)$. In fact, $1-P_{z}\{ B^*_{\tau_{\alpha}} > t\}$ is
the probability of exiting $\pa$ without ever exiting the ball
$B(0,t)$. These interpretations of the distribution functions of
$B_{\tau_{\alpha}}$ and $B_{\tau_{\alpha}}^*$ facilitate the use
of some well-known and quite precise estimates of harmonic
measure.

Our result for parabola-shaped regions in the plane provides a
complete solution to the problem and we state it separately.

\begin{thm}For the exit position $B_{\tau_{\alpha}}$ from the
parabola-shaped domain $\pa$ in the plane,
\begin{equation}
\lim_{t\to\infty} t^{\alpha-1}\log \left[
       P_{z}\{|B_{\tau_\alpha}|>t\}\right] = -\frac{\pi}{2A(1-\alpha)}.
\label{1.7}
\end{equation}
Furthermore,
\begin{equation}
E_{z}\left[\exp\big[{b\,\vert
B_{\tau_\alpha}\vert^{1-\alpha}}\big]\right]<\infty
\label{1.8}
\end{equation}
if and only if
$$
b < \frac{\pi}{2A(1-\alpha)}.
$$
Then, $E_z\left[\exp[{b\,\vert B_{\tau_\alpha} \vert^{a}}]\right]$
is integrable for each positive $b$ if $ a < 1-\alpha$, and is not
integrable for any positive $b$ if $a> 1-\alpha$.
\end{thm}
\noindent We note that we can determine whether or not $\exp[b\,
\vert B_{\tau_\alpha}\vert^a]$ is integrable in all cases,
including the critical case $a = 1-\alpha$, $b = \pi/\bigl[
2A(1-\alpha) \bigr]$.

\medskip\noindent For parabola-shaped regions in $\bR^n$, it is
proved in \cite{li} that
\begin{eqnarray}
-B_1&\leq& \liminf_{t\to\infty} t^{\frac{\alpha-1}{\alpha+1}}\log
\left[P_{z}\{\tau_{\alpha} > t\}\right]\\
& &\quad\leq\ \, \limsup_{t\to\infty}
t^{\frac{\alpha-1}{\alpha+1}} \log \left[P_{z}\{\tau_{\alpha} >
t\}\right]\ \leq\ -B_2\nonumber
\end{eqnarray}
for two positive constants $B_1$ and $B_2$ depending on $A$, on
$\alpha$ and on the dimension.  In \cite{lif}, the limit of
$t^{\frac{\alpha-1}{\alpha+1}}\log \left[P_{z}\{\tau_\alpha >
t\}\right]$ is shown to exist and its value is determined
explicitly. For comparison purposes we observe that if $n=2$,
$\alpha=1/2$ and $A=1$, the case of the parabola in the plane, the
results in \cite{lif} give
\begin{equation}
\lim_{t\to\infty} t^{-\frac{1}{3}} \log
\left[P_{z}\{\tau_{1/2}>t\}\right]
 = - \frac{3\pi^2}{8},
\end{equation}
while it follows from Theorem 2 that
\begin{equation}
\lim_{t\to\infty} t^{-\frac{1}{2}}\log\left[
P_{z}\{|B_{\tau_{1/2}}|>t\} \right] = - \pi.
\end{equation}

\noindent We stated Theorem 2 in terms of limits of logs of
distributions to draw a parallel with the previously cited work on
exit times. However, our results are  sharper than that, as we
obtain sharp estimates for the distribution itself, (see
Proposition 1 below).

We prove an extension of Theorem 2 to parabola-shaped regions in
higher dimensions.
\begin{thm}We let $\lambda_1$ be the smallest eigenvalue for the Dirichlet
Laplacian in the unit ball of $\bR^{n-1}$. For the exit position
$B_{\tau_{\alpha}}$ from the parabola-shaped region $\pa$ of (\ref{Palpha}), and
for $z \in \pa$,
\begin{equation}
\lim_{t\to\infty} t^{\alpha-1}\log \left[
       P_{z}\{|B_{\tau_\alpha}|>t\}\right] = -\frac{\sqrt{\lambda_1}}{A(1-\alpha)}.
\label{thm4}
\end{equation}
Furthermore,
$$
E_{z}\left[\exp \big[b\,\vert
B_{\tau_\alpha}\vert^{1-\alpha}\big]\right]
$$
is finite if $b<\sqrt{\lambda_1}/\bigl[A(1-\alpha)\bigr]$ and is
infinite if $b>\sqrt{\lambda_1}/\bigl[A(1-\alpha)\bigr]$.
\end{thm}

Our estimates on the distribution function of $B_{\tau_\alpha}$ in
higher dimensions are not sufficiently precise to determine
whether or not $\exp \big[b\,\vert B_{\tau_\alpha}\vert^a\bigr]$
is integrable in the critical case $a= 1-\alpha$, $ b =
\sqrt{\lambda_1}/\bigl[A(1-\alpha)\bigr]$. This is one reason why
we state the two dimensional result separately as Theorem 2.
Moreover, though the method we use to obtain the distribution
estimates for $B_{\tau_\alpha}$ in the planar case are relatively
standard in complex analysis, it forms the general outline of the
method used to obtain lower bounds for the distribution in higher
dimensions. For this reason also, it seems helpful to present the
two dimensional case separately.

A large part of this paper is devoted to adapting  the conformal mapping
techniques that have led to such precise harmonic measure estimates in planar
domains to the parabola-shaped regions ${\cal P}_{\alpha}$ in $\bR^n$. The
so-called  \lq Carleman  method\rq\  can be adapted to estimate harmonic measure
in these domains from above, as we do in Section 4.1. All else being equal,
harmonic measure is generally largest in the most symmetric case, and so one would
expect the upper bounds given by the Carleman  method to be reasonably precise.
However, we were not able to use the Carleman method to obtain the lower bounds for
harmonic measure that we need to determine the exact order of exponential
integrability of the exit position $B_{\tau_\alpha}$. For this, we adapt the
conformal mapping techniques used to prove Theorem 2. Our parabola-shaped domains
being symmetric, we can rewrite the distribution estimates for $B_{\tau_\alpha}$
as a distribution estimate in the corresponding planar parabola-shaped domain, but
at the cost of having to deal with a Bessel-type operator rather than the
Laplacian. Conformal invariance is lost, but the conformal mapping techniques can
still be made to work with considerably more effort. An outline of our method for
obtaining these relatively precise estimates from below for the distribution of
the exit position in parabola-shaped regions can be found at the beginning of
Section 4.2.

Burkholder's inequality (\ref{1.3}), as outlined earlier, allows one to deduce the
integrability properties of the random variable $B_{T_{\alpha}}^*$ for the cone
$\Gamma_\alpha$ from those for the exit time $T_{\alpha}$ for the cone, and vice
versa. In the case of parabola-shaped domains, we can partially pass from
integrability results for the exit position to integrability results for the exit
time and conversely. The proof provides a partial explanation of the connection
between the two critical exponents, $1-\alpha$ in the case of the exit position
and $(1-\alpha)/(1+\alpha)$ in the case of the exit time.
\begin{thm}We may deduce from
\begin{equation}
E_{z} \left[\exp
[{b_1\,({B^*_{\tau_\alpha}})^{1-\alpha}}]\right]<\infty,
\label{1.12}
\end{equation}
with $b_1$ positive, that there is some positive $b_2$ for which
\begin{equation}
E_{z} \left[\exp \bigl[{b_2\,
{\tau_\alpha}^{\frac{1-\alpha}{1+\alpha}}}\bigr]\right]<\infty.
\label{1.13}
\end{equation}
Conversely, we may deduce from (\ref{1.13}), with $b_2$ positive,
that (\ref{1.12}) holds for some positive $b_1$.
\end{thm}

\section{Proofs of Theorems 1 and 4}
\begin{proof}[Proof of Theorem 1] Clearly $|B_{\tau_\alpha}|\leq
B^*_{\tau_\alpha}$ and hence if $E_{z} \left[\exp
[{b\,({B^*_{\tau_\alpha}})^a}]\right]<\infty$, then $E_{z}
\left[\exp[b\,|{B_{\tau_\alpha}}|^a]\right]<\infty$, with the same
$a$ and $b$.

We turn to the proof of the converse. As observed in the
Introduction, $E_{z}\tau_{\alpha}^p<\infty$,  for each finite $p$.
By  (\ref{1.3}) we also have $E_{z}
({B^*_{\tau_\alpha}})^p<\infty$ for $p$ finite. By Doob's maximal
inequality,
$$
E_{z} ({B^*_{\tau_\alpha}})^p \leq \left(\frac{p}{p-1}\right)^p
E_{z} |{B_{\tau_\alpha}}|^p, \ \ {\rm for}\ 1<p<\infty.
$$
Thus there exists a $p_0$ such that for all $p>p_0$,
\begin{equation}
E_{z} ({B^*_{\tau_\alpha}})^p \leq 4
E_{z} |{B_{\tau_\alpha}}|^p.
\end{equation}
We choose an integer $k_0$, depending on $p_0$ and $a$, such that
$p_0<ka$ for all $k>k_0$.
Then
\begin{eqnarray*}
E_{z} \left[\exp [{b\;({B^*_{\tau_\alpha}})^a}]\right] &=
&1+\sum_{k=1}^{k_0} {\frac{b^k}{k!}}E_z
({B^*_{\tau_\alpha}})^{ak} +\sum_{k=k_0+1}^{\infty} {\frac{b^k}{k!}}E_{z} ({B^*_{\tau_\alpha}})^{ak}\\
&\leq& 1+\sum_{k=1}^{k_0} {\frac{b^k}{k!}}E_z ({B^*_{\tau_\alpha}})^{ak} +4\sum_{k=k_0+1}^{\infty}
{\frac{b^k}{k!}}E_{z} |{B_{\tau_\alpha}}|^{ak}\\
&\leq& 1+\sum_{k=1}^{k_0} {\frac{b^k}{k!}}E_{z} ({B^*_{\tau_\alpha}})^{ak}
+4\,E_{z} \left[\exp [{b\;|{B_{\tau_\alpha}}|^a}]\right]\\
&<&\infty,
\end{eqnarray*}
which proves the theorem.
\end{proof}

\begin{proof}[Proof of Theorem 4] We argue as in the proof of the classical
good-$\lambda$ inequalities, see \cite{bu} for example.  Suppose
(\ref{1.13}) holds for some $b_2>0$. We set $\beta=1+\alpha$ and
note that $0<\beta/2<1$. Then
\begin{eqnarray*}
P_{z}\{B^*_{\tau_\alpha}>t \}\ \leq\ P_{z}\{B^*_{\tau_\alpha}>t,
\; \tau_{\alpha}\leq t^{\beta}\} +
P_{z}\{\tau_{\alpha}>t^{\beta}\}.
\end{eqnarray*}
For the second term (\ref{1.13}), together with a Chebyshev style argument, gives
\begin{equation}
P_{z}\{\tau_{\alpha}>t^{\beta}\}\leq C\exp\left[-b_2t^{{\beta}
{(\frac{1-\alpha}{1+\alpha})}}\right]=C\exp\left[-b_2t^{1-\alpha} \right].
\label{2.2}
\end{equation}
Take $t>2|z|$, so that $t-|z|>t/2$. By our assumption on $t$ and
by scaling,
\begin{eqnarray}
P_{z}\{B^*_{\tau_\alpha}>t,  \; \tau_{\alpha}\leq t^{\beta}\}&
\leq
& P_{z}\left\{\sup_{0\leq
s<\tau_{\alpha}}|B_s-z|>\frac{t}{2}, \tau_{\alpha}\leq  t^{\beta}\right\}\label{2.3}\\
&\leq& P_{z} \left\{\sup_{0\leq
s<t^{\beta}}|B_s-z|>\frac{t}{2} \right\}\nonumber\\
&=& P_0 \left\{\sup_{0\leq
s<t^{\beta}}|B_s|>\frac{t}{2}\right\}\nonumber\\
&=& P_0 \left\{\sup_{0\leq
s<1}|B_s|>\frac{1}{2}\,t^{1-\beta/2}\right\}\nonumber\\
&\leq&
C\exp\left[- \frac{1}{8}\,t^{2-\beta}\right] =C\exp\left[-
\frac{1}{8}\,t^{1-\alpha}\right].\nonumber
\end{eqnarray}
In the  last inequality we used the  well-known fact (see
\cite{bu}, inequality (2.15)) that
$$
P_0\left\{\sup_{0\leq s<1}|B_s|>\lambda\right\}\leq
C\exp\left[-\lambda^2/2\right].
$$
It now follows from (\ref{2.2}) and (\ref{2.3}) that
$$
P_z\{ B_{\tau_\alpha}^* > t\} \leq C \exp\left[ -2b_1
t^{1-\alpha}\right],
$$
for all large $t$, where $2b_1 = \min\{1/8, b_2\}$. Since
$$
E_{z} \left[\exp [{b\,({B^*_{\tau_\alpha}})^{1-\alpha}}] \right] =
1 + b\int_0^\infty e^{bt} P_z\{ B_{\tau_\alpha}^* >
t^{1/(1-\alpha)}\}\,dt,
$$
we obtain (\ref{1.12}).

\smallskip Conversely, suppose that (\ref{1.12}) holds, so that
\begin{equation}
P_z\{ B_{\tau_\alpha}^* > t\} \leq C \exp\left[ -b_1
t^{1-\alpha}\right].
\label{tc}
\end{equation}
Let $\beta=1/(1+\alpha)$. We assume, as we may, that $t$ is
very large and set
$$
\tilde{\cal P}(t)={\cal P_{\alpha}}\cap B(0, t^{\beta})
$$
and write $\tilde\tau$ for its exit time. Clearly $\tilde\tau\leq \tau_{\alpha}$ and
\begin{eqnarray}
P_{z}\{\tau_{\alpha}>t\}\ \leq\ P_{z}\{\tau_{\alpha}>t,
\tau_{\alpha}=\tilde\tau\} +P_{z}\{\tau_{\alpha}>t,
\tau_{\alpha}>\tilde\tau\}.
\label{2.4}
\end{eqnarray}
Let us denote the ball centered at $0$ and of radius
$r$ in $\bR^{n-1}$ by $B_{n-1}(0, r)$ and the exit time of Brownian motion from it
 by $\tau_{B_{n-1}(0, r)}$.  We recall that for $t>1$,
$$
P_{0}\{\tau_{B_{n-1}(0, 1)} >t\}\leq C \exp\left[-\lambda_1
t\right],
$$
where $C$ is a constant independent of $t$ and $\lambda_1$ is the
first Dirichlet eigenvalue of  $B_{n-1}(0, 1)$.  By our definition
of the region $\tilde{\cal P}(t)$,
 \begin{eqnarray}
P_{z}\{\tau_{\alpha}>t, \tau_{\alpha}=\tilde\tau\}&\leq &
P_{z}\{\tilde\tau>t\}
\label{2.5}\\
&\leq & P_{0}\{\tau_{B_{n-1}(0, t^{\alpha\beta})}>t\}\nonumber\\
&=& P_0\{\tau_{B_{n-1}(0, 1)}>t^{1-2\alpha\beta}\} \nonumber\\
&\leq&  C\exp\left[-\lambda_1 \, t^{1-2\alpha\beta}\right], \nonumber\\
&=& C\exp\left[-\lambda_1\, t^{\frac{1-\alpha}{1+\alpha}}\right],\nonumber
\end{eqnarray}
where we used scaling for the first  equality above.   On the
other hand, using (\ref{tc}),
 \begin{eqnarray}
P_{z}\{\tau_{\alpha}>t, \tau_{\alpha}>\tilde\tau\}&\leq&
P_{z}\{B^*_{\tau_\alpha}>t^{\beta}\}
\label{2.6}\\
&\leq & C\exp\left[-b_1 t^{\beta(1-\alpha)}\right].\nonumber\\
&=&  C\exp\left[-b_1 t^{\frac{1-\alpha}{1+\alpha}}\right].\nonumber
\end{eqnarray}
The estimates  (\ref{2.4})--(\ref{2.6}) prove that (\ref{1.13})
may be deduced from (\ref{1.12}). This completes the proof of
Theorem 4.
\end{proof}

We note that our argument above proves that if $E_{z} \left[\exp
[b_2\, {\tau_\alpha}^{a}]\right]<\infty$ for some $0<a<1$ and $b_2>0$,
then $E_{z} \left[\exp [b_1\,({B^*_{\tau_\alpha}})^{2a/(1+a)} ]
\right]<\infty$ for some $b_1$ depending on $b_2$. Conversely, if
$E_{z} \left[\exp [b_1\,({B^*_{\tau_\alpha}})^{a}]\right] <\infty$
for some $a$ and $b_1$, then
$$E_{z} \left[\exp [b_2\,
{\tau_\alpha}^{a/(a+2\alpha)}]\right] < \infty$$ for some $b_2$
depending on $b_1$. Thus by Theorems 2 and 3 we see that
\begin{equation}
E_{z} \left[\exp [b\, {\tau_\alpha}^{p}]\right]<\infty
\end{equation}
for some $b>0$ if and only
$$
p\leq {\frac{1-\alpha}{1+\alpha}},
$$
 as was already proved in \cite{ba} for $\alpha=1/2$ and in
$\cite{li}$ for general $\alpha$ in $(0,\,1)$.

\section{Parabola-shaped domains in the plane}

\subsection{Harmonic measure estimates}

For the moment we restrict ourselves to two dimensions. The key
estimate for the distribution function for $B_{\tau_\alpha}$, of
which the limit (\ref{1.7}) is a direct consequence, is
Proposition 1.  As we shall see below, we may suppose without loss
of generality that $z=z_0$, where $z_0$ is the point $(1,0)$.
\begin{prop}There are constants $C_1$ and $C_2$ depending only on $\alpha$
and $A$ such that, as $t \to \infty$,
\begin{eqnarray*}
C_1\exp\left[-\frac{\pi}{2A(1-\alpha)} t^{1-\alpha}\right] &\leq&
P_{z_0}\{|B_{\tau_\alpha}|>t\}\cr
&\leq &\ \,
 C_2\exp\left[-\frac{\pi}{2A(1-\alpha)} t^{1-\alpha} +
{\rm o}(t^{1-\alpha})\right].
\end{eqnarray*}
\end{prop}
The first step in the proof of Proposition 1 is to show that there
is a negligible difference between $P_{z_0}\{ \vert
B_{\tau_\alpha} \vert > t\}$ and $P_{z_0}\{ B_{\tau_\alpha}^1  >
t\}$, where for $z = (x,Y) \in \bR \times \bR^{n-1}$ we write
$z^1$ for $x$, the projection onto the first coordinate. This will
follow from a simple estimate that holds in all dimensions.
\begin{lemma}
For the exit position $B_{\tau_\alpha}$ of Brownian motion from ${\cal P}_\alpha$
in $\bR^n$, and for sufficiently large $t$,
$$
P_{z_0}\{ B_{\tau_\alpha}^1  > t\}\ \leq\ P_{z_0}\{ \vert B_{\tau_\alpha} \vert >
t\}\ \leq\ P_{z_0}\{ B_{\tau_\alpha}^1  > t - A^2 t^{2\alpha-1} \}.
$$
\end{lemma}
\begin{proof}
For fixed $t$ positive, we write $S(0,t)$ for the sphere center
$0$ and radius $t$ and we write $x(t)$ for the common first
coordinate of the points of intersection of the boundary of $\cal
P_{\alpha}$ with $S(0,t)$. Then
\begin{equation}
P_{z_0}\left\{|B_{\tau_\alpha}|>t\right\}=P_{z_0}\{
B_{\tau_\alpha}^1 > x(t)\}.
\label{3.1}
\end{equation}
We also have, for all sufficiently large $t$,
\begin{equation}
t-A^2t^{2\alpha-1} < x(t) < t.
\label{3.2}
\end{equation}
The upper bound is clear. For the lower bound, we observe that any point $(x,Y)$
with $x = t-A^2t^{2\alpha-1}$ and $\vert Y \vert = A(t-A^2t^{2\alpha-1})^\alpha$
lies inside the sphere $S(0,t)$. Indeed,
\begin{eqnarray*}
(t-A^2t^{2\alpha-1})^2&+&A^2(t-A^2 t^{2\alpha-1})^{2\alpha}-t^2\cr
& = &A^4t^{4\alpha-2}-2A^2t^{2\alpha}+A^2(t-A^2 t^{2\alpha-1})^{2\alpha}\cr
&\leq &A^4t^{4\alpha-2}-2A^2t^{2\alpha}+A^2t^{2\alpha}\cr
&=&A^2t^{2\alpha}(A^2t^{2(\alpha-1)}-1)\cr
&<&0,
\end{eqnarray*}
where we use the assumption that $0<\alpha<1$ and note that $0< t-A^2
t^{2\alpha-1} < t$.
\end{proof}
\noindent Since $\left( t - A^2 t^{2\alpha-1} \right)^{1-\alpha} =
t^{1-\alpha}[1+{\rm o}\,(1)]$, Proposition 1 would follow from Lemma 1 together
with
\begin{eqnarray}
\quad C_1\exp\bigg[-\frac{\pi}{2A(1-\alpha)}\
t^{1-\alpha}\bigg]&\leq& P_{z_0} \{B_{\tau_\alpha}^1 > t\}
\label{3.1a}\\
&\leq & \,
C_2\exp\bigg[-\frac{\pi}{2A(1-\alpha)}\
t^{1-\alpha}\bigg].\nonumber
\end{eqnarray}
\subsection{Proof of (3.3) by conformal mapping}
We wish to use the Ahlfors--Warschawski estimates on conformal
mappings. Following the standard notation in this area, we
introduce the functions $\varphi_+(x) = -\varphi_-(x) = Ax^\alpha$
for $x>0$ and $\theta(x) = \varphi_+(x)-\varphi_-(x)$, so that
$\theta(x) = 2Ax^\alpha$ for $x \geq 0$. Then ${\cal P}_{\alpha}$
has the form
$$
{\cal P}_\alpha\ =\ \{z =x+iy\colon x>0 \mbox{ and }\varphi_- (x)
< y < \varphi_+(x)\}.
$$
We write $f(z)$ for the conformal mapping of ${\cal P}_\alpha$
onto the standard strip
$$
S = \{ w\colon \vert {\rm Im}\,w \vert < \pi/2\}
$$
for which $f(z_0) = 0$ and $f'(z_0) >0$. The mapping $f$ is then
symmetric in the real axis. We write $E_t$ for $\partial \pa \cap
\{{\rm Re}\,z>t\}$ and write $F_s$ for $\partial S \cap \{{\rm
Re}\,w > s\}$. For $t>0$, $E_t$ corresponds under the mapping $f$
to $F_s$, for some $s$ that we denote by $s=s(t)$.
%
\bigskip\bigskip\noindent

\newbox\parabola
\figinit{pt}
\figpt1:(40,38)
\figpt2:(100,38)
\figpt3:(94,72)
\pssetupdate{yes}
\figvisu{\parabola}{}{
\figinsert{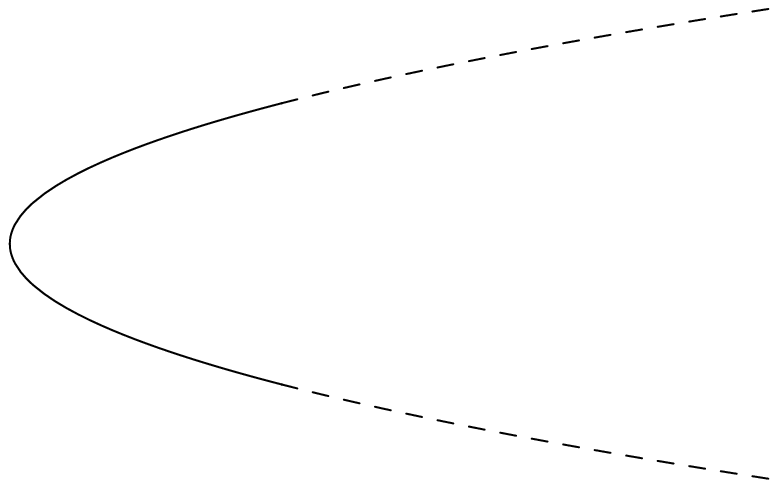,0.5}
\figsetmark{$\bullet$}\figwritene 1:$z_0$(4)
\figwritec[2]{${\cal P}_\alpha$}
\figwritec[3]{$E_t$} }
\newbox\arrows
\figinit{pt} \figpt1:(-30,18) \figpt2:(-10,15) \figpt3:(10,9)
\figpt4:(30,0) \figpt11:(-50,-12) \figpt12:(-30,-20)
\figpt13:(-10,-28) \figpt14:(10,-30) \pssetupdate{yes}
\pssetupdate{yes} \psbeginfig{arrows.ps}
\psarrowBezier[1,2,3,4]  
\psarrowBezier[14,13,12,11]  
\psendfig
\figvisu{\arrows}{}{ \figinsert{arrows.ps,0.7} \figwritenw 3:$f$(5)
\figwritese 12:$g$(9) }
%
%
\newbox\strip
\figinit{pt}
\figpt4:(-90, 30)
\figpt5:(-20,30)
\figpt6:(-90,-30)
\figpt7:(-20,-30)
\figpt14:(-20, 30)
\figpt15:(100,30)
\figpt16:(-20,-30)
\figpt17:(100,-30)
\figpt8:(-60,0)
\figpt9:(-20,0)
\figpt10:(35,31)
\pssetupdate{yes}
\psbeginfig{strip.ps}
\psline[4,5]
\psline[6,7]
\pssetdash{8}
\psline[14,15]
\psline[16,17]
\psendfig
\figvisu{\strip}{}{
\figinsert{strip.ps,0.8}
\figsetmark{$\bullet$}\figwritene 8:$0$(4)
\figwritec[9]{$S$}
\figwritec[10]{$F_s$ $\big(s=s(t)\big)$}}
\centerline{\quad\box\parabola\hfill\vbox{\box\arrows\bigskip}
                                 \hfill\vbox{\box\strip\bigskip}\quad}
\bigskip\bigskip\noindent

\noindent Then, by conformal invariance of harmonic measure, $P_{z_0}
\{B_{\tau_\alpha}^1 > t\}$ is the harmonic measure at 0 of
$F_{s(t)}$ with respect to $S$, that is
\begin{equation}
P_{z_0} \{B_{\tau_\alpha}^1 > t\} = \omega(0,F_{s(t)}; S).
\label{3.2a}
\end{equation}
Harmonic measure in the strip is easy to estimate, for example by
mapping the strip onto the unit disk where harmonic measure at the
origin coincides with normalized angular measure on the circle.
One may show that
$$
\omega(0,F_s; S) = \frac{1}{\pi} \arg \left[ e^{2s} - 1 +
2ie^s\right],
$$
whence, for $s$ large,
\begin{equation}
\frac{1}{\pi}\,e^{-s} \leq \omega(0,F_s; S) \leq \frac{4}{\pi}\,
e^{-s}.
\label{3.2b}
\end{equation}
Warschawski's estimates \cite[Theorem VII]{Wa} for the real part of the mapping
$f$ involve an error term $\int^\infty \theta'(x)^2 / \theta(x)\,dx$ for which, in
our case,
$$
\int_1^\infty \frac{[\,\theta'(x)\,]^2}{\theta(x)}\, dx =
A\int_1^\infty \frac{(2\alpha x^{\alpha-1})^2}{2x^\alpha}\, dx =
2\alpha^2 A\int_1^\infty x^{\alpha-2}\, dx < \infty.
$$
From \cite[Theorem VII]{Wa} we deduce that, for $z_1 = x_1 + iy_1$
and $z_2 = x_2 + iy_2$,
$$
{\rm Re}\, f(z_2) - {\rm Re}\, f(z_1) = \pi \int_{x_1}^{x_2}
\frac{dx}{\theta(x)} + {\rm o}\,(1),
$$
as $x_1$, $x_2 \to +\infty$, uniformly with respect to $y_1$ and
$y_2$. As a consequence we obtain that as $t \to \infty$,
$$
{\rm Re}\, f(t+iy) = \pi \int_1^t \frac{dx}{\theta(x)} + {\rm
O}(1) = \frac{\pi}{2A(1-\alpha)}\,t^{1-\alpha} + {\rm O}(1).
$$
Since this estimate is uniform in $y$, it follows that
\begin{equation}
s(t) = \frac{\pi}{2A(1-\alpha)}\,t^{1-\alpha} + {\rm O}\,(1),
\quad \mbox{as } t \to \infty.
\label{3.2c}
\end{equation}
We will use this estimate on the boundary correspondence under $f$
again in Section 4. Together (\ref{3.2a}), (\ref{3.2b}) and
(\ref{3.2c}) yield (\ref{3.1a}). Thus Proposition 1 is proved.

\bigskip\noindent{\bf Note }It is possible to shortcut the above
explicit calculations by using Haliste's estimates for harmonic
measure \cite[Chapter 1.2]{ha}, that are themselves based on the
Ahlfors-Warschawski approximations. Haliste formulates her
estimates in terms of the harmonic measure of a vertical cross cut
$\theta_t = [t-iAt^\alpha, t+iAt^\alpha]$ with respect to the
truncated domain $\pa(t) = \{ (x,y) \colon 0<x<t,\, \vert y \vert
< Ax^\alpha\}$. We set
$$
B_{\tau_\alpha}^{1,*} = \sup\{|B_t^1|\colon 0 \leq t <\tau_\alpha\}.
$$
Then,
$$
\omega(z_0,\theta_t; {\cal P}_\alpha(t)) = P_{z_0}
\{B_{\tau_\alpha}^{1,*} > t\},
$$
for which Haliste provides estimates similar to (\ref{3.1a}). A
version of Lemma 1 for the maximal functions
$B_{\tau_\alpha}^{1,*}$ and $B_{\tau_\alpha}^*$ yields Proposition
1 with $\vert B_{\tau_\alpha} \vert$ replaced by
$B_{\tau_\alpha}^*$. At this point, it suffices follow the
argument in the next section with $B_{\tau_\alpha}^*$ instead of
$\vert B_{\tau_\alpha} \vert$ and to recall that $\vert
B_{\tau_\alpha} \vert$ and $B_{\tau_\alpha}^*$ have the same
integrability properties (Theorem 1).

\subsection{Proof of Theorem 2}

As above, we assume for the moment that $z=z_0=(1, 0)$. First, the
limit (\ref{1.7}) is a direct consequence of Proposition 1. Next
we observe that
\begin{equation}
E_{z_0}\left[\exp[b\,|B_{\tau_\alpha}|^a]\right]=1 + b
\int_0^\infty e^{bt} P_{z_0} \{|B_{\tau_\alpha}|>t^{1/a}\} dt.
\label{3.7}
\end{equation}
Proposition 1 yields that
\begin{eqnarray*}
C_1\exp\bigg[\bigg( b - \frac{\pi}{2A(1-\alpha)}\bigg)t \bigg]
&\leq&
e^{bt} P_{z_0} \{|B_{\tau_\alpha}|>t^{1/(1-\alpha)}\}\\
&\leq& \, C_2\exp\bigg[\bigg(b - \frac{\pi}{2A(1-\alpha)}\bigg)
t+o(t)\bigg].
\end{eqnarray*}
We suppose first that $ b < \pi/[2A(1-\alpha)] $ and set
$$
\varepsilon = \frac{\pi}{2A(1-\alpha)} - b.
$$
Then $\varepsilon $ is positive and we observe that $-\varepsilon
t+o(t)\leq- \varepsilon t/ 2$ for sufficiently large $t$. Hence,
for $t$ large,
$$
e^{bt} P_{z_0}\{|B_{\tau_\alpha}|>t^{1/(1-\alpha)}\}\leq C_2
e^{-\frac{1}{2}\varepsilon t}.
$$
It follows from the case $a=1-\alpha$ of (\ref{3.7}) that
$E_{z_0}\bigl[ \exp [b\, |B_{\tau_\alpha}|^{1-\alpha} ] \bigr]
<\infty$ in this case.

In the case $b \geq \pi/[2A(1-\alpha)]$,
$$
e^{bt} P_{z_0}\{|B_{\tau_\alpha}|> t^{1/(1-\alpha)}\} \geq C_1
\exp \bigg[ \left(b -\frac{\pi}{2A(1-\alpha)}\right)\ t \bigg]
\geq C_1.
$$
This gives $ E_{z_0}\bigl[ \exp[b\, |B_{\tau_\alpha}|^{1-\alpha} ]
\bigr] =\infty$, for such $b$.

The cases $0 < a < 1-\alpha$, $b > 0$ and $a>1-\alpha$, $b>0$ may
be handled similarly, or one may compare the expected value with
that of $\exp[ b\, \vert B_{\tau_\alpha} \vert^{1-\alpha} ]$
where, say, $b = \pi/[4A(1-\alpha)]$ and $b= \pi/[A(1-\alpha)]$,
respectively.

\smallskip
We now remove the assumption that  $z=(1, 0)$.  We first deal with the upper
bound. We let $z=(x, y)$ and assume that $t$ is very large and certainly much
larger than 1.  By translation of paths it is clear that $\omega((x_1,
y),\theta_t; \pa(t))\leq \omega((x_2, y), \theta_t; \pa(t))$ whenever $x_1\leq
x_2$. Hence, we may assume that $x\geq 1$. Thus for $t>>x>>1$ we have, by
symmetry,
\begin{equation}
\sup_{z'\in \theta_x}\omega(z',\theta_t; \pa(t))=\omega(x,\theta_t; \pa(t)).
\end{equation}
From this and the more general upper bound of Haliste \cite{ha} we have that for
$z=(x, y)$ and  $t>> x>>1$,
\begin{eqnarray}
\omega(z,\theta(t);
\pa(t))&\leq & C_2\exp\left[-\pi\int_x^t \frac{du}{\theta(u)} \right]\nonumber\\
&=& C_2\exp\bigg[-\frac{\pi}{2A(1-\alpha)}\ t^{1-\alpha}-\frac{\pi}{2A(1-\alpha)}\
x^{1-\alpha}\bigg].
\end{eqnarray}
These two inequalities, (3.8) and (3.9), give the desired upper
bound estimate on $P_{z} \{B_{\tau_\alpha}^{1,*} > t\}$.

For the lower bound, we may assume by the above argument that
$z=(x, y)$ with $0<x<1$. As before, we may assume that $t$ is much
greater than $1$. By a standard Whitney chain argument and the
Harnack inequality we have
\begin{equation}
\omega(z,\theta_t;\pa(t))\geq C(z)\, \omega(z_0,\theta_t;\pa(t)),
\end{equation}
where $C(z)$ is a function of $z$ depending on the distance of $z$
to the boundary of ${\cal P}_{\alpha}$. The general lower bound
follows from this and the case of $z=(1, 0)$ which we have already
done. \qed

\section{Parabola-shaped regions in $\bR^n$}

The parabola-shaped regions with which we work have the form
$$
\pa =\{(x,Y)\in \bR\times\bR^{n-1}\colon x
> 0,\; \vert Y \vert  < Ax^\alpha\},
$$
for $0<\alpha<1$ and $A>0$. Our objective is to derive estimates for the
distribution function of the exit position of Brownian motion from such regions
or, equivalently, for the harmonic measure of the exterior of the ball of center
$0$ and radius $t$ with respect to such regions. In Section 4.1, we derive an
upper bound for the distribution function by means of the Carleman method and in
Section 4.2 we introduce a new conformal mapping technique to derive an equally
sharp lower bound.

\subsection{Carleman method: Upper bound for harmonic measure}

It will be more convenient to write the domain $\pa$ in this
section as
$$
\pa = \{(x,Y)\in \bR\times\bR^{n-1} \colon x>0,\ Y\in
B_{n-1}(0,Ax^\alpha)\},
$$
where $B_{n-1}(0,r)$ is the ball in $\bR^{n-1}$ centered at 0 and
of radius $r$. For convenience of notation we set
$\theta(x)=B_{n-1}(0,Ax^\alpha)$ and refer to $\theta(x)$ as a
cross cut of $\pa$ at $x$. For $t$ large we set
$$
\pa(t)=\{(x,Y)\colon 0<x<t, \, Y\in\theta(x)\}.
$$
This is the domain $\pa$ truncated to the right of $t$. For such a
$t$ and any $(x,Y)\in \pa(t)$, we denote by
$\omega\bigl((x,Y),\theta(t); \pa(t)\bigr)$ the harmonic measure
of $(t,0)+\theta(t)$ at the point $(x, Y)$  relative to the domain
$\pa(t)$. We wish to estimate $\omega\bigl((1, 0), \theta(t);
\pa(t)\bigr)$.
\begin{prop} There exit two constants $C_1$ and $C_2$, that depend on
$n$, $\lambda_1$, $A$ and $\alpha$, such that for $t>C_1$,
$$
\omega\bigl((1,0),  \theta(t); \pa(t)\bigr) \leq C_2\,
t^{\alpha(n-1)/2}\
\exp\bigg[-\frac{\sqrt{\lambda_1}}{A(1-\alpha)}\
t^{1-\alpha}\bigg]
$$
\end{prop}
Before we begin the proof of Proposition 2, we show how it leads to an estimate
for the distribution of $B_{\tau_\alpha}$. With the notation of Lemma 1 and of the
note at the end of Section 3.2,
$$
P_{(1,0)}\left\{|B_{\tau_\alpha}|>t\right\}=P_{(1,0)}\{
B_{\tau_\alpha}^1 > x(t)\} \leq P_{(1,0)}\{ B_{\tau_\alpha}^{1,*}
> x(t)\}.
$$
The distribution function for $B_{\tau_\alpha}^{1,*}$ is the harmonic measure that
is estimated in Proposition 2, so that
$$
P_{(1,0)}\{ B_{\tau_\alpha}^{1,*}
> t\} =  \omega\big((1,0), \theta(t); \pa(t)\big).
$$
From the estimate (\ref{3.2}), namely $t-A^2 t^{2\alpha-1} < x(t) < t$, it
follows, as in the proof of Lemma 1, that $x(t) ^{1-\alpha} = t^{1-\alpha}[1+{\rm
o}(1)]$.  Together with Proposition 2, this leads to
\begin{eqnarray*}
P_{(1,0)}\{ B_{\tau_\alpha}^{1,*}> x(t)\}&\leq& C_2\,
x(t)^{\alpha(n-1)/2}\
\exp\bigg[-\frac{\sqrt{\lambda_1}}{A(1-\alpha)}\
x(t)^{1-\alpha}\bigg]\cr
& \leq &C_2\,
t^{\alpha(n-1)/2}\
\exp\bigg[-\frac{\sqrt{\lambda_1}}{A(1-\alpha)}\
t^{1-\alpha}[1+{\rm o}(1)]\bigg]
\end{eqnarray*}
We may absorb the term $C_2\, t^{\alpha(n-1)/2}$ into the ${\rm
o}(1)$ term in the exponential to deduce an estimate for the
distribution function for $B_{\tau_\alpha}$ in the following form.
\begin{prop}Suppose that $\epsilon$ is positive. There
exists a constant $C_1$ depending on $\epsilon$, $n$, $\lambda_1$,
$A$ and $\alpha$ such that, for $t>C_1$,
$$
P_{(1,0)}\left\{|B_{\tau_\alpha}|>t\right\} \leq
\exp\bigg[-\frac{\sqrt{\lambda_1}}{A(1-\alpha)}\,[1-\epsilon]\,
t^{1-\alpha}\bigg]
$$
\end{prop}

\begin{proof}[Proof of Proposition 2] Our estimates follow those of Haliste \cite{ha}. We
take $t$ to have some large, fixed value and set
\begin{equation}
h(x)=\int_{\theta(x)}\omega^2(x,Y)\,dY, \quad 0<x<t,
\end{equation}
where for convenience we write $\omega(x,Y)$ for
$\omega((x,Y),\theta(t); \pa(t))$. Differentiating $h$ (see
Haliste \cite{ha} for the justification of this step), we obtain
\begin{equation}
h'(x)=\int_{\theta(x)} 2\,\omega_x (x,Y)\,\omega(x,Y)\,dY
\label{4.2}
\end{equation}
and
\begin{equation}
h''(x) = 2 \int_{\theta(x)} \omega_{xx}(x,Y)\,\omega(x,Y)\,dY + 2
\int_{\theta(x)} \vert \omega_x(x,Y)\vert^2\, dY
\end{equation}
We observe that, since $\omega(x,Y)$ is increasing in $x$ for each $Y$, the
derivative $\omega_x(x,Y)$ is non negative and hence $h'(x) \geq 0$. Since
$\omega(x,Y)$ is harmonic, we have
$$
\omega_{xx}(x,Y)+\Delta_Y\, \omega(x,Y) = 0.
$$
Thus,
\begin{equation}
h''(x) = - 2 \int_{\theta(x)}\Delta_Y\,\omega(x,Y)\,\omega(x,Y)\,
dY + 2 \int_{\theta(x)}\omega_x(x,Y)^2\, dY.
\end{equation}
Since the harmonic measure vanishes on the lateral boundary of the
domain, $\omega(x,Y)=0$ if $Y\in\partial\theta(x)$, with $0<x<t$.
Thus, integrating by parts, we obtain
\begin{equation}
h''(x)=2\int_{\theta(x)}\vert\nabla_Y\,\omega(x,Y)\vert^2\, dY +
2\int_{\theta(x)} \omega_x(x,Y)^2\, dY. \label{4.5}
\end{equation}
Writing $B(0,r)$ for $B_{n-1}(0,r)$, we now recall that for all
$u$ that are differentiable on $B(0,r)$ and vanish on $\partial
B(0,r)$,
\begin{equation}
\lambda_{B(0,r)}\leq \frac{\int_{B(0,r)}\vert \nabla
u\vert^2}{\int_{B(0,r)}\vert u\vert^2 },
\label{4.6}
\end{equation}
where $\lambda_{B(0,r)}$ is the first eigenvalue of $B(0,r)$ for
the Laplacian. By scaling,
$$
\lambda_{B(0,r)}= \frac{1}{r^2}\lambda_1,
$$
where $\lambda_1$ is the eigenvalue of the unit ball. In our case
this gives
\begin{equation}
\lambda_{\theta(x)}=\frac{1}{A^2 x^{2\alpha}}\ \lambda_1.
\label{4.7}
\end{equation}
From  (\ref{4.5}) and (\ref{4.6}) we deduce that
\begin{eqnarray}
h''(x) & \geq & 2 \lambda_{\theta(x)}
\int_{\theta(x)}\omega(x,Y)^2\, dY +
2 \int_{\theta(x)} \omega_x(x,Y)^2\, dY \nonumber\\
& = & 2\lambda_{\theta(x)} h(x) + 2\int_{\theta(x)} \omega_x(x,Y)^2\, dY\label{4.8}.
\end{eqnarray}
On the other hand, by (\ref{4.2}) and H\"older's inequality,
$$
h'(x)\leq 2\, \bigg( \int_{\theta(x)} \omega_x(x,Y)^2
\,dY\bigg)^{1/2} \bigg(\int_{\theta(x)} \omega(x,Y)^2\,
dY\bigg)^{1/2}
$$
or
$$
h'(x)^2 \leq 4\bigg( \int_{\theta(x)} \omega_x(x,Y)^2\, dY \bigg)
h(x)
$$
or
$$
\frac{h'(x)^2}{2 h(x)} \leq 2\int_{\theta(x)} \omega_x(x,Y)^2\,
dY.
$$
This and (\ref{4.8}) give
\begin{equation}
h''(x)\geq 2\lambda_{\theta(x)} h(x) + \frac{h'(x)^2}{2h(x)}.
\end{equation}
Since $ 2\sqrt{a}\sqrt{b}\leq a+b $, we conclude that
$$
h''(x) \geq 2 \sqrt{\lambda_{\theta(x)}} \, h'(x)
$$
which, by (\ref{4.7}), is the same as
$$
\frac{h''(x)}{h'(x)}\geq \frac{2\sqrt{\lambda_1}}{A}\, \frac{1}{x^\alpha}.
$$
Following Haliste, we consider the function $g$ on the interval
$(0,t)$ given by
\begin{eqnarray*}
g(x) &=& \int_0^x\exp\bigg(2\int_0^s\ \sqrt{\lambda_{\theta(r)}}\, dr\bigg)\, ds \\
& = & \int_0^x\exp\bigg(\frac{2\sqrt{\lambda_1}}{A}\,  \int_0^s\ \frac{dr}{r^{\alpha}}
\bigg)\, ds\\
&=& \int_0^x\exp\bigg(\frac{2\sqrt{\lambda_1}}{A(1-\alpha)}\, s^{1-\alpha}\bigg)\,
ds.
\end{eqnarray*}
This function satisfies
\begin{eqnarray*}
g'(x) & = & \exp\bigg(2\int_0^x\sqrt{\lambda_{\theta(r)}}\, dr\bigg)\\
g''(x) & = & g'(x)\, 2\, \sqrt{\lambda_{\theta(x)}}\\
&=&g'(x) \frac{2\sqrt{\lambda_1}}{A}\ \frac{1}{x^\alpha},
\end{eqnarray*}
and so
$$
\frac{d}{dx}(\log h'-\log g')\geq 0.
$$
From this it follows that the function $\frac{h'(x)}{g'(x)}$ is non--decreasing on $(0, t)$. 
Since $g(0)=h(0)=0$, the generalized mean value theorem gives that for any $0<x<t$ there is
a $\xi\in (0, x)$ such that 
$\frac{h(x)}{g(x)}=\frac{h'(\xi)}{g'(\xi)}$.  Hence,  

$$\frac{h(x)}{g(x)}\leq\frac{h'(x)}{g'(x)}.$$
Since $g'(x)\ge 0$, this shows that $(\frac{h}{g})'(x)$ is nonnegative and hence the function $\frac{h(x)}{g(x)}$
is non--decreasing. 
Thus, 
\begin{equation}
h(x)\leq h(\xi)\, \frac{g(x)}{g(\xi)}, \ \ {\rm for}\  0<x<\xi<t.
\label{4.10}
\end{equation}
Setting $\mu(r)=2\sqrt{\lambda_{\theta(r)}}$,
\begin{eqnarray}
g(x)&= &\exp\bigg(\int_0^x\mu(r)\, dr\bigg)\int_0^x\exp\bigg(-\int_s^x\ \mu(r)\, dr\bigg)\, ds\nonumber\\
&=&\exp\bigg(\int_0^x\mu(r)\, dr\bigg)\,  H(x).
\label{4.11}
\end{eqnarray}
We may estimate $H(x)$ by
\begin{eqnarray*}
H(x) &=& \int_0^x\exp\bigg(-\frac{2\sqrt{\lambda_1}}{A(1-\alpha)}
\big[x^{1-\alpha}-s^{1-\alpha}\big]\bigg)\, ds\\
&= &\exp\bigg[-\frac{2\sqrt{\lambda_1}}{A(1-\alpha)}\ x^{1-\alpha}\bigg]\int_0^x\
\exp\bigg[\frac{2\sqrt{\lambda_1}}{A(1-\alpha)} s^{1-\alpha} \bigg]\, ds\\
&\leq & x.
\end{eqnarray*}
Therefore,
\begin{equation}
g(x)\leq x\exp\bigg(\int_0^x\mu(r)\, dr\bigg),\quad 0<x<t.
\label{4.12}
\end{equation}
On the other hand, setting
$$
x_0 = x_0(\lambda_1, A, \alpha)=\bigg[
\frac{A(1-\alpha)}{2\sqrt{\lambda_1}}\ln 2 +
1\bigg]^{1/(1-\alpha)} \qquad {\rm and }\qquad K =
\frac{2\sqrt{\lambda_1}}{A(1-\alpha)},
$$
we find that, for $x \geq x_0$,
\begin{eqnarray*}
H(x) &=& \exp\big[- K\, x^{1-\alpha}\big] \int_0^x
\exp\big[K\, s^{1-\alpha} \big]\, ds\\
&= & \exp\big[-K\, x^{1-\alpha} \big]
\int_0^{x^{1-\alpha}}\frac{1}{1-\alpha}\,
r^{{\alpha}/(1-\alpha)}\, \exp\big[K\, r \big]\, dr\\
&\geq & \frac{1}{1-\alpha}\, \exp\big[-K\,
x^{1-\alpha} \big] \int_1^{x^{1-\alpha}}\
\exp\big[K\, r\big]\, dr\\
&=&\frac{1}{K(1-\alpha)} \exp\big[-K\,x^{1-\alpha} \big]
\big[\exp\big[K\, x^{1-\alpha}\big]
- \exp K\big]\\
&=&\frac{1}{K(1-\alpha)} \bigg[1-\exp
\big[K\, (1-x^{1-\alpha})\big]\bigg]\\
&\geq &\frac{1}{2K(1-\alpha)}.
\end{eqnarray*}
We have shown that,
$$
H(x)\geq \frac{A}{4\sqrt{\lambda_1}}\ \ {\rm for}\ \ x \geq x_0,
$$
and this, together with (\ref{4.11}), gives
\begin{equation}
g(x)\geq \frac{A}{4\sqrt{\lambda_1}}\, \exp\bigg(\int_0^x\mu (r)
\,dr\bigg), \ \ {\rm for }\ \ x_0 \leq x < t.
\label{4.13}
\end{equation}
From (\ref{4.10}), (\ref{4.12}) and (\ref{4.13}) we deduce that,
$$
h(x)\leq \frac{4 \sqrt{\lambda_1}}{A}\, x \,
h(\xi)\,\exp\bigg(-\int_x^\xi \mu(r)\, dr\bigg), \ \ {\rm for }\
x_0 \leq x < \xi < t.
$$
Taking $x=x_0$ and letting $\xi$ tend to $t$, we arrive at
\begin{eqnarray*}
h(x_0)&\leq & \frac{4 \sqrt{\lambda_1}}{A}\, x_0 h(t)\, \exp\bigg(-\int_{x_0}^t\mu(r)\, dr\bigg)\\
&= & C_1(\lambda_1, A, \alpha)\, h(t)\, \exp\bigg[-\frac{2\sqrt{\lambda_1}}{A(1-\alpha)}\
t^{1-\alpha}\bigg],
\end{eqnarray*}
for an appropriate constant $C_1$. By our definition of $h$,
$$
h(t)\leq\text{ vol}\big(B_{n-1}(0,At^\alpha)\big)=\gamma_n A^{n-1}
t^{\alpha(n-1)},
$$
where $\gamma_n$ is the volume of the unit ball in $\bR^{n-1}$.
Thus,
$$
h(x_0)\leq C_2(n,\lambda_1,A,\alpha)\, t^{\alpha(n-1)}\,
\exp\bigg[-\frac{2\sqrt{\lambda_1}}{A(1-\alpha)}\
t^{1-\alpha}\bigg],
$$
for all sufficiently large $t$.
We now choose an $r$, independent of $t$, such that the ball
$B_n\big((x_0,0),2r\big)$ is contained in $\pa$; clearly such an
$r$ exists. Then, by the Harnack inequality,
$$
\omega(x_0,0)\leq C_3(n)\, \omega(x,Y)\ \ {\rm for\ all}\ (x,Y)\in
B_n\big((x_0,0),r\big).
$$
Squaring and then integrating over the ball $B_{n-1}(0,r)$ leads
to
\begin{eqnarray*}
\omega^2 (x_0,0) &\leq & C_4(n,r)\,\int_{B_{n-1}(0,r)}\omega^2
(x_0,Y)\, dY\cr
 &\leq& C_4(n,r)
\int_{\theta(x_0)}\omega^2 (x_0,Y)\, dY = C_4(n, r)\, h(x_0).
\end{eqnarray*}
From this we  finally obtain
$$
\omega((x_0,0), \theta(t); \pa(t)) \leq C_5(n,\lambda_1,A,\alpha)
\,t^{\alpha(n-1)/2}\, \exp\bigg[-\frac{\sqrt{\lambda_1}}{A(1-\alpha)}\
t^{1-\alpha}\bigg],
$$
for $t>x_0$. One final application of the Harnack inequality to
move from $(x_0,0)$ to $(1,0)$, and the proposition is proved.
\end{proof}

\subsection{Conformal mapping method: Lower bound}
In this section we write $\pa^n$ for the region in (\ref{Palpha})
to emphasize the dimension. We observe that, because of the
cylindrical symmetry of ${\cal P}^n_\alpha$ (this region is
invariant under rotation about the $x$-axis) and because of the
symmetry of the boundary values of the harmonic measure, the value
of the harmonic measure at $(x,Y)$ in ${\cal P}^n_\alpha$ depends
only on $x$ and on $\vert Y \vert$.

We associate with $\pa^n$ the corresponding domain ${\cal P}_\alpha = {\cal
P}^2_\alpha$ in two dimensions. The technique we develop to obtain lower bounds
for the distribution function of the exit position in ${\cal P}_\alpha^n$ is the
following. We replace Laplace's equation in $\pa^n$ by the corresponding
Bessel-type partial differential equation in $\pa$ -- in the case $n=2$, this
reduces to the Laplacian. Mirroring the arguments in Section 3.2, we map the
parabola-shaped domain $\pa$ conformally onto the standard strip $S = \{w: \vert
{\rm Im}\,w \vert < \pi/2\}$ and determine the form the partial differential
equation takes in the strip after this change of variable. Adapting the
Ahlfors-Warschawski conformal mapping estimates to our purposes, and in some
instances refining them, we show that the partial differential equation in the
strip is but a small perturbation of the Bessel-type partial differential equation
that we started with -- it is almost conformally invariant. The solution of the
unperturbed Bessel-type p.d.e.\ in the strip may be easily estimated. We show how
to compare the solutions of the perturbed p.d.e.\ with those of the unperturbed
p.d.e., so as to obtain the estimates we require on the original harmonic measure
in the parabola-shaped region in $\bR^n$. As well as keeping track of how the
p.d.e.\ changes as we change from one domain to another, we also need to keep
track of the boundary conditions -- but here Warschawski's detailed conformal
mapping estimates are exactly what we need.

We break the proof into a number of subsections and lemmas.

\subsubsection{From a parabola-shaped region in
$\bR^n$ to a planar strip}

To begin with we compute how the Laplace operator changes as we
drop down from $n$ dimensions to $2$ dimensions.

\begin{lemma}Suppose that $H(x,Y)$ is a $C^2$-function on ${\cal
P}^n_\alpha$ that is invariant under rotation about the $x$-axis,
so that $H$ depends only on  $x$ and on $\vert Y \vert$. We
associate with $H$ a function $h(z)$ in the half parabola-shaped
planar domain
$$
{\cal P}_\alpha^{+}\ =\ \{z = x+iy \colon x > 0\ {\rm and }\ 0 < y
< Ax^\alpha\},
$$
defined by
$$
h(x+iy) = H(x,Y) \quad \mbox{ whenever }\ \vert Y \vert = y.
$$
Then
\begin{equation}
\Delta H(x,Y) = \Delta h(x+iy) + (n-2) \frac{h_y(x+iy)}{y}.
\label{1001}
\end{equation}
\label{Hh}
\end{lemma}

\begin{proof}
In this proof, we denote a point in ${\cal P}^n_\alpha$ by
$(x_1,\, x_2, \ldots, x_n)$. With this notation,
$$
H(x_1,\, x_2, \ldots, x_n) = h\left( x_1,\sqrt{x_2^2 + \ldots+
x_n^2}\right) = h(x,y)
$$
with $x=x_1$ and $y =\sqrt{x_2^2 + \ldots+ x_n^2}$. Then,
\begin{eqnarray*}
\Delta H & = & \sum_{i=1}^n \frac{\partial^2}{\partial
x_i^2}
      h(x_1,\sqrt{x_2^2 + \ldots+ x_n^2})   \cr
& = &\frac{\partial^2 h}{\partial x^2} + \sum_{i=2}^n
        \frac{\partial}{\partial x_i}\left[
          \frac{\partial h}{\partial y} \frac{x_i}{\sqrt{x_2^2 + \ldots+ x_n^2}}
          \right]\cr
& = & \frac{\partial^2 h}{\partial x^2} + \sum_{i=2}^n
        \left[ \frac{\partial h}{\partial y}\frac{(x_2^2 + \ldots+ x_n^2) -
        x_i^2}{(x_2^2 + \ldots+ x_n^2)^{3/2}} +
        \frac{\partial^2 h}{\partial y^2}
        \frac{x_i^2}{x_2^2 + \ldots+ x_n^2}
        \right] \cr
& = & \frac{\partial^2 h}{\partial x^2} + \frac{\partial h}{\partial y}
        \left[ \frac{n-1}{\sqrt{x_2^2 + \ldots+ x_n^2}} -
        \frac{1}{\sqrt{x_2^2 + \ldots+ x_n^2}}\right]
        + \frac{\partial^2 h}{\partial y^2} \cr
& = &\Delta h(z) + (n-2) \frac{h_y(z)}{y}.
\end{eqnarray*}
\end{proof}
From  now on we may work in two dimensions and have conformal
mapping at our disposal, at the expense of having to deal with the
more complicated Bessel-type differential operator appearing in
(\ref{1001}), rather than the Laplacian. The complication arises
because this operator is not conformally invariant.

As in Section 3.2, we denote by $w = f(z)$ the conformal mapping
from the domain $\pa$ onto the standard strip $S$, for which $f(1)
= 0$ and $f'(1) > 0$. Since $f$ is real on the real axis, the
upper half ${\cal P}_\alpha^+$ of the parabola-shaped domain
${\cal P}_\alpha$ is mapped to the upper half $S^+$ of the strip,
specifically, $S^+ = \{ w \colon 0< {\rm Im}\, w < \pi/2\}$.

We denote the inverse mapping of $f(z)$ by $g(w)$. We associate a
function $k(w)$ in $S^+$ with a function $h(z)$ in ${\cal
P}_\alpha^+$ according to
\begin{equation}
k(w) = h(g(w)), \quad \mbox{for}\ w \in S^+.
\label{1002}
\end{equation}
Then $h(z) = k(f(z))$, for $z \in {\cal P}_\alpha^+$. In the next
lemma, we compute how the differential operator on the right of
(\ref{1001}) changes under this change of variables.

\begin{lemma}Suppose that $h(z)$ is a $C^2$-function in the
domain ${\cal P}_\alpha^+$ and that $k$ is defined in the strip
$S^+$ by (\ref{1002}). Then, with $g(w) = z$,
\begin{equation}
\Delta h(z) + (n-2) \frac{h_y(z)}{y} = \frac{\Delta k(w)}{\vert
g'(w)\vert^2} -2 (n-2)\frac{{\rm Im}\, \left[ k_w(w)
\,/\,g'(w)\right]}{{\rm Im}\,\left[ g(w) \right]}. \label{1003}
\end{equation}
\label{pdefork}
\end{lemma}
\begin{proof}
We recall that $ h(z) = k(f(z))$, and use the formulae for change
of variable in the complex partial derivatives $\partial/\partial
z$ and $\partial/\partial\overline{z}$, as explained, for example,
in Krantz \cite[Section 1.2]{kr}. First,
$$
\Delta h(z) = (\Delta k)(f(z))\, \vert f'(z) \vert^2, \qquad z \in {\cal
P}_{\alpha}.
$$
In general,
$$
\frac{\partial }{\partial ({\rm Im}\,z)} =i \left(\frac{\partial
}{\partial z} -
                \frac{\partial }{\partial\overline{z}}
\right),
$$
whence {\openup 4pt
\begin{eqnarray*}
h_y(z) & = & i \left(\frac{\partial h}{\partial z} -
                \frac{\partial h}{\partial\overline{z}}\right) \cr
& = & i \left(\frac{\partial }{\partial z}[k(f(z))] -
                         \frac{\partial }{\partial\overline{z}}[k(f(z))]\right)\cr
& = & i \left(\frac{\partial k}{\partial w}(f(z))\, f'(z) -
                         \frac{\partial k }{\partial\overline{w}}(f(z))\,
                         \overline{f'(z)}\right).
\end{eqnarray*}

}\noindent Since $k(w)$ is real-valued, $k_{\overline w} = \overline{k_w}$. We
therefore obtain,
$$
h_y(z) = -2\,{\rm Im}\left[ k_w\bigl(f(z)\bigr)\, f'(z) \right].
$$
Thus,
$$\displaylines{
\Delta h(z) + (n-2) \frac{h_y(z)}{y}\hfill\cr
\qquad\qquad\qquad =  (\Delta k)(f(z))\, \vert f'(z)
\vert^2 - 2 (n-2)
     \frac{{\rm Im}\left[ k_w\bigl(f(z)\bigr)\, f'(z) \right]}{{\rm Im}\,\left[ g(w) \right]}.\hfill\cr
}$$
On substituting $w$ for $f(z)$ and $1/g'(w)$ for $f'(z)$, we obtain (\ref{1003}).
\end{proof}

\subsubsection{Asymptotic estimates for the conformal mapping $g$}

The success  of the transformations introduced in the previous subsection depends
on being able to simplify the expression on the right hand side of (\ref{1003}),
which is essentially an expression for the Laplacian in the domain ${\cal
P}_\alpha^n$ in $\bR^n$ transformed to the strip $S^+$ in the plane. This is
achieved using modifications of results of Warschawski \cite{Wa}, which give
asymptotic expressions for the conformal mapping $g$ and for its derivative. For
some of these we draw on \cite{CH}. We begin with an estimate for the imaginary
part of the mapping $g$. In this context it is helpful to keep in mind that
estimates for the real part of the mapping $g$, that is for the rate of growth of
$g$, are generally more difficult.

We will adopt the more general situation, as considered by Warschawski, of a
conformal mapping $f$ of a domain of the form
$$
D = \{z: \vert {\rm Im}\, z \vert < \phi({\rm Re}\,z) \},
$$
where $\phi(x)$ is continuous on the real line, onto the strip
$S$. In our case, $\phi(x) = Ax^{\alpha}$, for $x$ positive.
Warschawski's domains are not necessarily symmetric, but the
symmetric case is sufficiently general for our purposes here.
Warschawski writes $\theta(x) = 2\phi(x)$ for the width of the
domain $D$ at $x$. We assume that $D$ has {\sl boundary
inclination 0 at $x=\infty$\/} in that
$$
\phi(x_2) - \phi(x_1) = {\rm o}\,(x_2 - x_1) \quad {\rm as }\ x_1,\,x_2 \to
\infty.
$$
This is the case if $\phi(x)$ is continuously differentiable and
$\phi'(x) \to 0$ as $x \to \infty$, and includes our
parabola-shaped domains ${\cal P}_{\alpha}$.

We take the conformal mapping $f$ of $D$ onto $S$ for which $f$ is real and $f'$
is positive on the real axis.  We denote the inverse mapping of $f$ by $g$, in
agreement with earlier notation. We begin with an estimate for the derivative of
$g$.

\begin{lemalp}{\bf Warschawski, \cite[Theorem X(ii)]{Wa}} For each
$\rho$ with $0 < \rho < 1$, the conformal mapping $g$ of the strip $S$ onto $D$,
for which $g$ is real and $g'$ is positive on the real axis, satisfies
\begin{equation}
g'(w) = \left[ \frac{1}{\pi} + {\rm o}(1) \right] \theta( {\rm
Re}\,g(w) ),
\label{1005}
\end{equation}
uniformly as ${\rm Re}\,w \to \infty$ in the sub strip $S_\rho = \{ w: \vert
{\rm Im}\,w \vert < \rho \pi/2\}$.
\label{lemmaWar}
\end{lemalp}

In \cite[Theorem X(iii)]{Wa}, Warschawski obtains an asymptotic
expression for ${\rm Im}\,g(w)$ as ${\rm Re}\,w \to \infty$.
Taking advantage of the symmetry of the domain $D$ (so that ${\rm
Im}\,g(w)= 0$ when ${\rm Im}\,w = 0$) and adapting Warschawski's
proof, we prove

\begin{lemma}The conformal mapping $g$ of the strip $S$ onto $D$, for
which $g$ is real and $g'$ is positive on the real axis, satisfies
\begin{equation}
{\rm Im}\,g(w) = \left[ \frac{1}{\pi} + {\rm o}\,(1) \right]
\theta(g({\rm Re}\,w))\,{\rm Im}\,w,
\label{1006}
\end{equation}
as ${\rm Re}\,w \to \infty$ with $w \in S$. \label{imagpartg}
\end{lemma}

\begin{proof}By symmetry of the mapping $g$, it is enough to prove (\ref{1006})
in the case of ${\rm Im}\, w$ positive. First we show that
(\ref{1006}) holds in each sub strip $S_\rho$, ($0 < \rho < 1$).
By Theorem II(a) in \cite{Wa},
\begin{equation}
\lim_{u \to \infty} \arg g'(u+iv) = 0, \label{1007}
\end{equation}
uniformly in $v$, $\vert v \vert < \pi/2$. Combining (\ref{1007})
with Lemma \ref{lemmaWar}, we find that
$$
{\rm Re}\,g'(w) = \left[ \frac{1}{\pi} + {\rm o}\,(1)\right]
\theta( {\rm Re}\,g(w) ),
$$
for $w \in S_\rho$. Then, for $w = u+iv_0 \in S_\rho$,
\begin{eqnarray*}
{\rm Im}\,g(u+iv_0)
& = & {\rm Im}\,g(u+iv_0) - {\rm Im}\, g(u)\cr
& = & v_0\left.\left( \frac{\partial}{\partial v} {\rm Im}\,g(u+iv)
        \right)\right\vert_{v=v_1} \qquad \left[v_1\in (0,v_0)\right]\cr
& = & v_0 \,{\rm Re}\,g'(u+iv_1)\cr
& = & v_0\left[ \frac{1}{\pi} + {\rm o}\,(1)\right]
        \theta( {\rm Re}\,g(u+iv_1) ).
\end{eqnarray*}
We assert that, uniformly in $v$, $\vert v \vert < \pi/2$,
\begin{equation}
\theta( {\rm Re}\,g(u+iv) ) = \theta( g(u) ) + {\rm o}(1),
\label{1008}
\end{equation}
(which is why we need not distinguish between $\theta( g({\rm
Re}\,w) )$ and $\theta( {\rm Re}\,g(w) )$, up to ${\rm o}\,(1)$).
This then gives the stated expression for ${\rm Im}\,g(u+iv_0)$,
for $\vert v_0 \vert < \rho\pi/2$. To prove (\ref{1008}) we note
that, as at the bottom of Page 290 of \cite{Wa} and as a
consequence of the assumption that $D$ has boundary inclination
$0$ at $x = \infty$,
$$
{\rm Re}\,g(u+iv)  =  g(u)  + {\rm o}(1),\ \ \mbox{as } u \to
\infty.
$$
On using yet again the assumption that $D$ has boundary
inclination $0$ at $x = \infty$, (\ref{1008}) follows.

Finally, we need to show that (\ref{1006}) holds uniformly in $v$.
Given $\epsilon$ small and positive, we take $\rho = 1 -
\epsilon$, so that (\ref{1006}) holds for $w = u + iv$, with
$\vert v \vert \leq (1 - \epsilon)\pi/2$. In particular,
\begin{eqnarray*}
{\rm Im}\,g\left(u+i\frac{\pi}{2}(1 - \epsilon)\right)& = &
\left[\frac{1}{\pi} + {\rm o}\,(1)\right] \theta(g(u))
\frac{\pi}{2}(1 - \epsilon) \cr
& = & (1 - \epsilon)\frac{\theta(g(u))}{2} + {\rm o}\,(\theta(g(u))),
\end{eqnarray*}
as $u \to \infty$. Thus, using (\ref{1008}), the image of the sub
strip $S_\rho$ is a region of the form
$$
\left\{ z: \vert {\rm Im}\,z\vert \leq (1 - \epsilon)
\frac{\theta({\rm Re}\, z)}{2} + {\rm o}\,(\theta({\rm Re}\, z))
\right\}.
$$
Since $g(w)$ lies outside this region if $v > (1 -
\epsilon)\pi/2$, we find that
\begin{eqnarray}
{\rm Im}\,g(w)& \geq & (1 - \epsilon) \frac{\theta({\rm
Re}\,g(w))}{2} + {\rm o}\,(\theta({\rm Re}\,g(w)))  \label{v1}\\
& = & (1 - \epsilon)\frac{\theta(g(u))}{2}
+ {\rm o}\,(\theta(g(u)))\nonumber.
\end{eqnarray}
On the other hand, if $ v > (1 - \epsilon)\pi/2$ then, simply
because $g(w)$ lies in $D$,
\begin{equation}
{\rm Im}\,g(w) \leq \frac{\theta({\rm Re}\,g(w))}{2} =
\frac{\theta(g(u))}{2} + {\rm o}\,(1).
\label{v2}
\end{equation}
Together (\ref{v1}) and (\ref{v2}) imply that, for $u$
sufficiently large and $ v > (1 - \epsilon)\pi/2 $,
$$
\left\vert {\rm Im}\,g(w) - \frac{\theta(g(u))}{\pi} v \right\vert
\leq 2 \epsilon\,\theta(g(u)).
$$
It follows that (\ref{1006}) holds uniformly in $v$, $ -\pi/2 < v
< \pi/2$, with ${\rm o}\,(1)$ replaced by $C\epsilon$. Since
$\epsilon$ may be as small as we please, the proof of (\ref{1006})
is complete.
\end{proof}

We wish to remove the restriction in Warschawski's Lemma
\ref{lemmaWar} that $w$ lies in a fixed sub strip of the standard
strip $S$, at least when $D$ is a parabola-shaped domain. In this
situation,
$$
\theta(x) = 2A x^{\alpha}, \quad {\rm for }\ x >0.
$$
\begin{prop}We set $g$ to be the conformal mapping of
the strip $S$ onto ${\cal P}_{\alpha}$ for which $g$ is real and $g'$ is positive
on the real axis. Then the following estimate for the derivative of $g$ holds:
\begin{equation}
g'(w) = \left[ \frac{1}{\pi} + {\rm o}(1) \right] \theta( g({\rm
Re}\,w) ) = 2A \left[ \frac{1}{\pi} + {\rm o}(1) \right]
g(u)^{\alpha},
\label{1010}
\end{equation}
as $u$, the real part of $w$, tends to $\infty$, uniformly in
${\rm Im}\, w$.
\label{derivofg}
\end{prop}
\noindent Proposition 4 is part of the main result in
\cite{CH}, in which it is shown that the function
$$
h(z) =  -  \exp\left[ \frac{\pi}{2(1-\alpha)} \left[ 1 -
(1+Bz)^{1-\alpha}  \right]\right],\quad\ B = A^{-1/(1-\alpha)},
$$
is univalent in $\pa$ and maps $\pa$ onto the interior $D$ of a
Dini-smooth curve $C$ lying inside the unit circle $\vert z \vert
=1$, except for $z=-1$. A smooth Jordan curve $C$ is said to be
Dini-smooth if there is an increasing function $\omega(x)$, that
satisfies the Dini-condition
$$
\int_0^1 \frac{\omega(x)}{x}\,dx < \infty,
$$
for which the angle $\beta(s)$ of the tangent to $C$, considered
as a function of arclength, satisfies
$$
\vert \beta(s_2) - \beta(s_1) \vert < \omega(s_2 - s_1).
$$
The proof in \cite{CH} that $h(\pa)$ is indeed bounded by a
Dini-smooth curve, that lies inside the unit circle except for
$z=-1$, is quite involved. Proposition 4 follows from this in a
relatively straightforward manner. We repeat a version of the
argument here for the reader's convenience.

We first note that
\begin{equation}
h'(z)\, =\, - \frac{\pi }{2}\,B\,(1+Bz)^{-\alpha}\, h(z),\ \ z \in
\pa.
\label{gprime}
\end{equation}
We consider the mapping
$$
\Psi(z) = h\left[g\left(\log \frac{1+z}{1-z} \right)\right], \ \ z
\in \Delta.
$$
The unit disk $\Delta$ is first mapped onto the strip $S$, which
is mapped by $g$ onto the parabola-shaped domain $\pa$ and,
finally, this is mapped by $h$ onto the inside of a Dini-smooth
curve $C$. As explained in \cite{CH}, we are now in a position to
apply \cite[Theorem 10.2]{Pom} and may deduce that $\Psi'$ has a
continuous, non-zero extension to the closure of the unit disk. In
particular, $\Psi'(1) \not = 0$, and we note that $\Psi(1)=0$. We
can derive information on the behaviour of the derivative of $g$
from information on the derivative of $\Psi$. We write $w(z) =
\log[(1+z)/(1-z)]$, for $z \in \Delta$. Using the expression
(\ref{gprime}) for the derivative of $h$,
\begin{eqnarray*}\Psi'(z) & = & h'\bigl(g(w(z))\bigr)\,g'\bigl(w(z)\bigr)\,w'(z)\cr
& = & - \frac{\pi}{2}\, B\,
\bigl[1+B g(w(z))\bigr]^{-\alpha} \psi(z)\,
g'\bigl(w(z)\bigr)\frac{2}{1-z^2}\cr
& = &- \pi B\, \frac{1}{1+z}\,
\frac{g'\bigl(w(z)\bigr)}{\bigl[1+B g(w(z))\bigr]^{\alpha}}\,
\frac{\Psi(z)}{1-z}.
\end{eqnarray*}
We let $z \to 1$ from within the unit disk. Then
$$
\frac{\Psi(z)}{1-z}\, \to\, - \Psi'(1),
$$
which is non-zero. Hence,
$$
\frac{g'\bigl(w(z)\bigr)}{\bigl[1+B g(w(z))\bigr]^{\alpha}}\,
\to\, \frac{2}{\pi B}, \ \ {\rm as}\ z \to 1,\ z \in \Delta.
$$
It is not difficult to see that $({\rm Re}\,z)/(1+B z) \to B^{-1}$
as ${\rm Re}\,z \to \infty$ with $z \in \pa$. Substituting $z =
g(w)$ in this limit yields
$$
\frac{{\rm Re}\,g(w)}{1+B g(w)} \to \frac{1}{B}\quad\mbox{as }
{\rm Re}\,w \to \infty \mbox{ with } w \in S.
$$
Since the unrestricted limit as $z \to 1$ within the unit disk
corresponds to the unrestricted limit as  ${\rm Re}\,w \to \infty$
within the strip $S$, it follows from the previous two estimates
that
$$
\frac{g'(w)}{[{\rm Re}\,g(w)]^\alpha}\, \to\, \frac{2}{\pi
B}\,B^\alpha  = \frac{2A}{\pi}, \ \ {\rm as}\ {\rm Re}\, w \to
\infty,\ w \in S,
$$
which is (\ref{1010}).

\subsubsection{Asymptotic form of the differential operator}

Armed with the asymptotics for ${\rm Im}\, g(w)$ in Lemma
\ref{imagpartg} and those for $g'(w)$ in Proposition
\ref{derivofg}, we are now ready to derive the asymptotics for the
differential operator
\begin{equation}
\frac{\Delta k(w)}{\vert g'(w)\vert^2} -2 (n-2)\frac{{\rm Im}\,
\left[ k_w(w) \,/\,g'(w)\right]}{{\rm Im}\,\left[ g(w) \right]},
\label{diffop}
\end{equation}
which acts on $C^2$-functions defined in the strip $S$ and arises
in Lemma \ref{pdefork}. We keep in mind that, for a real-valued
function $k(w) = k(u+iv)$,
$$
k_v(w) = -2\, {\rm Im}\,\bigl(k_w(w)\bigr).
$$
First we compute, using Proposition \ref{derivofg},
\begin{eqnarray*}
- 2\,{\rm Im}\,\left[ \frac{k_w(w)}{g'(w)}\right] & = & - 2\,{\rm
Im}\,\left[\frac{k_w(w)}{\bigl[ 1/\pi + {\rm o}\,(1)\bigr]
\theta\bigl(  g(u)\bigr)}\right]\cr
& = & -2\, \frac{\pi}{\theta\bigl( g(u)\bigr)}{\rm Im}\,\left[
k_w(w)\bigl( 1 + {\rm o}\,(1)\bigr)\right]\cr
& = & -2\, \frac{\pi}{\theta\bigl( g(u)\bigr)}
\left[ -\frac{1}{2} k_v(w) + {\rm Im}\,\bigl[{\rm o}\,(1) k_w(w)\bigr] \right]\cr
& = &-2\, \frac{\pi}{\theta\bigl( g(u)\bigr)}
\left[ -\frac{1}{2} k_v(w) + {\rm o}\,(1) k_v(w)\right]\cr
& = & \frac{\pi}{\theta\bigl( g(u)\bigr)}\,k_v(w)
\bigl[ 1+{\rm o}\,(1)\bigr].
\end{eqnarray*}
Using Lemma \ref{imagpartg} to estimate the imaginary part of $g(w)$, we find that
\begin{eqnarray*}
-2\,\frac{{\rm Im}\, \left[ k_w(w) \,/\,g'(w)\right]}{{\rm
Im}\,\left[ g(w) \right]}
& = & \frac{\pi}{\theta\bigl(
g(u)\bigr)} \frac{k_v(w)[1 + {\rm o}(1)]}{[1/\pi + {\rm
o}(1)]\theta\bigl( g(u)\bigr)\,v}\cr
& = &
\frac{\pi^2}{\theta^2\bigl( g(u)\bigr)} \frac{k_v(w)}{v} \bigl[ 1
+ {\rm o}(1)\bigr].
\end{eqnarray*}
Similarly, and again using Proposition \ref{derivofg},
$$
\frac{1}{\vert g'(w)\vert^2} = \frac{\pi^2}{\theta^2\bigl(
g(u)\bigr)} \bigl[ 1 + {\rm o}(1)\bigr].
$$
In summary, the differential operator in (\ref{diffop}) becomes
\begin{equation}
\frac{\pi^2}{\theta^2\bigl( g(u)\bigr)} \left[ \bigl[ 1 + {\rm
o}(1)\bigr] \Delta k(w) + (n-2)\bigl[ 1 + {\rm
o}(1)\bigr]\frac{k_v(w)}{v}\right].
\label{xyz}
\end{equation}

\bigskip\noindent{\bf Summary to date: }The work in the foregoing
sections has been leading up to the following. Suppose that $H$ is harmonic in the
domain ${\cal P}^n_\alpha$ in $\bR^n$ and that $H(x,Y)$ is rotationally symmetric
about the $x$-axis. Suppose that the function $h$ in the planar domain $\pa^+$ is
constructed from $H$ according to $h(x+iy) = H(x,Y)$, with $\vert Y \vert = y$.
Then, by Lemma \ref{Hh}, $h$ satisfies
$$
\Delta h(z) + (n-2) \frac{h_y(z)}{y} = 0, \qquad z \in \pa^+.
$$
From $h(z)$ we construct the function $k(w)$ in the strip $S^+$
according to
$$
k(w) = h\bigl( g(w) \bigr), \qquad w \in S^+,
$$
where $g$ is a symmetric conformal mapping of the strip $S$ onto
the parabola-shaped domain $\pa$. The partial differential
equation satisfied by $k$ is given in the next proposition, which
follows directly from (\ref{xyz}) and Lemma \ref{pdefork}.
\begin{prop}There is a function $\epsilon(w)$ in the
strip $S$, with the properties that

\noindent
\item{(i)} $\epsilon(w) \to 0$ as $u \to \infty$, uniformly in $v$,

\noindent\item{(ii)} whenever the function $k(w)$ arises from a
rotationally symmetric harmonic function $H$ in $\pa^n$, as
described above, then $k$ satisfies the partial differential
equation
\begin{equation}
\Delta k(w) + [n + \epsilon(w)-2] \frac{k_v(w)}{v} = 0,\ \ w \in
S^+.
\label{pdek}
\end{equation}
\label{Section7thm}
\end{prop}
\begin{remark} If the harmonic function $H$ with
which we began had lived in a cylinder in $\bR^n$ (of radius
$\pi/2$ and with axis along the $x$-axis), then the associated
function $h$ would have the standard strip $S$ as its domain of
definition. The mapping $g$ would be the identity mapping and so
$k$ would simply satisfy $\Delta k(w) + (n-2)k_v(w)/v = 0$ in this
case. Proposition \ref{Section7thm} may be thought of as asserting
that $k$ behaves asymptotically as if it derived from a
cylindrical domain. One may also interpret Proposition
\ref{Section7thm} as asserting that while the differential
operator $\Delta h(z) + (n-2)h_y(z)/y$ is not conformally
invariant in the same way that the Laplacian is, it is
asymptotically conformally invariant. The conformal invariance of
the Laplacian was used in (\ref{3.2a}) in Section 3.2, and
Proposition \ref{Section7thm} is essentially an extension of this
to higher dimensions.
\end{remark}
\subsubsection{Sub solutions and a maximum principle}

We need to determine the boundary conditions satisfied by a
function $k$ that is constructed, as in Proposition
\ref{Section7thm}, from the rotationally symmetric harmonic function
$$
H(x,Y) = P_{(x,Y)}\{\vert B_{\tau_\alpha} \vert>t \},
$$
in the region ${\cal P}^n_\alpha$ in $\bR^n$. Thus $H$ is the
harmonic measure of the exterior of the ball of radius $t$ w.r.t.\
$\pa^n$. The gradient of $H$ w.r.t.\ $Y$ vanishes when $x=0$
because of the rotational symmetry. This translates into the
boundary condition $h_y(x,0) = 0$ for the associated function $h$
in the planar domain $\pa^+$ and, in turn, into the boundary
condition $k_v(u,0) = 0$ for the function $k(u,v)$ in the strip
$S^+$:
\begin{equation}
k_v(u) = 0, \qquad -\infty < u < \infty.  \label{bc1}
\end{equation}
The boundary values of $H$ lead to the condition on the boundary
of $\pa^+$ that $h(x,Ax^\alpha)=1$ if $\vert (x,Ax^\alpha) \vert >
t$ and $h(x,Ax^\alpha) = 0$ if $\vert (x,Ax^\alpha) \vert < t$.
Under the conformal mapping $f$ of the domain $\pa$ onto the strip
$S$, this becomes the following boundary condition for $k$:
\begin{equation}
 k(u+i\pi/2) = 0, \ \ -\infty < u < s;
\label{bc2}
\end{equation}
\begin{equation}
k(u+i\pi/2)= 1, \ \ s < u < \infty.
\label{bc3}
\end{equation}
Here $s$ depends on $t$, as specified in (\ref{3.2c}). We note
that the point $(1,0)$ in ${\cal P}^n_\alpha$ corresponds to the
point $0$ on the boundary of the strip $S^+$.

Let us therefore suppose that $k$ is a solution of the p.d.e.\
\begin{equation}
\Delta k(w) + [n + \epsilon(w)-2] \frac{k_v(w)}{v} = 0, \qquad w \in S^+,
\label{Pdek}
\end{equation}
where $\epsilon(w) \to 0$ as $u \to \infty$, uniformly in $v$,
with the boundary conditions (\ref{bc1}), (\ref{bc2}) and
(\ref{bc3}). We will show in the next section that $k(0)$ decays
at a slower exponential rate as $s \to \infty$ than solutions of
\begin{equation}
\Delta k(w) + [n + \delta - 2] \frac{k_v(w)}{v} = 0, \qquad w \in S^+,
\label{Pdedelta}
\end{equation}
when $\delta$ is positive, the boundary conditions being the same
as those satisfied by $k$.

This comparison between the solutions of (\ref{Pdek}) and
(\ref{Pdedelta}) has a heuristic interpretation that may be
helpful to keep in mind. In the limiting case $\delta = 0$, the
differential equation (\ref{Pdedelta}) becomes $\Delta k(w) +
(n-2)k_v(w)/v = 0$, the solutions of which, with the above
boundary conditions, may be thought of as deriving from harmonic
measure in a cylinder of radius $\pi/2$ in $\bR^n$. The solutions
of (\ref{Pdedelta}) may then be thought of as corresponding to
harmonic measure in such a cylinder in a slightly higher \lq
dimension\rq\ when $\delta$ is positive. Thus, our results will
show that the distribution function of the exit position of
Brownian motion from a parabola-shaped region in $\bR^n$ decays
like the distribution function of the exit position from a
cylinder in $\bR^n$, but with a time change that is given
explicitly by (\ref{3.2c}).

It is natural to consider solutions of (\ref{Pdedelta}) in the half strip
$$
S_s^+ = S^+ \cap \{u < s\} = \{ w = u+iv: u <s \mbox{ and }
0<v<\pi/2\}.
$$
In fact, by symmetry, a solution of the p.d.e.\ (\ref{Pdedelta})
in $S^+$ that satisfies the boundary conditions (\ref{bc1}),
(\ref{bc2}) and (\ref{bc3}) will take the constant value $1/2$ on
the vertical cross cut $u=s$. This boundary condition can be
satisfied by using separation of variables  to solve
(\ref{Pdedelta}) in the half strip $S_s^+$ and then taking a
series expansion. The rate of decay of the solution at 0 as $s$
becomes large is then determined by the first term in this Bessel
series. This is the term that is therefore of interest to us. For
each $m$, we write $J_m(v)$ for the Bessel function of order $m$,
we write $j_m$ for its smallest positive zero and we set
$$
\hat J_m(v) = v^{-m} J_m(v).
$$
The first term in the Bessel series for a solution of
(\ref{Pdedelta}) in $S_s^+$ is (a constant times)
\begin{equation}
\phi_\delta(w) = e^{2 j_m (u - s)/\pi} \hat J_m\left( \frac{2
j_m}{\pi}\,v \right),\ \ \mbox{where } m = \frac{1}{2}(n + \delta
-3).
\label{J}
\end{equation}
Since $\hat J_m$ satisfies the differential equation (see \cite[Section
17.22]{WW}, for example),
$$
\hat J_m''(v) + [2m+1] \frac{\hat J_m'(v)}{v} + \hat J_m(v) = 0,
$$
$\phi_\delta(w)$ satisfies the p.d.e.\ (\ref{Pdedelta}) in
$S_s^+$, as well as the boundary conditions (\ref{bc1}) and
(\ref{bc2}). On the vertical side of the half strip $S_s^+$, its
values are simply $\hat J_m\bigl( 2 j_m v/\pi \bigr)$. One needs
to take the entire series to have a solution which equals $1/2$ on
the vertical cross cut $u=s$.

We write $L$ for the operator
\begin{equation}
L[f] = \Delta f + \bigl[n + \epsilon(w) -2 \bigr] \frac{f_v}{v}.
\label{L[f]}
\end{equation}
Then $k(w)$ is a solution of $L[k]=0$ in the half strip $S_s^+$ with the boundary
conditions (\ref{bc1}) and (\ref{bc2}). In order to compare $k$ to solutions of
(\ref{Pdedelta}), we construct sub solutions for $L$ in $S_s^+$, and then use a
maximum principle. Of course, all our estimates need to be uniform in $s$. We show
how to obtain the sub solutions that we need in the next lemma.

\begin{lemma}We suppose that, for a fixed positive $\delta$, the
number $u_\delta$ is chosen so large that $2\vert \epsilon(w)
\vert < \delta $ for $u > u_\delta$. We suppose that a function
$k_\delta$ is defined in the rectangle $R_s = S_s^+ \cap
\{u>u_\delta\}$ and satisfies the partial differential equation
(\ref{Pdedelta}) there. Suppose further that $\partial
k_\delta/\partial v$ is negative in $R_s$. Then, $L[k_{\delta}]
\geq 0$ in $R_s$. In particular, $L[\phi_{\delta}] \geq 0$.
\label{subsuper}
\end{lemma}

\begin{proof}With $k_\delta$ as in the statement of the lemma,
\begin{eqnarray*}
L[k_\delta]
& = & \Delta k_\delta + \bigl[n + \epsilon(w) -2 \bigr]
    \frac{1}{v}\frac{\partial k_\delta}{\partial v}\cr
& = &\Delta k_\delta + \bigl[n + \delta  -2 \bigr]
    \frac{1}{v}\frac{\partial k_\delta}{\partial v} + \bigl[
    \epsilon(w) - \delta \bigr]
    \frac{1}{v}\frac{\partial k_\delta}{\partial v}\cr
& = & \bigl[ \epsilon(w) - \delta \bigr]
      \frac{1}{v}\frac{\partial k_\delta}{\partial v}\cr
\end{eqnarray*}
Since $2\vert \epsilon(w) \vert < \delta $ for $u > u_\delta$, it
follows that $\epsilon(w) - \delta$ has the same sign as $-\delta$
in the rectangle $R_s$. Since $\partial k_\delta/\partial v$ is
negative in $R_s$, we deduce that $L[k_\delta]$ has the same sign
as $\delta$ in $R_s$.

The statement about $\phi_\delta$ now follows from the facts that $\phi_\delta$
satisfies (\ref{Pdedelta}) and that  $\hat J_m$ is decreasing on the interval
$(0,j_m)$.

\end{proof}

The other ingredient we need is an appropriate form of the maximum
principle. While the version presented here is most probably not
new, we have been unable to find it in the literature.
Consequently, we outline the proof for completeness.

\begin{lemma}If f is a sub solution of L which is $C^2$ in the
closure of the rectangle $R_s$ (the second derivatives are
continuous up to the boundary) and which is non positive on the
top, left and right parts of the boundary, then $f(u, v)\leq 0$
for any $(u, v)\in R_s$.
\label{maxprin}
\end{lemma}

\begin{proof}Let $Z_t=(X_t, Y_t)$ be the diffusion associated with
the operator $L$.  Then $Y_t>0$ almost surely for all $t$. This is
true in the case $\epsilon(w) =0$,  since $Z_t$ is then a Bessel
process and as such it never hits zero (see \cite[Chapter XI]{RY}). If
$\epsilon$ is not zero, we still assume that $-\delta \leq
\epsilon(x, y) \leq \delta$ in the rectangle.  It follows by a
stochastic comparison theorem argument (as in the classical
Ikeda--Watanabe theorem, \cite{RW}) that $Y_t>0$ almost surely for
all $t>0$. Now, let $\tau$ be the first time that $Z_t$ hits the
boundary of the rectangle with the diffusion starting at
$z_0=(x_0, y_0)\in R_s$. This time is finite almost surely.  Of
course, by the above, $Z_{\tau}$ belongs only to the three sides
of the rectangle with probability 1.  Applying  It\^o's Lemma,
\begin{equation}
f(Z_{t \min \tau})-f(z_0)=M_t +A_t
\label{max}
\end{equation}
where $M_t$ is a martingale and $A_t=\int_0^t L[f](Z_s)ds$. Since
$f$ is a sub solution of $L$, we have $L[f]\geq  0$.  Taking
expectations of both sides of (\ref{max}), we conclude that
$$
 f(z_0)\leq E_{z_0}\left(f(Z_{t\min \tau})\right).
$$
We now let $t\to \infty$.  Since $f$ is bounded in the closure of
the rectangle,
$$
f(z_0)\leq E_{z_0}(f(Z_{\tau})),
$$
which proves $ f(z_0)\leq 0$.
\end{proof}

\subsubsection{Rate of exponential decay}

We now have all the ingredients necessary to prove the following
estimate of $k$.
\begin{prop}Suppose that the function $\epsilon(w)$, $w \in S$,
satisfies $\epsilon(w) \to 0$ as $u \to \infty$, uniformly in $v$.
Suppose that $k(w)$ is the solution of (\ref{pdek}) in $S^+$ with
the boundary conditions (\ref{bc1}), (\ref{bc2}) and (\ref{bc3}),
so that $k$ derives from harmonic measure in $\pa^n$ as in
Section~4.2.3. Then, given $\epsilon$ positive,
$$
k(0) \geq \exp\left[-\left(\frac{2 j_m}{\pi} +
\epsilon\right)s\right],
$$
for all sufficiently large $s$, where $m= (n-3)/2$.
\end{prop}
\begin{proof}
The function $k$ is bounded by $1$ on the vertical side $u=s$ of
$R_s$ (since it is but a certain harmonic measure in the
parabola-shaped region $\pa^n$ in disguise). More precisely,
$k(s,v) \to 1/2$ as $s\to \infty$, uniformly for $v$ in
$(0,\pi/2)$. In fact, the harmonic measure of $\partial \pa^n \cap
\{x>t\}$ w.r.t.\ the parabola-shaped region $\pa^n$ approaches
$1/2$ uniformly on the cross section $\{(t,Y): Y \in \bR^{n-1},\
\vert Y \vert < A t^ \alpha\}$. Hence, for all sufficiently large
$s$,
\begin{equation}
\frac{1}{4} \leq k(s,v) \leq 1,\ \ {\rm for}\ 0 < v < \pi/2.
\label{kcv}
\end{equation}

Given $\epsilon$ positive, we set $m_1 = \frac{1}{2}(n+\delta -3)$
and choose $\delta$ positive, but so small that $j_{m_1} \leq j_m
+ \pi\epsilon/4$. This is possible since the first positive zero
of the Bessel function depends continuously on the order of the
Bessel function and increases with this order. We suppose that
$u_\delta$ is as in Lemma \ref{subsuper} and that $s> u_\delta$.
Direct comparison of $k$ with the function $\phi_\delta$ of
(\ref{J}) doesn't quite work, as $\phi_\delta$ is positive on the
side $u=u_\delta$ of $R_s$ while we do not know $k$ there. We
consider the positive function
$$
k_\delta(w) = \left[e^{2 j_{m_1} (u - s)/\pi} - e^{2 j_{m_1}
[(u_\delta - s) - (u-u_\delta)]/\pi}\right] \hat J_{m_1}\left(
\frac{2 j_{m_1}}{\pi}\,v \right),\ \ w \in R_s,
$$
in which the second exponential term compensates for the positive
values of $\phi_\delta$ on the side $u=u_\delta$. This function is
a solution of (\ref{Pdedelta}) in $R_s$ and, moreover, $\partial
k_\delta/\partial v$ is negative in $R_s$ since $\hat J_{m_1}$ is
decreasing on $(0,j_{m_1})$. It follows from Lemma \ref{subsuper}
that $L[k_\delta] \geq 0$. The function $k_\delta$ satisfies zero
Dirichlet boundary conditions on the sides $u=u_\delta$ and $v =
\pi/2$ of $R_s$, and zero Neumann condition on the side $v=0$. On
the right side of the rectangle its boundary values satisfy
$$
k_\delta(s+iv) = \left[ 1 - e^{4j_{m_1}(u_\delta - s)}\right]\hat
J_{m_1}\left(2 j_{m_1} v/\pi \right) \leq \hat J_{m_1}\left(2
j_{m_1} v/\pi \right) \leq \hat J_{m_1}(0),
$$
for all sufficiently large $s$. Together with (\ref{kcv}), we see
that we can choose a fixed small, but positive, $b_1$ such that
$$
b_1 k_\delta(s+iv) \leq b_1 \hat J_{m_1}(0) \leq 1/4 \leq k(s+iv),
$$
for all sufficiently large $s$. Thus, $b_1 k_\delta - k \leq 0$ on
the three sides $u=u_\delta$, $u=s$, $v= \pi/2$ of $R_s$, while
$b_1 k_\delta - k$ satisfies a zero Neumann condition on the side
$v=0$. Since $L[b_1 k_\delta - k] = b_1 L[k_\delta] \geq 0$ in
$R_s$, we conclude from the Maximum Principle, Lemma
\ref{maxprin}, that $b_1 k_\delta - k \leq 0$ in $R_s$. This leads
to a lower bound for $k(u_\delta+1)$ since
\begin{eqnarray*}
k(u_\delta+1) & \geq & b_1 k_\delta(u_\delta+1)\cr
& = & 2 b_1 e^{2 j_{m_1} u_\delta/\pi} \sinh(2 j_{m_1}/ \pi) e^{-2 j_{m_1}
s/\pi}\hat J_{m_1}(0)\cr
& \geq & b_2 e^{-2 j_{m_1} s/\pi}.
\end{eqnarray*}
By the Harnack inequality, $k(0) \geq b_3 k(u_\delta+1)$, for a
constant that does not depend on $s$. Setting $b_4 = b_3 b_2$,
$$
k(0) \geq b_3 k(u_\delta+1) \geq b_4 \exp\left[-\frac{2 j_m}{\pi}
s - \frac{\epsilon}{2} s\right].
$$
As $b_4$ does not depend on $s$, we have $b_4 \geq e^{-\epsilon
s/2}$ for all sufficiently large $s$. The proof of Proposition~6
is complete.

\end{proof}

\subsubsection{Lower bound for harmonic measure}

The point $0$ in the strip $S$ corresponds to the point $(1,0)$ in the region
$\pa^n$ under the transformations in Subsection 4.2.1. Thus,
$$
P_{(1,0)}\left\{\vert B_{\tau_\alpha}\vert > t\right\} = k(0),
$$
where the function $k$ satisfies the partial differential equation
(\ref{pdek}) of Proposition 5 and the boundary conditions
(\ref{bc1}), (\ref{bc2}) and (\ref{bc3}) with $s = s(t)$ as given
by (\ref{3.2c}). Thus Proposition 6 leads directly to the
following lower bound for harmonic measure in the parabola-shaped
region $\pa^n$.
\begin{prop}Suppose that $\epsilon$ is positive. There
exists a constant $C_2$ depending on $\epsilon$, $n$, $A$ and
$\alpha$ such that, for $t>C_2$,
$$
P_{(1,0)}\left\{|B_{\tau_\alpha}|>t\right\} \geq
\exp\bigg[-\frac{\sqrt{\lambda_1}}{A(1-\alpha)}\,[1+\epsilon]\,
t^{1-\alpha}\bigg]
$$
\end{prop}

\subsection{Concluding remarks}The distributional inequalities in
Propositions 3 and 7 lead immediately to the limit (\ref{thm4}),
and to Theorem 3 by following the steps in the proof of Theorem 2
in Section 3.3.

\medskip It is natural to hope that the machinery constructed in
Section~4.2 would lead to an upper bound for harmonic measure, and
not only to a lower bound, thus rendering the use of the Carleman
method and Section~4.1 unnecessary. In fact, there is no mention
in Section~4.2 of bounds of any kind until Section~4.2.5. However,
we have been unable to prove a counterpart for Proposition~6
involving an upper bound for $k(0)$.

\noindent {\it Acknowledgment. We are very grateful to an anonymous referee for 
the careful reading of this paper and the many useful comments.}


\begin{thebibliography}{99}

\bibitem{ba} R.~Ba\~nuelos, R.~D.~DeBlassie, R.~Smits, The
first exit  time of  planar  Brownian motion from the interior of
a parabola, Ann.\ Prob.\ {\bf 29} (2001), 882--901.

\bibitem{van} M.~van den Berg, Subexponential behavior of the Dirichlet heat
kernel, {\sl J.\ Funct.\ Anal.\/} {\bf 198} (2003), no. 1, 28--42.

\bibitem{bu} D.~L.~Burkholder, Exit times of Brownian motion, harmonic
majorization and Hardy spaces, Adv.\ Math.\ {\bf 26} (1977),
182--205.

\bibitem{CH} T.\ Carroll and W.\ K.\ Hayman, Conformal mapping of
parabola-shaped domains, to appear in {\sl Computational Methods
and Function Theory\/}.

\bibitem{DS} D.\ DeBlassie and R.\ Smits, Brownian motion in twisted domains,
preprint 2003.

\bibitem{ess}  M. Ess\'en and K. Haliste, A problem of Burkholder and the
existence of harmonic majorants of $|x|^p$ in certain domains in
$\bR^d$, Ann.\ Acad.\ Sci.\ Fenn.\ Ser.\ A.\ I.\ Math., {\bf 9}
(1984), 107--116.

\bibitem{ha} K. Haliste, Some estimates of harmonic majorants,
Ann.\ Acad.\ Sci.\ Fenn.\ Ser.\  A.\ I.\ Math., {\bf 9} (1984),
117--124.

\bibitem{kr}S.\ G.\ Krantz, {\sl Complex Analysis:\ the Geometric Viewpoint\/},
Carus Mathematical Monographs 23, Mathematical Association of America, 1990.

\bibitem{li} W.~Li, The first exit time of Brownian motion  from an unbounded convex domain,
Annals of Probability, Vol. 31, 1078--1096,  (2003).

\bibitem{lif} M.~Lifshits and Z.~Shi, The first
exit time of Brownian motion from a parabolic domain. {\sl
Bernoulli\/} {\bf 8} (2002), no. 6, 745--765.

\bibitem{RY} D.~Revuez and M.~Yor, {\sl Continuous Martingales and
Brownian Motion\/}, Springer-Verlag, Berlin Heidelberg, 1991.

\bibitem{RW} L.~C.~G.~Rogers and D.~Williams, {\sl Diffusions, Parkov Processes, and
Martingales, Volume 2, It\^o Calculus\/}, Wiley Series in
Probability and Statistics, 1987.

\bibitem{Pom} Ch.\ Pommerenke, {\sl Univalent Functions\/},
Vandenhoeck and Ruprecht, Gottingen, 1975.

\bibitem{Sp}  F. ~Spitzer,  Some theorems concerning
two--dimensional Brownian motion, Trans.\ A.M.S.\  {\bf 87} (1958), 187--197.

\bibitem{Wa} Warschawaski, S.E., On conformal mapping of infinite strips,
Trans.\ A.M.S.\ {\bf 51} (1942) 280--335.

\bibitem{WW}E.~T.~Whittaker and G.~N.~Watson, {\sl A Course
of Modern Analysis\/}, Fourth edition, Cambridge University Press
(1952).

\end{thebibliography}
\end{document}